\documentclass[a4paper,11pt,headinclude=true]{scrartcl}

\usepackage[utf8]{inputenc}
\usepackage{ifthen}
\usepackage{amsmath,amsthm,amssymb,amsfonts}
\usepackage{graphicx}
\usepackage[percent]{overpic}
\usepackage{enumitem}
\usepackage[ruled,vlined]{algorithm2e}
\usepackage{bm}
\usepackage{accents}
\usepackage[unicode,colorlinks=true,pagebackref=false]{hyperref}
\usepackage{caption}

\newtheorem{theorem}{Theorem}
\newtheorem{remark}[theorem]{Remark}
\newtheorem{corollary}[theorem]{Corollary}
\newtheorem{proposition}[theorem]{Proposition}

\numberwithin{equation}{section}
\numberwithin{theorem}{section}

\newcommand{\Kry}{\mathcal{K}}
\newcommand{\RKry}{\mathcal{Q}}
\newcommand{\N}{\mathbb{N}}
\newcommand{\T}{\mathbb{T}}

\newcommand{\R}{\mathbb{R}}
\renewcommand{\C}{\mathbb{C}}

\newcommand{\dd}{\mathrm{d}}

\def\ee{\ensuremath{\mathrm{e}}}
\def\ii{\ensuremath{\mathrm{i}}}

\DeclareMathOperator{\real}{Re}
\DeclareMathOperator{\imag}{Im}

\DeclareMathOperator{\vspan}{span}

\DeclareMathOperator{\rank}{rank}

\newcommand\setopen{\mathrm{o}}
\newcommand\tbar[1]{\accentset{\rule{.5em}{0.8pt}}{#1\hskip 1pt}}
\newcommand\ttbar[1]{\accentset{\rule{.5em}{0.8pt}}{#1}}

\newcommand\Minr{\mbox{\scriptsize$\mathsf{M}$}}
\newcommand\tMinr{\mbox{$\mathsf{M}$}}
\newcommand\Hast{\mathsf{H}}

\newcommand{\figref}[1]{Figure~\ref{#1}}

\newif\ifthesis
\thesisfalse

\begin{document}

\title{A review of the Separation Theorem of Chebyshev-Markov-Stieltjes
for polynomial and some rational Krylov subspaces}
\author{Tobias Jawecki}
\date{
\textit{\small Institute for Theoretical Physics, Vienna University of Technology,
Vienna, Austria}\\
{\small tobias.jawecki@tuwien.ac.at}\\
{\small \today}}
\maketitle

\section*{Abstract}

The accumulated quadrature weights of Gaussian quadrature formulae
constitute bounds on the integral over the intervals between the quadrature nodes.
Classical results in this concern date back to works
of Chebyshev, Markov and Stieltjes and are referred to as
Separation Theorem of Chebyshev-Markov-Stieltjes (CMS Theorem).
Similar separation theorems hold true for some classes
of rational Gaussian quadrature.
The Krylov subspace for a given matrix and initial vector
is closely related to orthogonal polynomials
associated with the spectral distribution of the initial vector
in the eigenbasis of the given matrix,
and Gaussian quadrature for the Riemann-Stieltjes integral
associated with this spectral distribution.
Similar relations hold true for rational Krylov subspaces.
In the present work,
separation theorems are reviewed in the context of Krylov subspaces
including rational Krylov subspaces
with a single complex pole of higher multiplicity
and some extended Krylov subspaces.
For rational Gaussian quadrature
related to some classes of rational Krylov subspaces with a single pole,
the underlying separation theorems are newly introduced here.
\newline\newline\noindent\textit{Keywords: Separation Theorem of Chebyshev-Markov-Stieltjes,
intertwining property, Gaussian quadrature, rational Gaussian quadrature, Krylov subspace techniques}
\newline\noindent\textit{2020 MSC: 15A23 15A42 15B57 26A42 42C05 65D32}



\newpage
\tableofcontents


\section{Introduction and historical context}\label{sec.0mo}
In the 
\ifthesis
present chapter
\else
present work
\fi
we consider an Hermitian matrix $A\in\C^{n\times n}$ and a given vector~$u\in\C^n$.
The coefficients of $u$ in the orthonormal eigenbasis of $A$ are referred to as spectral coefficients, see~\eqref{eq.defwj} below.
These coefficients rely on an underlying inner product on $\C^n$ which is specified in~\eqref{eq.Minnerdef4},
and denoted as \tMinr-inner product in the sequel.
Furthermore, Krylov subspaces in the sequel also rely on the \tMinr-inner product.

\subsection{Historical context and previous works}\label{subsec.introhistoric}
For a polynomial or rational Krylov subspace of a matrix $A$ with starting vector $u$, the spectral coefficients of $u$ play a crucial role:
The linear functional $f \mapsto (u,f(A)\,u)_{\Minr}$ can be understood as a Riemann-Stieltjes integral
associated with a non-decreasing step function~$\alpha_n$.
This step function is defined by the eigenvalues of $A$ and the spectral coefficients of $u$,
and many results concerning the theory of polynomial Krylov subspaces have their origin in the theory of
orthogonal polynomials, namely,~polynomials on the real axis which are orthogonal
w.r.t.~the Riemann-Stieltjes integral associated with~$\alpha_n$;
see also~\cite{GM10} for a survey.
We also refer to these polynomials as orthogonal polynomials associated
with the distribution $\dd\alpha_n$.
In a similar manner, orthogonal rational functions describe rational Krylov subspaces.

For the polynomial case, the Lanczos method~\cite{La50} is used in practice
to generate an \tMinr-orthonormal basis of the Krylov subspace and the associated Jacobi matrix,
which corresponds to the representation of $A$ in the respective Krylov subspace, see also~\cite{Sa11}.
The respective \tMinr-orthonormal basis vectors satisfy a three-term recursion
which conforms to the three-term recursion of the underlying orthogonal polynomials~associated with $\dd\alpha_n$;
the Krylov basis and the orthogonal polynomials exist in an equivalent manner.

The Jacobi matrix associated with orthogonal polynomials for a given distribution
plays a crucial role for Gaussian quadrature formulae for the respective Riemann-Stieltjes integral,
which also goes by the name Gauss-Christoffel quadrature, cf.~\cite{Ga81}.
For the Gauss-Christoffel quadrature formula with $m$ quadrature nodes which integrates polynomials of degree $\leq 2m-1$ exactly,
the quadrature nodes are given by the zeros of the $(m+1)$-th orthogonal polynomial
and the quadrature weights are given by so called Christoffel numbers.
Early works on quadrature formulae~\cite{Wi62,GW69}
(historical remarks in~\cite{Ga81} also refer to earlier works of Goertzel)
show that these quadrature nodes and weights can be computed via the Jacobi matrix;
the zeros of the $(m+1)$-th orthogonal polynomial coincide with eigenvalues of the respective Jacobi matrix,
and the Christoffel numbers correspond to entries of its eigenvectors.
In these works, the underlying distribution is not necessarily based on a matrix-vector pair
as it is the case when considering a polynomial Krylov subspace;
a reference to the Krylov setting is made later in~\cite{FF94,FH93}
and also discussed in detail (including historical remarks) in~\cite{GM10,LS13}.
In this context, the eigenvalues of the Jacobi matrix are also referred to as Ritz values,
and $m$ denotes the dimension of the Krylov subspace.
Furthermore, the Christoffel numbers, which are given by entries of the eigenvectors
of the Jacobi matrix, can be written as spectral coefficients of a vector~$x\in\R^m$.
In particularly, the vector $x$~corresponds to the representation of
the starting vector~$u$ in the Krylov subspace,~i.e.,~the first unit vector scaled by the norm of~$u$.
Here, the spectral coefficients of~$x\in\R^m$ denote its
coefficients in the $\ell^2$-orthonormal\footnote{
       The notation `$\ell^2$-orthonormal'
       refers to a basis orthonormal w.r.t.~the Euclidean inner product.
}~eigenbasis of~$J_m$.

The Separation Theorem of Chebyshev-Markov-Stieltjes (CMS Theorem)
states that accumulated quadrature weights of a Gaussian quadrature formula (i.e.,~the accumulated Christoffel numbers)
are bounded by Riemann-Stieltjes integrals over the interval
between the left integral limit and the quadrature nodes.
For details and historical remarks see~\cite{Sze85,Ak65,vA93}.
In an equivalent manner, this statement can be formulated in a Krylov setting:
The accumulated entries of eigenvectors of the Jacobi matrix
(spectral coefficients of $x$)
are bounded by Riemann-Stieltjes integrals associated with $\alpha_n$ over
the interval between the left-most eigenvalue of $A$ 
and the Ritz values.
The step function $\alpha_n$ corresponds to accumulated spectral coefficients of $u$,
and as a corollary, accumulated spectral coefficients of $x$ yield bounds on sums of spectral coefficients of $u$
and vice versa.
Analogously, this statement can be formulated as an intertwining property
of the distribution $\dd\alpha_n$ and a distribution $\dd\alpha_m$
associated with the step function $\alpha_m$ which is defined by Ritz values and spectral coefficients of $x$:
Similar to $f \mapsto (u,f(A)\,u)_{\Minr}$ and $\alpha_n$, the functional $f\mapsto (x, f(J_m)\, x)_2$ can be understood
as a Riemann-Stieltjes integral associated with the step function $\alpha_m$.
The underlying Gaussian quadrature formula implies $ (u,p(A)\,u)_{\Minr} = (x, p(J_m)\, x)_2 $ for polynomials $p$ of degree $\leq 2m-1$,
and therefore, the distributions $\dd\alpha_n$ and $\dd\alpha_m$ have the same moments up to degree $2m-1$.

For distributions with the same moments, an intertwining property is stated in~\cite[Theorem 22.2]{KS53}, see also~\cite[Theorem 2.2.5]{Fi96} and~\cite[Theorem 3.3.4]{LS13}.
Indeed, this intertwining property coincides
with the result of the CMS Theorem.
In the context of Krylov subspaces, the intertwining property of
the distributions $\dd\alpha_n$ and $\dd\alpha_m$ already appeared earlier in~\cite{FF94,Fi96,LS13}.
For further remarks (including many historical remarks) on the moment problem we particularly refer to~\cite{LS13}.
The identity $ (u,p(A)\,u)_{\Minr} = (x, p(J_m)\, x)_2 $ above corresponds to quadrature properties,
and results from well-known identities for polynomials in the Krylov subspace in an equivalent manner;
$p(A)u=V_m\,p(J_m)x$ for polynomials $p$ of degree $\leq m-1$
and $p(A)u- V_m\,p(J_m)x \perp_{\Minr} V_m$ for polynomials $p$ of degree $m$
where $V_m\in\C^{m\times n}$ denotes the \tMinr-orthonormal Krylov basis written in matrix form.

The related theory in~\cite{Sze85,Ak65} 
applies in a slightly more general setting, namely, for Gaussian quadrature formulae which
integrate polynomials of degree $\leq 2m-2$ exactly (where $m$ is the number of quadrature nodes).
This includes Gauss-Radau quadrature formulae where one of the $m$ quadrature nodes
is preassigned.
The quadrature nodes and weights of Gauss-Radau quadrature formulae
can be represented by the zeros of a so-called quasi-orthogonal polynomial
and the Christoffel numbers associated with this polynomial, respectively.
Similar to the Jacobi matrix, the recursion of the underlying set of polynomials
constitutes a tridiagonal structure in matrix form,
and the respective quadrature nodes and weights correspond to the eigenvalues
and entries of eigenvectors, respectively, of this tridiagonal matrix.
This relation between Gauss-Radau quadrature formulae and the eigendecomposition
of this Jacobi-like tridiagonal matrix
goes back to~\cite{Go73} and is reviewed in detail in~\cite[Section 6.2]{GM10}.

In the present work we also consider rational Krylov subspaces, namely, 
subspaces spanned by $\{r(A)u\}$ where $r=p/q$ for polynomials $p$ of degree $\leq m-1$
and a preassigned denominator polynomial $q$ of degree $\leq m-1$
(here, $p$ and $q$ are complex polynomials and $m$ again denotes the dimension of the Krylov subspace).
For early works on rational Krylov subspaces we refer to~\cite{ER80,Ru84},
and for a review we refer to~\cite{Gu10}.
The zeros of the denominator $q$ are also referred to as poles in this context.
Rational Krylov techniques using a single pole of multiplicity $m-1$ yield the most prominent cases,
the resulting rational Krylov subspaces are also referred to as \textit{Shift-and-Invert (SaI)} Krylov subspaces.

The rational Krylov subspace with preassigned denominator polynomial~$q$ and starting vector~$u$
is identical to the polynomial Krylov subspace with starting vector~$q(A)^{-1}u$.
The respective orthogonal polynomials
(particularly, orthogonal polynomials associated with a scaled distribution~$\dd\widehat{\alpha}_n$)
divided by the denominator polynomial~$q$ yield rational functions which are orthogonal w.r.t.~the
Riemann-Stieltjes integral associated with~$\alpha_n$ (as given previously), cf.~\cite{DB07}.
These orthogonal rational functions, evaluated at~$A$ as a matrix function and applied to~$u$,
provide an \tMinr-orthonormal basis of the rational Krylov subspace.
Furthermore, results regarding Gaussian quadrature formulae carry over to the rational setting:
The orthogonal rational functions which span the rational Krylov subspace
of dimension $m$ with a preassigned denominator $q$ constitute a
rational quadrature formula for the Riemann-Stieltjes integral
associated with $\alpha_n$,
which integrates rational functions~$r=p/|q|^2$ exactly for
polynomials~$p$ of degree~$\leq 2m-1$.
For an overview on rational Gaussian quadrature see also~\cite{Ga93},
and for the relation between rational Krylov subspaces
and rational Gaussian quadrature we also refer to~\cite{LRW08, De09, JR11}.

The relation between a rational Krylov subspace with denominator~$q$
and starting vector~$u$,
and the polynomial Krylov subspace with starting vector~$q(A)^{-1}u$
is more of a theoretical nature.
In practice, various algorithms, covering different settings,
are relevant to construct a rational Krylov subspace,
and result in different sequences of \tMinr-orthonormal basis vectors of this subspace.
To keep our results general, we do no restrict ourselves to a specific algorithm
or an underlying recursion for the basis vectors
in that concern.
Assuming an \tMinr-orthonormal basis of a rational Krylov subspace is given,
we refer to 
the representation of~$A$ in this basis
as \textit{Rayleigh quotient}~$A_m\in\C^{m \times m}$.
Furthermore, we reuse the notation~$x\in\C^m$ for the representation of~$u$ in the given basis.
As stated above, a rational Krylov subspace is closely related
to orthogonal rational functions which constitute a rational Gaussian quadrature formula.
In particular, the quadrature nodes and weights for this rational Gaussian quadrature formula
correspond to the eigenvalues of the Rayleigh quotient $A_m$
and the spectral coefficients of $x$, respectively.
We remark that the eigenvalues of $A_m$ (also referred to as rational Ritz values, which are real due to $A_m$ being Hermitian)
and the spectral coefficients of $x$,
which refer to the coefficients of $x$ in the $\ell^2$-orthonormal
eigenbasis of $A_m$, are independent of the choice of the basis.
Furthermore, the respective quadrature formula
conforms to the identity $(u,r(A)\,u)_{\Minr} = (x, r(A_m)\, x)_2$
for rational functions $r=p/|q|^2$ as above.

Similar to the polynomial case, the functional $ f \mapsto (x, f(A_m)\, x)_2 $
can be understood as a Riemann-Stieltjes integral associated with $\alpha_m$,
which is now defined by eigenvalues of $A_m$ and the spectral coefficients of $x$.
The rational quadrature properties imply that $\dd\alpha_n$ and $\dd\alpha_m$ have $2m-1$ identical rational moments.

For rational Gaussian quadrature formulae, CMS type results
depend on the choice of the denominator, and do not seem to be as popular as for the polynomial case.
In~\cite{Li98} a separation theorem is
given for a class of Laurent polynomials and an integral defined on the positive real axis.
Here, Laurent polynomials correspond to rational functions with denominator~$q(\lambda)=\lambda^{m/2}$ for even~$m$. 
In a Krylov setting, this class of rational functions is related to some extended Krylov subspaces~\cite{DK98}
for a matrix~$A$ with positive eigenvalues (i.e.,~the step function~$\alpha_n$ is defined on the positive real axis).
However, the results of~\cite{Li98} have not been applied in a Krylov setting yet.

More recently,~\cite{ZTK19} computes piecewise estimates on~$\alpha_n$ based
on a Shift-and-Invert Krylov subspace with a pole of multiplicity~$m-1$ at zero
(i.e.,~$q(\lambda)=\lambda^{m-1}$), for a matrix~$A$ with positive eigenvalues.
In this work, a Shift-and-Invert representation is used instead of the Rayleigh quotient (see also~\cite[Subsection 5.4.3]{Gu10}).
The given estimates are based on an intertwining property of~$\dd\alpha_n$ and a distribution given by spectral properties
of the Shift-and-Invert representation; the intertwining property goes back to the polynomial case, referring to~\cite[Theorem 22.2]{KS53}.

In the present work, we also consider Krylov techniques related to rational Gauss-Radau quadrature formulae.
These quadrature formulae integrate rational function $r=p/|q|^2$ exactly,
where $p$ is a polynomial of degree $\leq 2m-2$, $q$ is the given denominator,
and one of the $m$ quadrature nodes is preassigned, see also~\cite{LRW08,JR13}.
For rational Gauss-Radau quadrature formulae in a more general setting see also~\cite{Ga04,DB10,DB12a}.
Analogously to the Gauss-Radau quadrature formulae in the polynomial case,
this slightly generalizes the previously discussed rational quadrature properties
but can be treated similarly concerning the intertwining properties of the underlying distributions $\dd\alpha_n$ and $\dd\alpha_m$.

\subsection{Applications}


\paragraph{Computable estimates on spectral coefficients of $u$.}

A direct computation of eigenvalues and spectral coefficients requires access to the
eigenbasis of the given matrix $A$
which is not practical for problems of a large problem size $n$ in general;
typically, the full spectrum of $A$ is not available.
However, information on partly accumulated spectral coefficients,
namely,~the step function $\alpha_n$ on subsets of the spectrum of $A$,
can be sufficient for some applications.
CMS type results provide suitable estimates for this purpose,
which can be evaluated using Krylov techniques.
In particularly, this yields piecewise estimates on $\alpha_n$ covering the full spectrum of $A$.
These estimates hold true
independently of the convergence of individual (rational) Ritz values.
However, more detailed information is provided for parts of the spectrum
which are well resolved by (rational) Ritz values.
We proceed to give some applications based on estimates on $\alpha_n$.

The eigenvalues of $A$ together with the spectral coefficients of $u$
have some relevance for the approximation of the action of a matrix function $f(A)u$,
e.g.,~the matrix exponential function or the matrix inverse.
Polynomial Krylov methods yield good approximations on matrix functions 
without any a~priori information on the spectrum
\ifthesis
of~$A$
(such approximants are discussed in Chapter~\ref{ch.chapterBIT} and~\ref{ch.NA} for the matrix exponential function).
\else
of~$A$.
\fi
However, further knowledge on the spectrum of $A$ can help to improve
the quality of the approximation (here we also refer to the introduction of~\cite{FH93}).
In~\cite{FF94}, piecewise estimates on $\alpha_n$ are applied
to construct a polynomial preconditioner for the conjugate gradient method.
This approach is based on the intertwining property of the distributions $\dd\alpha_n$ and $\dd\alpha_m$,
where the latter is computed using a small number of Lanczos iterations in the progress
(thus, $\alpha_m$ is based on a polynomial Krylov subspace here).

In~\cite{HPS09}, the authors consider iterative bidiagonalization methods
to solve ill-posed linear systems. 
In this work, effects of a noisy
right-hand side on the projected problem are discussed.
The ill-posed problems therein are associated with an underlying distribution
(similar to $\dd\alpha_n$ given previously in the present introduction),
and due to problem assumptions and noise on the initial data this distribution
is of a special structure which carries over to the projected problem.
This process is closely related to the
intertwining property of the distributions $\dd\alpha_n$ and $\dd\alpha_m$ in the Lanczos case,
and results in criteria to detect the noise level on the run, as introduced in~\cite{HPS09}.

In~\cite{ZTK19}, an inhomogeneous differential equation,
arising in applications of dynamic analysis of structure,
 is diagonalized using eigenvectors of a large matrix.
This requires computation of a moderate number of eigenvectors, namely,
eigenvectors such that the external force vector is resolved with sufficient accuracy.
The spectral decomposition of this vector is associated with a distribution $\dd\alpha_n$,
and estimates on this distribution allow to determine intervals which cover eigenvalues corresponding to the required eigenvectors.
In~\cite{ZTK19}, estimates on $\alpha_n$ are based on a Shift-and-Invert Lanczos method,
and yield a pole selection strategy and stopping criteria for an eigenproblem solver based on rational Krylov methods.

In future works, estimates on $\alpha_n$ will be applied to
design special rational approximations to the action of the exponential of skew-Hermitian matrices. 

\paragraph{The structure of $\alpha_n$ roughly carries over to $\alpha_m$.}
\ifthesis
In Chapter~\ref{ch.localbestapprox}, we consider
\else
In~\cite{Ja22b}, the authors consider
\fi
a localized best approximation property of rational Krylov approximants
to the action of a matrix exponential.
In particularly, we consider the exponential of a skew-Hermitian matrix applied to a vector
which is subject to some assumptions.
Namely, strict increases of $\alpha_n$ are, up to a small perturbation,
located in an interval.
For some rational Krylov subspaces we illustrate that such properties carry over
to the associated step function $\alpha_m$.
These ideas are based on theoretical results derived in the
\ifthesis
present chapter;
\else
present work;
\fi
and
\ifthesis
in Chapter~\ref{ch.localbestapprox}
\else
in~\cite{Ja22b}
\fi
this approach motivates a localized best approximation
result which can show a mesh-independent convergence
(in a setting where the matrix exponential 
arises from a spatial discretization of a PDE (evolution equation)).
In contrast to previously mentioned applications, computable estimates on $\alpha_n$
are not topical
\ifthesis
for the approach of Chapter~\ref{ch.localbestapprox}.
\else
for~\cite{Ja22b}.
\fi

\paragraph{Other applications.}
Apart from the Krylov setting, the CMS Theorem
has applications in various fields, e.g., for a work
on discretization of quantum systems see~\cite{Re79}.

Furthermore, bounds on distribution functions have some importance
in probability theory and statistics;
and various bounds are referenced to Chebyshev, Markov, Stieltjes and others.
This includes variants of the CMS Theorem formulated in terms of moments,
e.g.,~\cite{Ze54} or more recently~\cite{Hue15}.
Moment-matching methods also appear in the context of system theory~\cite{An05}.

\!\!\!Krylov methods also have applications in the
approximation to bilinear forms~$(u,f(A)\,u)_{\Minr}$, where $f$ is a given function,
see also~\cite{LRW08, GM10, JR11, JR13}.
Due to the relation between $(u,f(A)\,u)_{\Minr}$ and a Riemann-Stieltjes integral associated with $\alpha_n$,
estimates on this bilinear form are directly related to quadrature formulae.
However, these applications will not be further discussed in the present work.


\subsection{Main contributions and overview of present work}

We proceed to highlight the main contributions of the
\ifthesis
present chapter,
\else
present work,
\fi
including results or remarks which are considered to be new by the author.

\begin{itemize}

\item We introduce a new CMS type result for a class of rational Gaussian quadrature formulae,
namely, quadrature formulae based on rational functions with a single real pole of higher multiplicity,
see Theorem~\ref{thm.saixxx} in Subsection~\ref{subsec.SepTheoremRatKrylov}.
To prove this result, we introduce rational majorants and minorants on Heaviside type functions
in Proposition~\ref{prop.thatr2}.
In a Krylov setting, this theorem applies to the SaI Krylov subspace with a real shift.
Our results include the case that the shift is located in the contour of the matrix spectrum;
we consider a more general setting compared to~\cite{ZTK19}.
An intertwining property of the distributions $\dd\alpha_n$ and $\dd\alpha_m$ 
holds true up to a constant, see Proposition~\ref{prop.SaIRip}.

\item
For the setting of rational functions with a single complex pole of higher multiplicity,
we introduce a new CMS type result which yields an upper bound on the Riemann-Stieltjes integral
over the interval between neighboring quadrature nodes and at the boundary,
see Proposition~\ref{prop.SaIC1} in Subsection~\ref{subsec.CMSSaIcomplex}.
This result applies to the SaI Krylov subspace with a single {\em complex} shift of higher multiplicity.
To prove this upper bound,
we make use of polynomial majorants on Heaviside type functions
on the unit circle given in~\cite{Go02}.
Furthermore, we propose the use of an isometric Arnoldi method
to compute the Rayleigh quotient of the SaI Krylov subspace with complex shift
in a cost efficient way (comparable to the Lanczos method which applies when the shift is real), see Remark~\ref{rmk.SaIRayleighcomplexpole}.

\item  Applying a CMS type result given in~\cite{Li98},
we present an intertwining property for $\dd\alpha_n$ and $\dd\alpha_m$ in the setting of an extended Krylov subspace in Subsection~\ref{subsec.STCMSxKry}.

\end{itemize}

Recalling results of~\cite{GM10} and others, we also apply the theory
of quasi-orthogonal polynomials in a polynomial Krylov setting.
This results in an Arnoldi-like decomposition where the residual is provided
by a quasi-orthogonal polynomial; we refer to the respective representation
as a \textit{quasi-orthogonal residual (qor-)} Krylov representation
for which one of the eigenvalues can be preassigned.

\begin{itemize}
\item The CMS Theorem is known to apply to Gauss-Radau quadrature formulae.
In the present work, we specify these results in a Krylov setting;
results in Section~\ref{sec.bounds3} for the polynomial case include the qor-Krylov setting,
e.g.,~the intertwining property of $\dd\alpha_n$ and $\dd\alpha_m$ holds true
when $\alpha_m$ is based on the qor-Krylov representation.
This potentially leads to refined estimates on $\alpha_n$ in practice.

\item Furthermore, we introduce a qor-Krylov approximation to the action of matrix functions
in Subsection~\ref{sec.polynomialqorKry},
comparable to the corrected Krylov scheme for the matrix exponential function given in~\cite{Sa92}.
\end{itemize}

Various results for the polynomial case carry over to the rational case, and
we introduce a rational qor-Krylov representation
where one of the eigenvalues is preassigned, similar to~\cite{LRW08,JR13}.
\begin{itemize}
\item For the rational case, we introduce an efficient procedure to compute
a rational qor-Krylov representation in Subsection~\ref{sec.ratqor}.
\item The CMS type result given in Subsection~\ref{subsec.SepTheoremRatKrylov}
and further estimates in Subsection~\ref{subsec.CMSSaIcomplex}
include the rational qor-Krylov case.
Considering these CMS theorems,
for some cases bounds on quadrature weights related to
quadrature nodes at the right boundary (of the spectrum of $A$)
are affected by $\alpha_n$ at the left boundary
(of the spectrum of $A$) and vice versa,~e.g., as in~Corollary~\ref{cor.piecewiseboundsrat};
$\alpha_n$ affects the bounds in a cycled sense at the boundaries.
This is no longer the case when one of the nodes is preassigned at the boundary of the spectrum,
see also Remark~\ref{rmk.CMSqorsplit},
and this potentially results in refined bounds.

\item We introduce a rational qor-Krylov approximation to the action of matrix functions
in Subsection~\ref{sec.ratqor}.

\end{itemize}

\paragraph{Overview of present work}

In Section~\ref{sec.intro1}
we first recall some theory of orthogonal polynomials 
and the relation between orthogonal polynomials and the polynomial Krylov subspace.
Here, polynomials are orthogonal w.r.t.~an inner product on the vector space,
which can be written as a Riemann-Stieltjes integral associated with
a non-decreasing step function $\alpha_n$.
Furthermore, we recall some known results for rational Krylov subspaces
based on the polynomial case.
In Subsection~\ref{subsec.remarksSaI}
we provide some remarks on the SaI Krylov subspace.
This includes a new approach to compute the SaI Krylov subspace with a complex shift
based on the isometric Arnoldi method -- a short-term recursion.
In Section~\ref{sec.2qorpol} we recall some theory on quasi-orthogonal polynomials
which results in a polynomial and rational qor-Krylov representation
in Subsection~\ref{sec.polynomialqorKry} and~\ref{sec.ratqor}, respectively.
Here, we also include some algorithmic details.

The main results of the
\ifthesis
present chapter
\else
present work
\fi
concerning
CMS theorems and intertwining properties
of distributions are stated in Section~\ref{sec.bounds3}.
We first recall quadrature properties in Subsection~\ref{subsec.bounds31}
concerning polynomial and Gaussian quadrature formulae for the
Riemann-Stieltjes integral associated with the step function $\alpha_n$.
Quadrature nodes and weights for these quadrature formulae
are provided by the Jacobi matrix or the Rayleigh quotient
of the respective Krylov subspace.
The following results in Subsection~\ref{subsec.bounds3polSepTh}--\ref{subsec.STCMSxKry}
are stated for quadrature nodes and weights of respective quadrature formulae,
and as such apply to eigenvalues and spectral coefficients for
representations in the respective Krylov subspaces.
In Subsection~\ref{subsec.bounds3polSepTh} we recall the classical CMS Theorem
which applies to the polynomial Krylov setting.
Besides other remarks in this subsection,
we also specify the step function $\alpha_m$ and recall the intertwining
property of the distributions $\dd\alpha_n$ and $\dd\alpha_m$.
In Subsection~\ref{subsec.SepTheoremRatKrylov} we introduce new results
concerning rational Gaussian quadrature formulae for a class of rational
functions with a single pole $s\in\R$ of higher multiplicity.
This result applies to the SaI Krylov setting with a real shift, and
the distributions $\dd\alpha_n$ and $\dd\alpha_m$
(whereat, $\dd\alpha_m$ is now provided by the rational Krylov subspace)
satisfy an intertwining property up to a constant shift.
In Subsection~\ref{subsec.CMSSaIcomplex} we proceed with a similar upper bound
for the rational case with a single pole $s\in\C$ of higher multiplicity,
which corresponds to a SaI Krylov setting with a complex shift.
In Subsection~\ref{subsec.STCMSxKry} we apply CMS theorems
given in~\cite{Li98} in the setting of an extended Krylov subspace,
which yields results similar to the polynomial case.
Previously discussed intertwining properties which correspond to
CMS theorems
are verified by numerical examples in Section~\ref{sec.experiments4}.


\section{Krylov subspace techniques and orthogonal polynomials}\label{sec.intro1}

A basis of a Krylov subspace obtained by the Lanczos method is closely related
to the theory of orthogonal polynomials. This relationship is explained in~\cite{GM10} and others
and is reviewed here.

In the sequel, let $A\in\C^{n\times n}$ be a given Hermitian matrix,
and let $u\in\C^n$ be a given initial vector.
The polynomial Krylov subspace, with~$m\leq n$, is denoted by
\begin{equation}\label{eq.KryAu}
\Kry_m(A,u) = \vspan\{u,A\,u,\ldots,A^{m-1}u \} ~\subset \C^n.
\end{equation}
Krylov subspace techniques rely on an inner product.
Although the Euclidean inner product on the underlying vector space
is practical in many cases, we consider a more general notation:
For two vectors $x,y\in\C^n$ we define the \tMinr-inner product
by\footnote{
        The \tMinr-inner product given in~\eqref{eq.Minnerdef4} induces
        a vector norm,~i.e.,~$\|x\|_{\Minr} = \sqrt{(x,x)_{\Minr}}$,
        which is equivalent to the Euclidean norm.}
\begin{equation}\label{eq.Minnerdef4}
(x,y)_{\Minr} = x^{\Hast} M y,
\end{equation}
where $M \in\C^{n\times n}$ is
an Hermitian\footnote{
         The matrix $M$ is Hermitian
         w.r.t.~the Euclidean inner product,~i.e., $M=M^{\Hast}$.
}~positive definite matrix which is given by the underlying problem setting.
This notation includes the Euclidean inner product, namely,
the case $M=I$ with\footnote{
        By $(\cdot,\cdot)_2$ and $\|\cdot\|_2$
        we denote the Euclidean inner product and norm, respectively.
}~$(x,y)_{\Minr} = (x,y)_2 $.
In the current work, the motivation behind the \tMinr-inner product
lies in problems which are based on
discretized Hilbert spaces,~e.g.,~for a FEM discretization
of the Hilbert space $L^2$ (on a spatial domain)
the inner product $(x,y)_{\Minr} = x^{\Hast} M y$
with $M$ representing the mass matrix of the finite element space
is a natural choice.

In the sequel we assume that $A$ is Hermitian (self-adjoint) w.r.t.~the \tMinr-inner product,
$$
(A\,x,y)_{\Minr}=(x,A\,y)_{\Minr},~~~x,y\in\C^n.
$$
Let~$\lambda_1,\ldots,\lambda_n \in\R$ denote the eigenvalues
and $q_1,\ldots,q_n\in\C^n$ the \tMinr-orthonormal eigenvectors
of~$A\in\C^{n\times n}$,~i.e.,~$Aq_j=\lambda_j q_j$
with $(q_j,q_k)_{\Minr}=\delta_{jk}$,
and let
\begin{equation}\label{eq.defwj}
w_j=(q_j,u)_{\Minr}\in\C
\end{equation}
denote the corresponding spectral coefficients
of the initial vector $u\in\C^n$, i.e.,
$$
u=\sum_{j=1}^n\,w_j\,q_j.
$$

In practice, the Lanczos method (cf.~\cite{Sa03}) delivers an
\tMinr-orthonormal\footnote{
    For two vectors $x,y\in\C^m$ an \tMinr-orthonormal basis $V_m$
    satisfies $(V_m\,x, V_m\,y)_{\Minr}=(x,y)_2$.}
basis~$V_m=(v_1,\ldots,v_m)\in\C^{n\times m}$
of the Krylov subspace $\Kry_m(A,u)$,~i.e.,
$$
\vspan(V_m)=\Kry_m(A,u),~~~\text{and}~~~(V_m,V_m)_{\Minr}=I,
$$
for which the starting vector $u$ satisfies
$$
(V_m,u)_{\Minr}=\beta_0\,e_1,~~~\beta_0 = \| u \|_{\Minr}~~~\text{and}~~e_1=(1,0,\ldots,0)^\Hast\in\R^m.
$$
`Full rank' of $ \Kry_m(A,u) $ means
\begin{equation}\label{eq.polAfullrank}
\rank\big( u,A\,u,\ldots,A^{m-1}u \big) = m.
\end{equation}
To proceed with the construction of~$\Kry_{m+1}(A,u)$
at the~$m$-th Lanczos iteration step we require~\eqref{eq.polAfullrank} to hold also for~$m+1$.
Otherwise $\Kry_m(A,u)$ is an invariant subspace of~$A$, and we refer
to this case as a {\em lucky breakdown} after~$m$ steps.
We remark that only if there exist at least $m$
coefficients~$w_j\neq 0$
with distinct eigenvalues~$\lambda_j$, then~\eqref{eq.polAfullrank}
holds true for a respective\footnote{
   See Proposition~\ref{rmk.wlamrank}, Appendix~\ref{appendix.A2}.
}~$m$.
In the sequel we will assume that no lucky breakdown occurs:
I.e., without loss of generality we assume that~\eqref{eq.polAfullrank} holds true for~$m\leq n$,
hence we consider~$w_j\neq 0$ with distinct eigenvalues~$\lambda_j$ for~$j=1,\ldots,n$.
We further assume the ordering
$$
\lambda_1 < \lambda_2 < \ldots < \lambda_n.
$$

With~\eqref{eq.KryAu} there exist polynomials $p_0,\ldots,p_{m-1}$
which satisfy
$$
v_\ell=p_{\ell-1}(A)u,~~~\ell=1,\ldots,m.
$$
For these polynomials the orthonormal property of
$V_m$,~i.e.,~$(v_\ell,v_k)_{\Minr}=\delta_{\ell k}$,
yields
\begin{equation}\label{eq.innereuclideanfromVm}
(p_{\ell-1}(A)u,p_{k-1}(A)u)_{\Minr}=\delta_{\ell k},~~~~\ell,k=1,\ldots,m.
\end{equation}

Various properties of Krylov subspaces have their origin
in the theory of orthogonal polynomials for which
we mainly refer to~\cite{Sze85,Ak65}.
The theory therein
can be formulated in terms of an integral-based inner product:
Following~\cite{GM10},
depending on~$u$
we consider the step function 
\begin{subequations}\label{eq.alpnanRSint}
\begin{equation}\label{eq.measureAv}
\alpha_n(\lambda)=\left\{
\begin{array}{ll}
0,&\lambda<\lambda_1, \\
\sum_{j=1}^\ell |w_j|^2, & \lambda_\ell \leq \lambda < \lambda_{\ell+1},~~\ell=1,\ldots,n-1,\\
\sum_{j=1}^n |w_j|^2, & \lambda_n \leq \lambda.
\end{array}
\right.
\end{equation}
We choose an interval $(a,b)$ which includes $\lambda_1,\ldots,\lambda_n$.
For $f\colon \R\to\C$ we have
\begin{equation}\label{eq.innerrsintintro}
\sum_{j=1}^n |w_j|^2 f(\lambda_j)
= \int_a^b f(\lambda)\, \dd\alpha_n(\lambda),
\end{equation}
where the right-hand side is to be understood as a Riemann-Stieltjes integral.
For the corresponding inner product
we introduce the notation 
\begin{equation}\label{eq.defalphaninnerprod}
(f,g)_{\alpha_n} = \int_a^b \tbar{f}(\lambda)g(\lambda)\, \dd\alpha_n(\lambda).
\end{equation}
\end{subequations}
In the eigenbasis of $A$ the vector $p(A)u$,
where $p$ is a polynomial,
has the representation
$$
p(A)u=\sum_{j=1}^n \,p(\lambda_j) \,w_j\,q_j.
$$
For two complex polynomials $p$ and $g$
the \tMinr-inner product of $p(A)u$ and $g(A)u$ reads
\begin{equation}\label{eq.innereuclideanpq}
{(p(A)u,g(A)u)}_{\Minr} = \sum_{j=1}^n |w_j|^2\,\tbar{p}(\lambda_j)g(\lambda_j).
\end{equation}
With~\eqref{eq.innerrsintintro},~\eqref{eq.defalphaninnerprod}
and~\eqref{eq.innereuclideanpq} we have the equivalent formulations
\begin{equation}\label{eq.innerrseuclidean}
(p,g)_{\alpha_n} = \int_a^b \tbar{p}(\lambda)g(\lambda)\, \dd\alpha_n(\lambda)
= \sum_{j=1}^n |w_j|^2\,\tbar{p}(\lambda_j)g(\lambda_j)
= {(p(A)u,g(A)u)}_{\Minr}.
\end{equation}
Thus, polynomials which satisfy~\eqref{eq.innereuclideanfromVm}
are indeed `$ \alpha_n $\,-\,orthonormal', i.e.,
\begin{equation}\label{eq.orthonormalpol}
(p_{\ell},p_{k})_{\alpha_n}  = \delta_{\ell k},~~~\ell, k=0,\ldots,m-1.
\end{equation}

We remark that the normalization factor $\beta_0$ as given previously
satisfies the identities
\begin{equation}\label{eq.beta0recall}
\beta_0 =  \big((u,u)_{\Minr}\big)^{1/2} = \big((1,1)_{\alpha_n}\big)^{1/2}
= \Big( \int_a^b 1 \,\dd \alpha_n(\lambda) \Big)^{1/2}.
\end{equation}


%
\paragraph{Three-term recursion, zeros of orthogonal polynomials, and the Jacobi matrix.}
Our assumption that a lucky breakdown does not occur for any $m<n$
corresponds to~$w_j\neq 0$ and~$\lambda_j$ being distinct for~$j=1,\ldots,n$
and entails
that the step function~$\alpha_n$ has~$n$ points of strict increase.
Following~\cite[Section 2.2]{Sze85} the respective inner product yields orthonormal polynomials
$p_0,\ldots,p_{n-1}$ of degree $0,\ldots,n-1$, respectively.
These polynomials enjoy a three-term recursion,
see also~\cite[Section 2.2]{GM10} or~\cite{Ak65,Sze85}:
\begin{proposition}\label{prop.threetermorthopol}
Let~$\beta_0=(\int_a^b 1 \,\dd \alpha_n)^{1/2}$ as in~\eqref{eq.beta0recall}.
With~$p_0 = 1/\beta_0$,~$p_{-1}= 0$ and~$m<n$
there exist~$a_1,\ldots,a_m\in\R$,~$\beta_1,\ldots,\beta_m>0$
and $ \alpha_n $-orthonormal polynomials~$p_0,\ldots,p_m$
for which the three-term recursion
\begin{equation}\label{eq.threetermrecursionp}
\lambda\,p_{j-1}(\lambda) = \beta_{j-1} p_{j-2}(\lambda) + a_{j} p_{j-1}(\lambda) + \beta_{j} p_{j}(\lambda),~~~j=1,\ldots,m,
\end{equation}
holds.
Here, $a_j = (p_{j-1},\lambda p_{j-1})_{\alpha_n}$,
and $\beta_j>0$ is fixed such that $(p_j,p_j)_{\alpha_n}=1$.
\end{proposition}

In the sequel the notation $p_0,\ldots,p_m$ refers to the orthonormal polynomials
from Proposition~\ref{prop.threetermorthopol},
where $p_j$ is of degree $j$ for $j=0,\ldots,m$ due to the recursion~\eqref{eq.threetermrecursionp}.


\begin{proposition}[See Section 3.3 in~\cite{Sze85}]\label{prop.ritzvaluespol}
We recall the following well-known properties of the zeros of~$p_m$;
\begin{enumerate}
\item The zeros~$\theta_1,\ldots,\theta_{m}\in\R$ of~$p_m$ are distinct.
Assume
$$
\theta_{1}<\theta_{2}<\ldots<\theta_{m}.
$$
\item The zeros of~$p_m$ and the
eigenvalues~$\lambda_1,\ldots,\lambda_n$ of $A$ are interlacing.
This means~$\lambda_1 < \theta_1$,~$ \theta_m < \lambda_n~$,
and for~$k=1,\ldots,m-1$ there exists at least one~$\lambda_{j(k)}$ with
$$
\theta_{k} < \lambda_{j(k)} < \theta_{k+1}.
$$
\end{enumerate}
\end{proposition}
The three-term recursion~\eqref{eq.threetermrecursionp}
can be represented in terms of the so-called symmetric {\em Jacobi matrix}~$J_m$,
whose eigenvalues coincide with the zeros of~$p_m$:
With~$a_1,\ldots,a_m\in\R$ and~$\beta_1,\ldots,\beta_m>0$,
\begin{equation}\label{eq.defJm}
J_m =
\begin{pmatrix}
a_1&\beta_1&&&\\
\beta_1&a_2&\beta_2&&\\
  &\ddots&\ddots&\ddots&\\
  &&\beta_{m-2}&a_{m-1}&\beta_{m-1}\\
  &&&\beta_{m-1}&a_{m}
\end{pmatrix}\in\R^{m\times m}.
\end{equation}
Denoting $P(\lambda)=(p_0(\lambda),\ldots,p_{m-1}(\lambda))^\Hast\in\C^m$,
the recursion~\eqref{eq.threetermrecursionp} can be written in matrix form,
\begin{equation}\label{eq.threetermrecursionpmatrix}
\lambda P(\lambda) = J_m\,P(\lambda) + \beta_m\,p_{m}(\lambda) e_m.
\end{equation}
From~\eqref{eq.threetermrecursionpmatrix} we observe that the zeros~$\theta_1,\ldots,\theta_m$ of~$p_m$
are eigenvalues of~$J_m$ with non-normalized eigenvectors~$P(\theta_j)=(p_0(\theta_j),\ldots,p_{m-1}(\theta_j))^\Hast$,
$$
\theta_j P(\theta_{j}) = J_m\,P(\theta_j),~~~j=1,\ldots,m.
$$
We conclude that the matrix~$J_m$ has
$m$ distinct eigenvalues~$\theta_1,\ldots,\theta_m\in\R$
which are indeed identical to the zeros of~$p_m$
and for which the properties from
Proposition~\ref{prop.ritzvaluespol} hold true.
We refer to~$\theta_1,\ldots,\theta_m\in\R$ as {\em Ritz values.}

\paragraph{Polynomial Krylov subspace.}
We recall the usual denotation~$V_m=(v_1,\ldots,v_m)\in\C^{n\times m}$
for the \tMinr-orthonormal basis of $\Kry_m(A,u)$
provided by the Lanczos method. We have $\vspan\{V_m\}=\Kry_m(A,u)$
and $(V_{m+1}, V_{m+1})_{\Minr} = I$
where $V_{m+1}$ includes the subsequent basis vector $v_{m+1}$.
The basis~$\{v_1,\ldots,v_{m+1}\}$ satisfies a three-term recursion
according to the Lanczos algorithm~\cite[Section 6.6]{Sa03}.
(The Lanczos algorithm in~\cite[Section 6.6]{Sa03} relies on the
Euclidean inner product but can be generalized in a direct manner.)
Substituting $A$ for $\lambda$ in~\eqref{eq.threetermrecursionp} and applying $u$
yields a recursion for $p_0(A)u,\ldots,p_{m}(A)u$ which coincides
with the Lanczos three-term recursion.
Hence,~$v_j=p_{j-1}(A)u$ for~$j=1,\ldots,m+1$
with the orthonormal polynomials $p_0,\ldots,p_m$ from Proposition~\ref{prop.threetermorthopol}.
Analogously to~\eqref{eq.threetermrecursionpmatrix}
the three-term recursion defining $V_m$ can be written in matrix form,
\begin{equation}\label{eq.c3.arnoldidecom}
A\,V_m = V_m\,J_m + \beta_m v_{m+1} e_m^\Hast.
\end{equation}
We refer to~$\beta_m\,v_{m+1}$ as a {\em residual}. 
With~\eqref{eq.c3.arnoldidecom} and the \tMinr-orthogonality property of~$V_m$ the Jacobi matrix satisfies
$$
J_m=(V_m,A\,V_m)_{\Minr}.
$$

The tridiagonal structure of $J_m$ 
implies $A^j\,u = \beta_0 V_m\,J_m^j\,e_1$ for $j=0,\ldots,m-1$
and $\beta_0 V_m e_1 = u$ where~$\beta_0=\|u\|_{\Minr}$~\cite{DK89,Sa92}.
Thus,\footnote{
         The denotation $\Pi_j$ refers to the
          class of complex polynomials of degree $\leq j$.
}
\begin{subequations}\label{eq.pmequalities}
\begin{equation}\label{eq.pofAequalKm}
p(A)u =  \beta_0 V_m\,p(J_m) e_1,~~~p\in \Pi_{m-1}.
\end{equation}
Furthermore, the corresponding deviation for a polynomial~$p\in \Pi_{m}$ of exact degree~$m$ is in the span of the residual,
\begin{equation}\label{eq.pmperpKm}
\beta_0 V_m\,p(J_m) e_1 - p(A)u\in\vspan\{v_{m+1}\}\perp_{\Minr} \Kry_m(A,u).
\end{equation}
\end{subequations}
\begin{proposition}\label{prop.inprodpol}
With respect to the \tMinr-inner product the identity
\begin{equation}\label{eq.inprodpol}
(u,p(A)u)_{\Minr} = \beta_0^2\,(e_1,p(J_m)e_1)_2,~~~~p\in \Pi_{2m-1}
\end{equation}
holds true.
\end{proposition}
\begin{proof}
For $p\in\Pi_{2m-1}$ we can write $p=g_1\,g_2$
with~$g_1\in\Pi_{m-1}$ and~$g_2\in\Pi_m$, and
\begin{equation}\label{eq.inprodpolp1}
(u,p(A)u)_{\Minr} = (\tbar{g}_1(A) u,g_2(A) u)_{\Minr}
~~~\text{and}~~(e_1,p(J_m)e_1)_2 = (\tbar{g}_1(J_m) e_1,g_2(J_m)e_1)_2.
\end{equation}
For $\tbar{g}_1(A) u$ and $g_2(A) u$
we apply~\eqref{eq.pofAequalKm} and~\eqref{eq.pmperpKm}, respectively, to conclude
\begin{equation}\label{eq.inprodpolp2}
(\tbar{g}_1(A) u,g_2(A) u)_{\Minr} =
\beta_0^2 (V_m\,\tbar{g}_1(J_m) e_1, V_m\,g_2(J_m) e_1)_{\Minr}.
\end{equation}
With $(V_m,V_m)_{\Minr}=I$ we recall
\begin{equation}\label{eq.inprodpolp2x}
(V_m\,\tbar{g}_1(J_m) e_1, V_m\,g_2(J_m) e_1)_{\Minr}
= (\tbar{g}_1(J_m) e_1, g_2(J_m) e_1)_2.
\end{equation}
Combining~\eqref{eq.inprodpolp1},~\eqref{eq.inprodpolp2}
and~\eqref{eq.inprodpolp2x}
implies~\eqref{eq.inprodpol}.
\end{proof}



\paragraph{Rational Krylov subspace.}
For rational Krylov subspaces we consider rational functions~$r=p/q$ with
a preassigned denominator~$q$.
The zeros of~$q$ are also referred to as the poles of~$r$.
Using the notation~$s_1,s_2,\ldots\in \C\cup \pm\infty$ for the poles of~$r$, for which we define
\begin{equation}\label{eq.defqm}
q_{m-1}(\lambda) = \prod_{j=1,\,s_j\neq\pm\infty}^{m-1} (\lambda-s_j).
\end{equation}
Here we admit~$s_j=\pm \infty$ in order to include cases
for which the denominator of~$r$ is of a smaller degree than its numerator.
This can be used to constitute the so called extended Krylov subspace, see also~\cite{DK98}
and will further be relevant in Subsection~\ref{sec.ratqor} below.\footnote{
      For~$0\neq s_1,s_2,\ldots \in\C~$ we can exchange the factors of~$q_{m-1}$,
      i.e.,~$(\lambda-s_j)$,
      with~$(1-\lambda/s_j)$ to obtain a definition of~$q_{m-1}$ which is equivalent to~\eqref{eq.defqm}.
      This clarifies the convention~$s_j=\pm \infty$ in~\eqref{eq.defqm}.
}

We assume that the poles $s_j$ are distinct from the eigenvalues~$\lambda_{1},\ldots,\lambda_{n}$ of~$A$,
such that $q_{m-1}^{-1}(A)$ is well-defined.
The rational Krylov subspace~$\RKry_m(A,u)$ with poles $s_1,\ldots,s_{m-1}$
and $q_{m-1}$ from~\eqref{eq.defqm},
is defined by the 
span of~$\{r(A)u\colon r=p/q_{m-1}~~\text{for}~ p\in \Pi_{m-1} \}$,~i.e.,\footnote{
     In the sequel, we also use $q_{m-1}^{-1}(\lambda)$ for $q_{m-1}(\lambda)^{-1}=1/q_{m-1}(\lambda)$ to shorten the denotation.}
\begin{subequations} \label{eq.RKry_uq}
\begin{equation}\label{eq.defQm}
\begin{aligned}
\RKry_m(A,u)&:= \vspan\{q_{m-1}^{-1}(A) u, A\,q_{m-1}^{-1}(A) u,\ldots, A^{m-1}q_{m-1}^{-1}(A) u\}\\
&\;=\Kry_m(A, q_{m-1}^{-1}(A)u).
\end{aligned}
\end{equation}
To simplify the notation we write
\begin{equation}
u_q=q_{m-1}^{-1}(A)u.
\end{equation}
\end{subequations}
With~\eqref{eq.RKry_uq}, the rational Krylov subspace
$ \RKry_m(A,u) $ is identical to the
polynomial Krylov subspace~$\Kry_m(A, u_q)$.
Let~$w_j$ be the spectral coefficient of $u$ w.r.t.\ the eigenvalue $\lambda_j$,
then $q^{-1}_{m-1}(\lambda_j)w_j$ is the corresponding spectral coefficient of~$u_q$.
Analogously to $\alpha_n$ in~\eqref{eq.measureAv},
we introduce the step function 
\begin{equation}\label{eq.measureAvrat}
\widehat{\alpha}_n(\lambda)=\left\{
\begin{array}{ll}
0,&\lambda<\lambda_1\\
\sum_{j=1}^\ell |q^{-1}_{m-1}(\lambda_j)w_j|^2, & \lambda_\ell \leq \lambda < \lambda_{\ell+1},~~~\ell=1,\ldots,n-1,\\
\sum_{j=1}^n |q^{-1}_{m-1}(\lambda_j)w_j|^2, & \lambda_n \leq \lambda.
\end{array}
\right.
\end{equation}
Analogously to~\eqref{eq.defalphaninnerprod}, the Riemann-Stieltjes
integral associated with $\widehat{\alpha}_n$
defines an inner product,
$$
(f,g)_{\widehat{\alpha}_n} = \int_a^b \tbar{f}(\lambda)g(\lambda)\, \dd\widehat{\alpha}_n(\lambda).
$$
The $ \widehat{\alpha}_n $-orthonormal polynomials given
by Proposition~\ref{prop.threetermorthopol}
constitute a basis of~$\Kry_m(A, u_q)$.
For the existence of these orthonormal polynomials,
analogously as before we assume that
we have $n$~coefficients $w_j\neq 0$ for distinct eigenvalues $\lambda_j$,
together with $0\neq q_{m-1}^{-1}(\lambda_j)\in\C$.

Let~$J_m$ and~$V_m$ be the Jacobi matrix and the \tMinr-orthonormal basis for~$\Kry_m(A, u_q)$.
For the eigenvalues~$\theta_1,\ldots,\theta_m$ of~$J_m$ the
results of Proposition~\ref{prop.ritzvaluespol} remain valid.

The Jacobi matrix $J_m=(V_m, A\,V_m)_{\Minr}$
corresponds to a representation of $A$ in the 
underlying rational Krylov subspace $\RKry_m(A,u)=\Kry_m(A,u_q)$.
However, $V_m$ and $J_m$
are more of a theoretical nature in this context.
In practice, $u_q=q_{m-1}^{-1}(A)u$, is not directly available
and the rational Krylov subspace is not constructed
via its polynomial counterpart, but in an iterative manner.
While the Lanczos method is by far the most prominent approach
to construct a polynomial Krylov subspace,
various iterative algorithms are relevant for the
rational case.
Choosing a proper algorithm to construct a rational Krylov subspace
depends on the setting, e.g.,~the choice of poles.
For computational details we also refer to~\cite{DB07,BG15,Gu10,BR09}.
Unlike the polynomial case,
where $V_m$ refers to the orthonormal basis constructed by the Lanczos method,
we choose the notation for the rational Krylov subspace
independent of the underlying algorithm:
We assume $U_m\in\C^{n\times m}$
is a given \tMinr-orthonormal basis of~$\RKry_m(A,u)$, i.e.,
\begin{equation*}
U_m\in\C^{n\times m},~~~
\vspan\{U_m\} = \RKry_m(A,u)~~~\text{and}~~(U_m,U_m)_{\Minr}=I,
\end{equation*}
and we let $A_m$ refer to the respective {\em Rayleigh quotient}
\begin{equation*}
A_m = (U_m,A\,U_m)_{\Minr}\in\C^{m\times m}.
\end{equation*}
For instance, this notation covers rational Krylov bases
and representations constructed as in Subsection~\ref{subsec.remarksSaI}.
The matrix~$A_m$ is Hermitian w.r.t.~the Euclidean inner product
but in general not tridiagonal and does not coincide with~$J_m$.
\begin{subequations} \label{eq.Km}
Let us denote
\begin{equation}\label{eq.defKm}
K_m = (V_m,U_m)_{\Minr} \in\C^{m\times m}.
\end{equation}
$U_m$ and~$V_m$ represent orthonormal bases of the same subspace, thus,
\begin{equation} \label{eq.UmtoVmKm}
U_m = V_m(V_m,U_m)_{\Minr} = V_m K_m.
\end{equation}
By definition of the \tMinr-inner product  we have $K_m^\Hast K_m = U_m^\Hast M V_m K_m$,
and together with $V_m K_m=U_m$~\eqref{eq.UmtoVmKm} this yields
\begin{equation}\label{eq.Kmortho1}
K_m^\Hast\,K_m
= (U_m,U_m)_{\Minr} = I.
\end{equation}
\end{subequations}
Furthermore, $A_m$ and~$J_m$ are orthogonally similar matrices,
\begin{equation}\label{eq.RKAmtoJm}
A_m = (U_m,A\,U_m)_{\Minr} = K_m^\Hast\,(V_m,A\,V_m)_{\Minr}\,K_m
= K_m^\Hast\,J_m\,K_m,
\end{equation}
therefore, the eigenvalues of~$A_m$ are equal
to the Ritz values~$\theta_1,\ldots,\theta_m$ corresponding to~$\Kry_m(A, u_q)$.

We proceed with some identities in the rational Krylov subspace,
a rational counterpart to~\eqref{eq.pmequalities}.
Assume that $q_{m-1}^{-1}(A_m)$ is well-defined, and
let
$$ x := (U_m, u)_{\Minr}\in\C^m. $$
Then,
\begin{subequations}
\begin{equation}\label{eq.pmqeqRKm}
r(A) u = U_m\,r(A_m) x~~~\;\text{for $r=p/q_{m-1}$ with $p\in\Pi_{m-1}$.}
\end{equation}
This result was given earlier in~\cite[Lemma 4.6]{Gu10} and others.
Furthermore, let~$r=p/q_{m-1}$ for a polynomial~$p\in \Pi_{m}$
of degree exactly~$m$, then\footnote{
     A proof of~\eqref{eq.pmqeqRKm} and~\eqref{eq.pmqperpRKm}
     is also provided in~Proposition~\ref{prop.equalitiesratKry}, Appendix~\ref{appendix.A2}.
}
\begin{equation}\label{eq.pmqperpRKm}
( U_m\,r(A_m)x - r(A)u) \perp_{\Minr} \vspan\{ U_m\} = \RKry_m(A,u).
\end{equation}
\end{subequations}

Following~\cite[Remark 3.2]{Gu13b} we conclude:
\begin{proposition}\label{prop.ratkryeqinnerprod}
For $x = (U_m, u)_{\Minr}\in\C^m$ and rational functions~$r=p/|q_{m-1}|^2$ with~$p\in \Pi_{2m-1}$,
\begin{equation}\label{eq.ratkryeqinnerprod}
(u,r(A)u)_{\Minr} = (x,r(A_m) x)_2.
\end{equation}
\end{proposition}
\begin{proof}
For~$r\in \Pi_{2m-1}/|q_{m-1}|^2$ we write $r=r_1\,r_2$
with $r_1\in\Pi_{m-1}/\tbar{q}_{m-1}$ and $r_2\in\Pi_{m}/q_{m-1}$, and
\begin{equation}\label{eq.ratkryeqinnerprodp1}
(u,r(A)u)_{\Minr} = (\tbar{r}_1(A)u,r_2(A)u)_{\Minr},
~~~\text{and}~~(x,r(A_m)x)_2
= (\tbar{r}_1(A_m) x,r_2(A_m) x)_2.
\end{equation}
For $\tbar{r}_1\in\Pi_{m-1}/q_{m-1}$ and $r_2\in\Pi_{m}/q_{m-1}$ we apply~\eqref{eq.pmqeqRKm} and~\eqref{eq.pmqperpRKm}, respectively,
to conclude
\begin{equation}\label{eq.ratkryeqinnerprodp2}
(\tbar{r}_1(A)u,r_2(A)u)_{\Minr}
= (U_m\,\tbar{r}_1(A_m) x, U_m\,r_2(A_m) x)_{\Minr}.
\end{equation}
Combining~\eqref{eq.ratkryeqinnerprodp1} with~\eqref{eq.ratkryeqinnerprodp2}
and making use of $(U_m,U_m)_{\Minr}=I$ implies~\eqref{eq.ratkryeqinnerprod}.
\end{proof}

\subsection{Some remarks on the Shift-and-Invert (SaI) Krylov subspace}\label{subsec.remarksSaI}
The poles~$s_j$ are not required to be distinct.
A prominent example is the
\begin{center}
{\em Shift-and-Invert \,(\textbf{\em SaI})\, Krylov subspace,}
\end{center}
with~$q_{m-1}(\lambda)=(\lambda-s)^{m-1}$
for a single pole~$s\in\C$ 
of multiplicity~$m-1$.

\begin{remark}
\label{rmk.SaIRayleigh}
The rational Krylov subspace~$\RKry_m(A,u)$ with a single pole $s\in\C$
of multiplicity $m-1$ is identical
to the polynomial Krylov subspace~$\Kry_m(X,u)$
with~$X=(A-s\,I)^{-1}$,~i.e.,
$$
\Kry_m(X,u)=\vspan\{ u, (A-sI)^{-1}u, \ldots, (A-sI)^{-(m-1)}u\}.
$$
Note that~$\RKry_m(A,u)\subset\Kry_m(X,u)$
via the partial fraction decomposition for rational functions
with denominator~$q_{m-1}(\lambda)=(\lambda-s)^{m-1}$,
and~$\Kry_m(X,u)\subset\RKry_m(A,u)$
by normalizing.
Thus, the rational Krylov subspace~$\RKry_m(A,u)$
can be constructed analogously as the
polynomial Krylov subspace~$\Kry_m(X,u)$.
The matrix~$X$ is no longer Hermitian for~$ \imag s\neq 0$, and in this case the construction of
the Krylov subspace~$\Kry_m(X,u)$ requires the Arnoldi method,
the counterpart of the Lanczos method for general matrices.
Further computational details for the case~$ \imag s\neq 0$ are
given in Remark~\ref{rmk.SaIRayleighcomplexpole} below.
The Lanczos or Arnoldi method for~$\Kry_m(X,u)$ generates
an orthonormal basis~$U_m$
and an upper Hessenberg matrix~$X_m = (U_m,X U_m)_{\Minr}$.
With the subsequent basis vector $u_{m+1}$ and
$x_{m+1,m}=(X_{m+1})_{m+1,m}$, the Arnoldi decomposition of~$\Kry_m(X,u)$ (similar to~\eqref{eq.c3.arnoldidecom}) gives
\begin{equation}\label{eq.ArnoldiDecomposition}
(A-s\,I)^{-1} U_m = U_m\,X_m + x_{m+1,m}\,u_{m+1}\,e_m^\Hast.
\end{equation}
With~\eqref{eq.ArnoldiDecomposition} and using the notation~$y_m^\Hast=e_m^\Hast\,X_m^{-1}$, we obtain
\begin{equation}\label{eq.SAIKrylovid}
A\,U_m = U_m(X_m^{-1} + s\,I) - x_{m+1,m}  (A-s\,I) u_{m+1}y_m^\Hast.
\end{equation}
For the Rayleigh quotient~$A_m = (U_m,A\,U_m)_{\Minr} $,
identity~\eqref{eq.SAIKrylovid} implies
\begin{equation}\label{eq.Am0}
A_m  = X_m^{-1} + s\,I - x_{m+1,m}\,(U_m,A\,u_{m+1})_{\Minr}\, y_m^\Hast.
\end{equation}
This identity can be further simplified
in view of numerical efficiency and stability
(similar to~\cite[eq.~(5.7)]{DK98} for~$s=0$ or~\cite[eq.~(5.8)]{Gr12}
for~$s\in\R$):
With $A$ being Hermitian and the identity~\eqref{eq.SAIKrylovid} we have
\begin{equation}\label{eq.Am0inner}
(U_m,A\,u_{m+1})_{\Minr} = (A\,U_m,u_{m+1})_{\Minr}
= -x_{m+1,m}\big(  (A\,u_{m+1},u_{m+1})_{\Minr} - \tbar{s}\big)y_m.
\end{equation}
Combining~\eqref{eq.Am0} and~\eqref{eq.Am0inner}
together with $\kappa=(u_{m+1},A\,u_{m+1})_{\Minr}\in\R$
yields
\begin{equation}\label{eq.Am1}
A_m  = X_m^{-1} + s\,I + x_{m+1,m}^2 (\kappa-\tbar{s}) y_m\,y_m^\Hast.
\end{equation}
With~$A_m$ and~$y_m\,y_m^\Hast\in\C^{m\times m}$ being Hermitian we take the Hermitian part of~\eqref{eq.Am1} to obtain
\begin{equation*}
A_m = (X_m^{-1}+(X_m^{-1})^\Hast)/2 + \real(s)\,I + x_{m+1,m}^2 (\kappa-\real(s)) y_m\,y_m^\Hast.
\end{equation*}
This representation for $ A_m $ 
is equivalent to~\eqref{eq.Am0}
but it is better suited for numerical computation.
A shift of the inverse of the Hessenberg matrix~$X_m$, i.e.,~$X_m^{-1}+s\,I$, is closely related to the Rayleigh quotient~$A_m$, see also~\cite[Subsection 5.4.3]{Gu10},
but it does not conserve orthogonality. 
E.g., for~$s \notin \R $ the matrix~$X_m^{-1}+s\,I$ is not necessarily Hermitian.

Note that $x = \beta_0\,e_1$ for the SaI Krylov subspace.
\end{remark}

\begin{figure}
\centering
\begin{algorithm}[H]
\caption{An algorithm to compute an orthonormal basis~$U_m$
and the Rayleigh quotient~$A_m$ of~$\RKry_m(A,u)$
for a single pole $s\in\C$ of multiplicity~$m-1$,
the SaI case.}
\label{alg.SaI}
\SetAlgoLined
 $X=(A-sI)^{-1}$\;
 if $s\in\R$~~apply the Lanczos method for~$\Kry_{m-1}(X,u)$\;
 ~~else ~~apply the Arnoldi method for~$\Kry_{m-1}(X,u)$\;
 in both cases this returns~$\beta_0$,~$U_{m},X_{m}=(U_m,XU_m)_{\Minr},\beta_{m},u_{m+1}$\;
 $\kappa=(u_{m+1},A\,u_{m+1})_{\Minr}\in\R$\;
 $y_m^\Hast=e_m^\Hast\,X_m^{-1}$\;
 $A_m = (X_m^{-1}+(X_m^{-1})^\Hast)/2 + \real(s)\,I + \beta_m^2 (\kappa-\real(s)) y_m\,y_m^\Hast$\;
 set~$x=\beta_0 e_1$\;
 return~$x, U_m, A_m$\;
\end{algorithm}
\end{figure}

The procedure which is stated in Remark~\ref{rmk.SaIRayleigh} is summarized in Algorithm~\ref{alg.SaI}.

In some works concerning the SaI Krylov subspace, the matrix $X_m^{-1}+sI$ appears in place of the
Rayleigh quotient, e.g.~\cite{EH06,ZTK19}; for a comparison see also~\cite[Subsection 5.4.3]{Gu10}.
In the following remark we show that for $s\in\R$ the matrix $X_m^{-1}+sI$
satisfies an identity similar to~\eqref{eq.ratkryeqinnerprod}.
\begin{remark}\label{rmk.quadSaIKry}
Let $X=(A-sI)^{-1}$ for a given shift $s\in\R$.
Thus, $X$ is Hermitian.
Then, the matrix $X_m=(V_m,X\,V_m)_{\Minr}$ associated
with the polynomial Krylov subspace $\Kry_m(X,u)$ satisfies
$ (u,p(X)u)_{\Minr} = (x,p(X_m)x)_2 $
for $p\in\Pi_{2m-1}$ due to Proposition~\ref{prop.inprodpol}.
Polynomials of $X$ can be rewritten as rational functions of $A$,
see also Remark~\ref{rmk.gcircxisrat} in Appendix~\ref{appendix.A4}.
A polynomial in $X_m$ can be rewritten in an analogous manner:
We recall $q_{m-1}(\lambda) = (\lambda-s)^{m-1}$
for the given shift $s\in\R$.
For a given $p\in\Pi_{2m-2}$
we have $r\in\Pi_{2m-2}/q_{m-1}^2$ with $p(X)=r(A)$ and
$p(X_m) = r(X_m^{-1}+sI) $.
Thus, similar to~\eqref{eq.ratkryeqinnerprod} we have the identity
\begin{equation}\label{eq.XmSaIeqinnerprod}
(u,r(A)u)_{\Minr} = (x,r(X_m^{-1}+sI)x)_2,~~~~r\in\Pi_{2m-2}/|q_{m-1}|^2.
\end{equation}
Here, we remark $ |q_{m-1}| = q_{m-1}$ for $s\in\R$.
In~\eqref{eq.XmSaIeqinnerprod},
the numerator is of degree $2m-2$ instead of $2m-1$ as 
in~\eqref{eq.ratkryeqinnerprod}.
\end{remark}

\noindent We proceed with some additional remarks on the SaI Krylov subspace
with a complex shift~$s\in\C\setminus\R$.
\begin{remark}
\label{rmk.SaIRayleighcomplexpole}
As stated in Remark~\ref{rmk.SaIRayleigh},
the rational Krylov subspace with a single pole $s\in\C$
of multiplicity $m-1$
corresponds to the polynomial Krylov subspace $\Kry_m(X,u)$ with
$X=(A-s\,I)^{-1}\in\C^{n\times n}$.
Let us consider the case~$s\in\C\setminus\R$.

In contrast to the case $s\in\R$, the matrix $X$ is not Hermitian
for $s\in\C\setminus\R$,
and thus, the Lanczos three-term recursion fails to construct
the Krylov subspace $\Kry_m(X,u)$.
The Arnoldi method can be applied in this case but results in additional
computational cost compared to the Lanczos method.
However, to preserve some favorable properties of the Lanczos method in the case of $s\in\C\setminus\R$,
we can construct the Krylov subspace by applying an isometric Arnoldi method
on a transformed matrix, using a Cayley transform:
We recall that $A\in\C^{n\times n}$ is Hermitian w.r.t.~the \tMinr-inner product.
Then, the matrix
$$
Z=(A-\tbar{s}\,I)(A-s\,I)^{-1}\in\C^{n\times n}
$$
is unitary w.r.t.~the \tMinr-inner product,~%
i.e.,~$(Z\,v,Z\,w)_{\Minr}=(v,w)_{\Minr}$ for $v,w\in\C^n$.
We introduce the notation $\tau$ for the corresponding scalar Cayley transform
\begin{equation}\label{eq.cayleyscalarSaI}
\tau(\lambda)=(\lambda-\tbar{s})(\lambda-s)^{-1},
~~~~ \tau\colon \R \to \T\setminus\{1\},
\end{equation}
where $\T\subset\C$ denotes the unit circle.
The matrix $Z$ has eigenvalues $\tau(\lambda_j)$ and eigenvectors $q_j$,
where $\lambda_j$ and $q_j$ denote the eigenvalues and
eigenvectors of $A$, respectively.
The function $\tau$ as given in~\eqref{eq.cayleyscalarSaI} is bijective,
which implies that $A$ and $Z$ have the same number of distinct eigenvalues
with nonzero spectral coefficients $w_j=(q_j,u)_{\Minr}$.
From remarks stated previously in the current section,
and Proposition~\ref{rmk.wlamrank} in Appendix~\ref{appendix.A2},
we conclude that the rank of $\RKry_m(A,u)$
and the rank of $\Kry_m(Z,u)$ are identical.
For the polynomial Krylov subspace $\Kry_m(Z,u)$ we observe
$\Kry_m(Z,u)\subset\RKry_m(A,u)$ by normalizing.
Due to having the same rank, the rational Krylov subspace~$\RKry_m(A,u)$
and the polynomial Krylov subspace~$\Kry_m(Z,u)$ are identical.

For $\Kry_m(Z,u)$ we consider the following setting:
Let $V_m$ denote an \tMinr-orthonormal basis of the Krylov subspace $\Kry_m(Z,u)$
with $(V_m,u)_{\Minr}=\beta_0 e_1$
and an upper Hessenberg matrix $Z_m=(V_m,Z\,V_m)_{\Minr}\in\C^{m\times m}$,
and let $v_{m+1}$ denote the
subsequent basis vector with normalization factor $z_{m+1,m}=(Z_{m+1})_{m+1,m}>0$,
$\|v_{m+1}\|_{\Minr}=1$ and $(V_m,v_{m+1})_{\Minr}=0$,
such that
\begin{equation}\label{eq.KryDecompositionCayley}
Z\, V_m = V_m Z_m + z_{m+1,m}\,e_m^\Hast v_{m+1}.
\end{equation}
Such a representation can be generated by a short term Arnoldi method,
e.g.,~the isometric Arnoldi method~\cite[Algorithm~3.1, eq.~(3.4) and (3.5)]{JR94}
introduced in~\cite{Gra93,JR94}.
For further details we also refer to~\cite{BF97,Schae08,BMV18}.
We also recapitulate the isometric Arnoldi method in Algorithm~\ref{alg.isometricArnoldi}.
In contrast to the standard Arnoldi method,
the isometric Arnoldi method is more efficient in terms of computational cost,
comparable to the Lanczos algorithm for Hermitian matrices.

Let the decomposition~\eqref{eq.KryDecompositionCayley} be given
and set $U_m:=V_m$,
then $U_m$ conforms to an orthonormal basis of the rational
Krylov subspace $\RKry_m(A,u)$ with
denominator~$q_{m-1}(\lambda)=(\lambda-s)^{m-1}$,
and $x=\beta_0 e_1$.
Substituting $Z$ and $V_m$ in~\eqref{eq.KryDecompositionCayley} yields
$$
(A-\tbar{s}\,I)(A-s\,I)^{-1}\, U_m = U_m Z_m + z_{m+1,m}\,e_m^\Hast v_{m+1}.
$$
Similar to Remark~\ref{rmk.SaIRayleigh}, this provides a computable formulation
for the Rayleigh quotient $A_m=(U_m, AU_m)_{\Minr}$. With
$\kappa=(v_{m+1},A\,v_{m+1})_{\Minr}\in\R$ and $y_m^{\Hast}:=e_m^\Hast (I-Z_m)^{-1}$
we have
$$
A_m = (\tbar{s}I - s Z_m)(I - Z_m)^{-1}  + z_{m+1,m}^2 \big( \kappa-\tbar{s}\big) y_m y_m^{\Hast}.
$$
This procedure is summarized in Algorithm~\ref{alg.isometricSaI}.

As an alternative approach to compute the SaI Krylov subspace with $s\in\C\setminus\R$,
we also remark that the matrix $X = (A-s\,I)^{-1}$ is in the class of so called
$\mathrm{normal}(1,1)$~matrices
(cf.~\cite{BM00}),~i.e.,~
the \tMinr-adjoint of $X$ corresponds to
a rational function $p(X)q(X)^{-1}$ with $p,q\in\Pi_1$,~namely,
$$
X^\ast=(A-\tbar{s}\,I)^{-1} = (X^{-1} + (s-\tbar{s})\,I)^{-1}
= X\,(I + (s-\tbar{s})\,X)^{-1},
$$
due to $ X^{-1} = A-s\,I $.
For $\mathrm{normal}(1,1)$~matrices a short Arnoldi recurrence exists,
see~\cite{BM00,BMV18},
but we do not further discuss this approach in the current work.

\end{remark}

\begin{figure}
\centering
\begin{algorithm}[H]
 \caption{
An isometric Arnoldi method to compute an orthonormal basis~$U_m$
and the~$Z_m=(U_m,Z\,U_m)_{\Minr}$ of
$\Kry_m(Z,u)$ for a unitary matrix $Z$,~e.g.,
a Cayley transform $ Z = (A-\overline{s}\,I)(A-s\,I)^{-1}$ where $A$ is
an Hermitian matrix and $s\in\C\setminus\R$.
See Remark~\ref{rmk.SaIRayleighcomplexpole} and references therein.}
\label{alg.isometricArnoldi}
\SetAlgoLined
 $\beta_0 = \|u\|_{\Minr}, ~~ v_1 =u/\beta_0, ~~ \widehat{v}=v_1, ~~ Z_m = I_{m\times m}$\;
 for $k=1:m$\;
 ~~ $w = Z v_k $\;
 ~~ $\gamma = -(\widehat{v},w)_{\Minr}$\;
 ~~ $v_{\text{next}} = w + \gamma \widehat{v}$\;
 ~~ $\sigma = \|v_{\text{next}}\|_{\Minr}$; ~~~ // $= (1 - | \gamma|^2 )^{1/2}$ in exact arithmetic\; 
 ~~ $v_{k+1} = v_{\text{next}}/\sigma$\;
 ~~ if $k<m$\;
 ~~ ~~ $(Z_m)_{:,[k:k+1]} \leftarrow (Z_m)_{:,[k:k+1]}\cdot
          \left(\begin{array}{cc} -\gamma & \sigma \\ \sigma & \ttbar{\gamma} \end{array}\right)$\;
 ~~ ~~ $\widehat{v} \leftarrow \sigma\widehat{v} + \ttbar{\gamma} v_{k+1}$\;
 ~~ ~~ $\widehat{v} \leftarrow \widehat{v}/\|\widehat{v}\|_{\Minr}$;  ~~~ // not required in exact arithmetic\; 
 ~~ else // $k=m$\;
 ~~ ~~ $(Z_m)_{:,k} \leftarrow -\gamma (Z_m)_{:,k} $\;
 ~~ ~~ $ z_{m+1,m}=\sigma$\;
 return $\beta_0,U_m=(v_1,\ldots,v_m),v_{m+1},Z_m,z_{m+1,m}$\;
\end{algorithm}
\end{figure}

\begin{figure}
\centering
\begin{algorithm}[H]
 \caption{
An optimized algorithm to compute an orthonormal basis~$U_m$ and the Rayleigh quotient~$A_m$ of~$\RKry_m(A,u)$
for a single pole $s\in\C\setminus\R$ of multiplicity~$m-1$,
see Remark~\ref{rmk.SaIRayleighcomplexpole}.}
\label{alg.isometricSaI}
\SetAlgoLined
 $Z=(A-\tbar{s}\,I)(A-s\,I)^{-1}$\;
 apply the isometric Arnoldi method for $\Kry_m(Z,u)$, see Algorithm~\ref{alg.isometricArnoldi}\;
 this returns~$\beta_0$,~$U_{m},Z_{m}=(U_m,ZU_m)_{\Minr},z_{m+1,m},u_{m+1}$\;
 $\kappa=(u_{m+1},A\,u_{m+1})_{\Minr}\in\R$\;
 $y_m^\Hast=e_m^\Hast\,(I-Z_m)^{-1}$\;
 $A_m = (\tbar{s}I-s Z_m)(I-Z_m)^{-1}  + z_{m+1,m}^2 \big( \kappa-\tbar{s}\big) y_m y_m^{\Hast}$\;
 set~$x=\beta_0 e_1$\;
 return~$x, U_m, A_m$\;
\end{algorithm}
\end{figure}


\section{A review on quasi-orthogonal polynomials}\label{sec.2qorpol}

The theory of quasi-orthogonal polynomials is for instance covered in~\cite{Sze85,Ak65,GM10}.
We will refer to a 
special linear combination $ \widehat{p}_m $ of $p_{m-1}$ and $p_m$
as a {\em quasi-orthogonal polynomial} of degree $m$,
where $p_0,\ldots,p_m$ denote the ortho\-normal polynomials from the previous section.
In the class of quasi-orthogonal polynomials we impose an additional condition,
i.e., we require that
\begin{equation} \label{eq.pdach(xi)=0}
\text{the quasi-orthogonal polynomial~$\widehat{p}_m$ vanishes
at a given $\xi\in\R$, \,i.e.,\, $\widehat{p}_m(\xi)=0$.}
\end{equation}
Quasi-orthogonal polynomials also appear in the theory of Gauss-Radau quadrature formulae.
Similar to the three-term recursion of the orthogonal polynomials,
the underlying recursion of the polynomials $p_1,\ldots,p_{m-1}$ and $\widehat{p}_m$
constitutes a matrix $T_m$ which coincides with the Jacobi matrix $J_m$ up to one entry.
It was already shown in~\cite{Wi62,GW69}, that $T_m$ provides quadrature nodes and weights
of Gauss-Radau quadrature formulae associated with the underlying distribution
(i.e.,~$\dd\alpha_n$ in the present setting).
In the context of Gauss-Radau quadrature formulae, the preassigned zero $\xi$
corresponds to a preassigned quadrature node, see also~\cite{Ga04,GM10}.

At the beginning of the present section we recall some theory on quasi-orthogonal polynomials.
In Subsection~\ref{sec.polynomialqorKry} this theory will be applied to the polynomial Krylov subspace~$\Kry_m(A,u)$.
While keeping the ortho\-normal basis~$V_m$ of~$\Kry_m(A,u)$ as before,
we consider the modified matrix $T_m$ (given by the underlying recursion; see~\eqref{eq.defJmxi} below)
as a representation of~$A$ in~$\Kry_m(A,u)$.
This results in the matrix decomposition~\eqref{eq.qormatrixdecomposition},
where 
$\widehat{p}_m$ provides the residual.
Thus, we also refer to $T_m$ as a
\begin{center}
{\em quasi-orthogonal residual \,(\textbf{\em qor-})\,Krylov representation.}
\end{center}
The zero $\xi\in\R$ of $\widehat{p}_m$ which is preassigned
constitutes an eigenvalue of the modified matrix $T_m$.
The spectrum of $T_m$ constitutes a step function $\alpha_m$
which is introduced properly in Section~\ref{sec.bounds3} below.
Based on the CMS Theorem, the distributions
$\dd\alpha_n$ and $\dd\alpha_m$ satisfy some intertwining property
(in general, this result is known for the Gauss-Radau quadrature rule;
in the Krylov setting we specify this result in Section~\ref{sec.bounds3} below).
In the qor-Krylov setting, we can make use of the preassigned zero $\xi$
to modify computable bounds on $\alpha_n$,
which potentially result in refined bounds.
Furthermore, we consider the matrix $T_m$ to approximate a matrix function $f(A)u$.
This is referred to as qor-Krylov approximation, see~\eqref{eq.qorfappr} below.
The qor-Krylov approximation can be understood as a corrected Krylov approximation,
comparable to the corrected Krylov scheme for the matrix exponential in~\cite{Sa92}.
In the context of approximating matrix functions,
making use of quasi-orthogonal polynomials is a
new idea.

Later on in this section the theory of quasi-orthogonal polynomials will be applied to the
case of a rational Krylov subspace.
We also refer to~\cite{Ga04} for rational Gauss-Radau quadrature formulae,
which are also applied in a Krylov setting in~\cite{LRW08,JR13}.
As in the polynomial case, we aim to refine estimates on $\alpha_n$ in the sequel,
and we also introduce a rational qor-Krylov approximation.
In Remark~\ref{rmk:BmforratKry} below, we introduce a new procedure
to efficiently compute the rational qor-Krylov representation, i.e.,
we rewrite a rational Krylov subspace with arbitrary complex poles
as an extended Krylov subspace with a modified initial vector,
and then construct the rational qor-Krylov representation
based on results for the polynomial case.

\smallskip

We proceed to recall some theory on quasi-orthogonal polynomials.
Let~$p_0,\ldots,p_m$ be the sequence of orthonormal polynomials
from Proposition~\ref{prop.threetermorthopol}.
Let $\beta_{m-1}>0$ be given as in Proposition~\ref{prop.threetermorthopol},
and let~$\omega_m\in\R$ (to be fixed in the sequel,
see~\eqref{eq.qorgetomega}). We define
\begin{subequations} \label{eq.defphat}
\begin{equation}\label{eq.threetermrecursionphat}
\widehat{p}_{m}(\lambda) = (\lambda - \omega_{m}) p_{m-1}(\lambda) - \beta_{m-1} p_{m-2}(\lambda) .
\end{equation}
The polynomial $p_m$ satisfies the recursive identity~\eqref{eq.threetermrecursionp} (for $ j=m $).
Thus, $\widehat{p}_{m}$ can be expressed as a linear combination of~$p_{m-1}$ and~$p_m$,
\begin{equation}\label{eq.defquasiorthop}
\widehat{p}_m = \beta_m\,p_m + (a_m-\omega_m){p}_{m-1},~~~\text{hence,~~\;$\widehat{p}_m\perp p_0,\ldots,p_{m-2}$}.
\end{equation}
\end{subequations}
With the orthogonality property~\eqref{eq.defquasiorthop}
we refer to~$\widehat{p}_{m}$
as quasi-orthogonal polynomial of degree~$m$.\footnote{
   In the case $a_m=\omega_m$ the polynomial $\widehat{p}_m$
   in~\eqref{eq.threetermrecursionphat} is identical to $\beta_m p_m$,
   thus, $\widehat{p}_m$ is an orthogonal polynomial.}

According to the requirement $\widehat{p}_{m}(\xi)=0$ imposed above
(see~\eqref{eq.pdach(xi)=0}) for a given
$\xi\in\R$ with $p_{m-1}(\xi)\neq 0$,
definition~\eqref{eq.threetermrecursionphat} implies
\begin{subequations} \label{eq.0=pmdach(xi)-omega_m}
\begin{equation}\label{eq.pmhatispmpluspm1}
0=\widehat{p}_{m}(\xi) =  (\xi-\omega_{m}) {p}_{m-1}(\xi) - \beta_{m-1} {p}_{m-2}(\xi).
\end{equation}
This fixes the value of $ \omega_m $,
\begin{equation}\label{eq.qorgetomega}
\omega_{m} = \xi - \beta_{m-1} \frac{{p}_{m-2}(\xi)}{p_{m-1}(\xi)}.
\end{equation}
\end{subequations}

We now reuse the denotation~$\theta_1,\ldots,\theta_m$
in a modified way:
In the context of quasi-orthogonal polynomials, $\theta_1,\ldots,\theta_m\in\R$
denote the zeros of~$\widehat{p}_m$.
We assume the ordering $ \theta_1<\theta_2<\ldots<\theta_m $.

\begin{proposition}[See also Section 3.3 in~\cite{Sze85}]\label{prop.ritzqorpol}
Let~$\widehat{p}_m$ be the quasi-orthogonal polynomial
defined in~\eqref{eq.threetermrecursionphat},
with~$\omega_m$ from~\eqref{eq.qorgetomega} for a given~$\xi\in\R$ with $p_{m-1}(\xi)\neq 0$.
\begin{enumerate}[label=(\roman*)]
\item
The zeros $ \theta_1,\ldots,\theta_m $ of $ {\widehat p}_m $ are distinct.
\item
Interlacing property of eigenvalues~$\lambda_1,\ldots,\lambda_n$
and zeros of~$\widehat{p}_m$:
For~$k=1,\ldots,m-1$ there exists at least one~$\lambda_{j(k)}$ with
$$
\theta_{k} < \lambda_{j(k)} < \theta_{k+1}.
$$
\item
At most one of the zeros~$\theta_1,\ldots,\theta_m$ is located outside
of~$[\lambda_1,\lambda_n]$. 
E.g., in the case $\xi<\lambda_1$ we have~$\theta_1<\lambda_1<\theta_2<\ldots<\theta_m < \lambda_n$.
\end{enumerate}
\end{proposition}


As a slight modification of the Jacobi matrix~$ J_m $ from~\eqref{eq.defJm}
we now define the symmetric tridiagonal matrix
\begin{equation}\label{eq.defJmxi}
T_m =
\begin{pmatrix}
a_1&\beta_1&&&\\
\beta_1&a_2&\beta_2&&\\
  &\ddots&\ddots&\ddots&\\
  &&\beta_{m-2}&{a}_{m-1}&\beta_{m-1}\\
  &&&\beta_{m-1}& \bm{\omega_{m}}
\end{pmatrix}\in\R^{m\times m},
\quad \text{with \,$\omega_m$~\,from~\eqref{eq.qorgetomega}.}
\end{equation}
With the recursion~\eqref{eq.threetermrecursionp} and identity~\eqref{eq.defquasiorthop}
the sequence of orthonormal polynomials~$P(\lambda)=(p_0(\lambda),\ldots,p_{m-1}(\lambda))^\Hast\in\R^m$
and~$\widehat{p}_m$ satisfy
\begin{equation}\label{eq.qorpmatrixform}
\lambda P(\lambda) = T_{m}\,P(\lambda) + \widehat{p}_{m}(\lambda)e_m.
\end{equation}
Thus, the eigenvalues of~$T_m$ are exactly the zeros~$\theta_1,\ldots,\theta_m$ of~$\widehat{p}_m$.

\subsection{Krylov methods and quasi-orthogonal polynomials}\label{sec.polynomialqorKry}

Let~$p_0,\ldots,p_m$ be the orthonormal polynomials from
Proposition~\ref{prop.threetermorthopol},
which provide the \tMinr-orthonormal Krylov basis vectors~$v_{j}=p_{j-1}(A)u$
for~$j=1,\ldots,m+1$,
and let~$\widehat{v}_{m+1}=\widehat{p}_{m}(A)u$ with
the quasi-orthogonal polynomial~$\widehat{p}_m$ from~\eqref{eq.defphat}.
Analogously to~\eqref{eq.qorpmatrixform} we have the matrix decomposition
\begin{equation}\label{eq.qormatrixdecomposition}
A\,V_m = V_m\,T_m + \widehat{v}_{m+1}\,e_m^\Hast.
\end{equation}
We refer to~$\widehat{v}_{m+1}\in\C^n$
as the residual of~\eqref{eq.qormatrixdecomposition}, with
$$
\widehat{v}_{m+1} \in \vspan\{v_m,v_{m+1} \} \perp_{\Minr} \Kry_{m-1}(A,u).
$$
\begin{proposition}\label{prop.eqforqorpolyKry}
For~$p\in \Pi_{m-1}$,
\begin{equation}\label{eq.quasiKrylovidentity}
\beta_0 V_m\,p(T_m) e_1 = p(A)u.
\end{equation}
\end{proposition}
\begin{proof}
We prove $\beta_0 V_m T_m^j\,e_1 = A^j\,u$ for $j=0,\ldots,m-1$ by induction.
This holds true for $j=0$.
Assuming that it also holds true for some $j<m-1$, then
\begin{equation*} 
A^{j+1} u = A\,A^j u = \beta_0\,A\,V_m T_m^j\,e_1.
\end{equation*}
Together with identity~\eqref{eq.qormatrixdecomposition} this gives
\begin{equation*}
A^{j+1} u = \beta_0 V_m T_m^{j+1} e_1 + \beta_0 \widehat{v}_{m+1} e_m^\Hast T_m^j\,e_1.
\end{equation*}
Due to the tridiagonal structure of $T_m$ we have
$e_m^\Hast\,T_m^j\,e_1=0$ for $j=0,\ldots,m-2$.
Altogether, this implies
$\beta_0 V_m\,T_m^j\,e_1 = A^j\,u$
for $j=0,\ldots,m-1$, which completes the proof
of~\eqref{eq.quasiKrylovidentity}.
\end{proof}
In addition to Proposition~\ref{prop.eqforqorpolyKry}
we note that for~$p\in\Pi_m$ exactly of degree~$m$,
\begin{equation*}
\beta_0 V_m\,p(T_{m}) e_1 - p(A)u \in \vspan\{v_{m},v_{m+1} \} \perp_{\Minr} \Kry_{m-1}(A,u).
\end{equation*}

The following proposition is associated with identities of Gauss-Radau quadrature formulae, see also~\cite{GM10} or~\cite[Subsection 3.1.4]{Ga04}.
This relation is discussed in more detail in Section~\ref{sec.bounds3} below.
\begin{proposition}\label{prop.eqforqorpolyKryinner}
For~$p\in \Pi_{2m-2}$,
\begin{equation}\label{eq.quasiKrylovidentityinner}
(u,p(A)u)_{\Minr} = \beta_0^2(e_1,p(T_m)e_1)_2.
\end{equation}
\end{proposition}
\begin{proof}
We write $p=g_1\,g_2$ with $g_1,g_2\in\Pi_{m-1}$
and apply Proposition~\ref{prop.eqforqorpolyKry} to both terms,
$$
(u,p(A)u)_{\Minr}=(\tbar{g}_1(A)u,g_2(A)u)_{\Minr}
=\beta_0^2(V_m\,\tbar{g}_1(T_m)e_1,V_m\,g_2(T_m)e_1)_{\Minr}.
$$
With $(V_m,V_m)_{\Minr}=I$ this implies~\eqref{eq.quasiKrylovidentityinner}.
\end{proof}

We proceed by recapitulating
results from~\cite[Subsection 6.2.1]{GM10} and~\cite[Section 7]{Go73}
which reveal an algorithm to construct~$T_m$.
\begin{remark}[\cite{GM10,Go73}]\label{rmk.GolubcomputeTm}
Let $J_{m-1}$ be the Jacobi matrix constructed by $m-1$ steps of the Lanczos method.
After substituting $\xi$ for $\lambda$ in~\eqref{eq.threetermrecursionpmatrix},
the Jacobi matrix~$J_{m-1}$
and~$P(\xi)=(p_0(\xi),\ldots,p_{m-2}(\xi))^\Hast\in\R^{m-1}$ satisfy
$$
(J_{m-1}-\xi I)P(\xi) = -\beta_{m-1}\,p_{m-1}(\xi) e_{m-1}.
$$
The solution~$\delta=(\delta_1,\ldots,\delta_{m-1})\in\R^{m-1}$ of
the linear system
\begin{equation}\label{eq.compomegamgetdelta}
(J_{m-1} - \xi I)\delta = \beta_{m-1}^2 e_{m-1}
\end{equation}
is given by
$$
\delta_\ell = -\beta_{m-1} \frac{p_{\ell-1}(\xi)}{p_{m-1}(\xi)},
~~~\ell=1,\ldots,m-1.
$$
The eigenvalues of~$J_{m-1}$ are identical to the zeros of $p_{m-1}$,
hence, with $p_{m-1}(\xi)\neq 0$ the matrix~$(J_{m-1} - \xi I)$ is invertible.
The solution $\delta\in\R^{m-1}$ of~\eqref{eq.compomegamgetdelta}
yields a computable formula for $\omega_m$ via~\eqref{eq.qorgetomega},
$$
\omega_m = \xi + \delta_{m-1}.
$$
\end{remark}
Algorithm~\ref{alg.qorKrypol} represents a summary on Remark~\ref{rmk.GolubcomputeTm}.
In~\figref{fig.quasiortwmoverxi} we show values of $\omega_m$
over $\xi$ for a given example.

\begin{figure}
\centering
\begin{algorithm}[H]
\caption{An algorithm to compute $V_m$ and
the qor-Krylov representation $T_m$ for a given $\xi\in\R$
which is distinct to the eigenvalues of $J_{m-1}$.}
\label{alg.qorKrypol}
\SetAlgoLined
 apply the Lanczos method for~$\Kry_{m-1}(A,u)$: this returns~$\beta_0$,~$V_{m-1},J_{m-1},\beta_{m-1},v_m$\;
 set~$\omega_m = \xi + \beta_{m-1}^2 e_{m-1}^\Hast (J_{m-1} - \xi I)^{-1} e_{m-1}$ and define~$T_m$ via~\eqref{eq.defJmxi}\;
 set~$V_m=(V_{m-1},v_m)$\;
 return~$\beta_0, V_m, T_m$\;
\end{algorithm}
\end{figure}

\begin{figure}
\centering
\begin{overpic}
[width=0.8\textwidth]{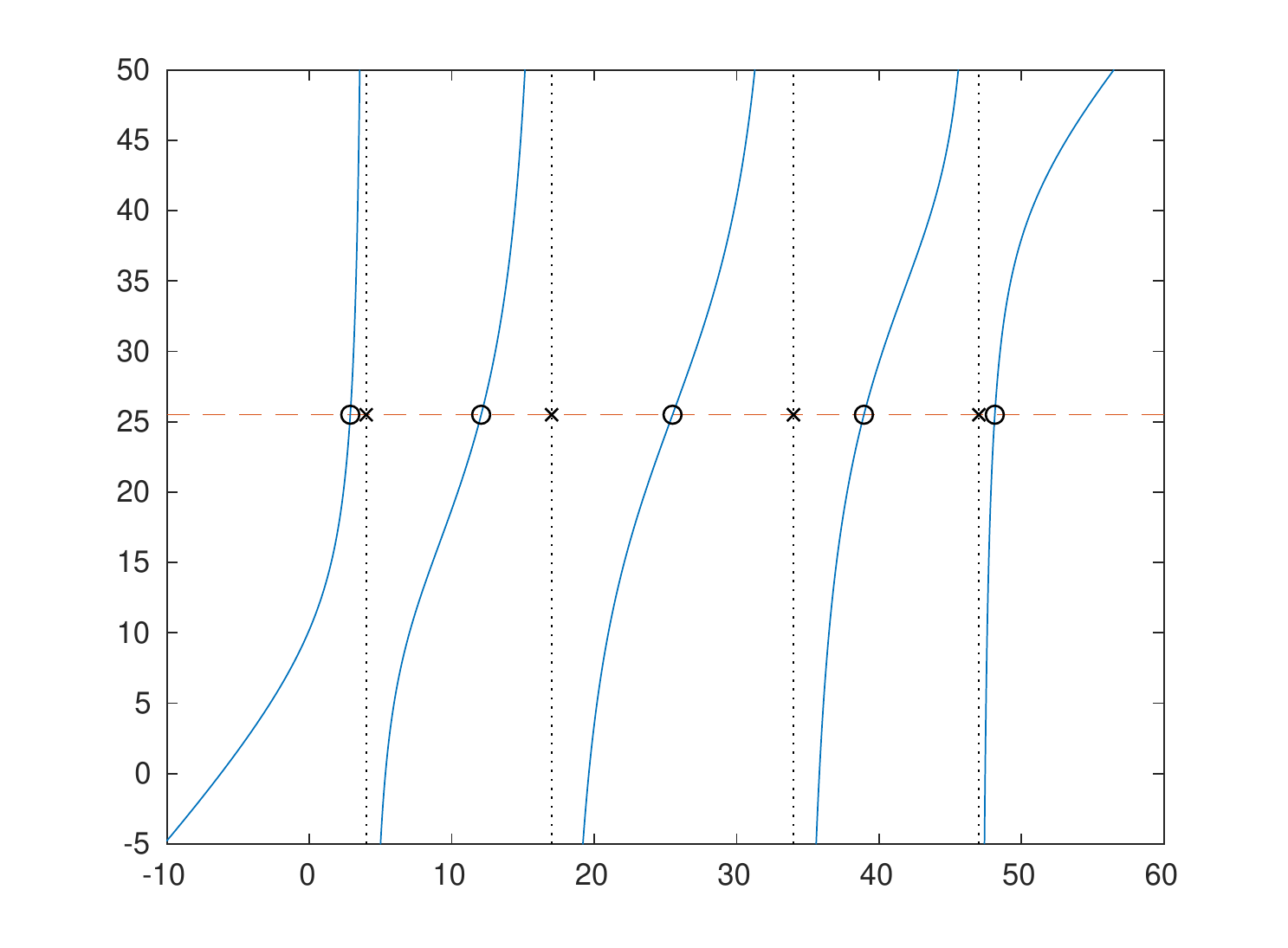}
\put(1,38){$\omega_m$}
\put(14.5,43){\small$a_m\!\!=\!\omega_m$}
\put(50,0){$\xi$}
\end{overpic}
\caption{This figure shows the matrix entry $\omega_m$ of the qor-Krylov
representation~$T_m$, computed for different values of~$\xi\in\R$
and~$m=5$.
To compute the entries~$\omega_m$ we follow Algorithm~\ref{alg.qorKrypol}.
As an example we choose~$A$ to be a~$n\times n$ diagonal matrix
with~$n=50$ and diagonal entries~$(1,\ldots,n)$,
and we choose~$u=(1,\ldots,1)^\Hast\in\R^{n}$.
When the choice of~$\xi$ matches one of the eigenvalues of~$J_m$
(marked by~('$\circ$')),
then the matrices~$T_m$ and~$J_m$ coincide (the matrix entry~$a_m$ of~$J_m$
is illustrated by the dashed horizontal line).
On the other hand, when $\xi$ coincides with an eigenvalue of~$J_{m-1}$
(marked by~('$\times$') and dotted vertical lines),
then~$\omega_m$ is undefined and
Algorithm~\ref{alg.qorKrypol} fails.
We remark that two neighboring eigenvalues of~$J_m$
enclose exactly one eigenvalue of~$J_{m-1}$,
cf.~\cite[Theorem~3.3.2]{Sze85}.
This property carries over to the eigenvalues of~$T_m$
via~\eqref{eq.defquasiorthop}
(indeed, the sign of~$\widehat{p}_m$ corresponds to the sign of~$p_m$
at the zeros of~$p_{m-1}$ and at the boundary of~$\R$).
Thus, two neighboring eigenvalues of~$T_m$
enclose exactly one eigenvalue of~$J_{m-1}$ for any valid choice of~$\xi$.}
\label{fig.quasiortwmoverxi}
\end{figure}

\paragraph{A quasi-orthogonal residual (qor-)Krylov approximation
           to matrix functions $ f(A)u $.}
$~$\newline
We refer to
\begin{equation} \label{eq.qorfappr}
\beta_0 V_m f(T_m)e_1 \approx f(A)u
\end{equation}
as quasi-orthogonal residual (qor-)Krylov approximation,
based on the construction of $ V_m $ and $ T_m $ according to
Algorithm~\ref{alg.qorKrypol}.
We recall that only $m-1$ steps of the Lanczos iteration are required.
This provides
the orthonormal basis $V_{m-1}$, the subsequent basis vector $v_m$
and the Jacobi matrix $J_{m-1}$.
The qor-Krylov approximation makes use of
the orthonormal basis $V_m=(V_{m-1},v_m)$, where the polynomial 
Krylov approximation,~i.e,~$\beta_0 V_{m-1} f(J_{m-1})e_1 \approx f(A)u$,
provides an approximation in the basis $V_{m-1}$.

The idea to `correct' the Krylov approximation
by including the subsequent basis vector
(which is $v_m$ at the {\small$(m-1)$}-th step) also appears in~\cite{Sa92},
namely, the corrected Krylov scheme for the matrix exponential
which is widely used in the Expokit package~\cite{Si98} and others.
Compared to the corrected Krylov scheme,
the qor-Krylov approximation can be favorable if spectral properties of $f(A)$
are relevant,~e.g.,~the mass conservation of $\ee^{-\ii tA}u$
carries over to the qor-Krylov approximation $\beta_0 V_m \ee^{-\ii t T_m}e_1$
due to $T_m$ being Hermitian.


\subsection{Rational Krylov methods and the theory of
quasi-orthogonal polynomials}\label{sec.ratqor}

A rational Krylov subspace satisfies~$\RKry_m(A,u)=\Kry_m(A,u_q)$
for~$u_q=q_{m-1}^{-1}(A) u$, where $q_{m-1}$ denotes
the denominator given by preassigned poles.
Let $J_m$ and $V_m$ denote the Jacobi matrix and \tMinr-orthonormal basis of~$\Kry_m(A,u_q)$.
The procedure of Subsection~\ref{sec.polynomialqorKry} applies to the polynomial Krylov subspace~$\Kry_m(A,u_q)$:
For a given~$\xi\in\R$ the matrix $T_m$ is defined in~\eqref{eq.defJmxi}
and satisfies the matrix decomposition~\eqref{eq.qormatrixdecomposition}
together with $A$ and $V_m$.

In a practical setting, $u_q=q_{m-1}^{-1}(A) u$ is not directly available
to construct the rational Krylov subspace via~$\Kry_m(A,u_q)$.
We proceed to generalize the qor-Krylov representation for the rational Krylov subspace:
Let $U_m$ be a given \tMinr-orthonormal basis of~$\RKry_m(A,u)$,
i.e.,~$\vspan\{U_m\}=\RKry_m(A,u)$ and $(U_m,U_m)_{\Minr}=I$.
The respective Rayleigh quotient is~$A_m=(U_m,A\,U_m)_{\Minr}$.
With the orthonormal transformation~$K_m = (V_m,U_m)_{\Minr} \in\C^{m\times m}$
we have~$A_m =  K_m^\Hast\,J_m\,K_m$ as given in~\eqref{eq.RKAmtoJm}.
For a representation of $T_m$ in the basis $U_m$ we introduce the notation
\begin{equation}\label{eq.defBm}
B_m =  K_m^\Hast\,T_m\,K_m.
\end{equation}
The eigenvalues of $B_m$
are equal to the eigenvalues $\theta_1,\ldots,\theta_m\in\R$ of $T_m$
and satisfy Proposition~\ref{prop.ritzqorpol}.
The Hermitian structure of $T_m$ carries over to $B_m$.

\begin{proposition}\label{prop.idforqorratKry}
With $x=(U_m,u)_{\Minr}$ we have
\begin{equation}\label{eq.ratkryqoridr}
r(A)u = U_m\,r(B_m) x~~~\text{for~\;$ r\in\Pi_{m-1}/q_{m-1}$}.
\end{equation}
\end{proposition}
\begin{proof}
Let $\zeta_0=\|u_q\|_{\Minr} $, let $V_m$ be the \tMinr-orthonormal basis of $\Kry_m(A,u_q)$,
and let $T_m$ be the respective qor-Krylov representation for a given $\xi\in\R$.
Then Proposition~\ref{prop.eqforqorpolyKry} w.r.t.\ $\Kry_m(A,u_q)$ implies
\begin{equation}\label{eq.ratkryqoridp1}
p(A)u_q =\zeta_0\,V_m\,p(T_m) e_1,~~~p\in \Pi_{m-1}.
\end{equation}
This implies $q_{m-1}(A)u_q=\zeta_0 V_m\,q_{m-1}(T_m) e_1$,
and with the identities~$q_{m-1}(A)u_q=u$ and $(V_m,V_m)_{\Minr}=I$ we arrive at
\begin{equation}\label{eq.ratkryqoridp2}
\zeta_0\,e_1 = q_{m-1}^{-1}(T_m) (V_m,u)_{\Minr}.
\end{equation}
Let $r=p/q_{m-1}$ for $p\in\Pi_{m-1}$ then $r(A)u=p(A)u_q$, and with~\eqref{eq.ratkryqoridp1} we have
\begin{equation}\label{eq.ratkryqoridp3}
r(A)u = \zeta_0\,V_m\,p(T_m) e_1.
\end{equation}
Inserting~\eqref{eq.ratkryqoridp2} into~\eqref{eq.ratkryqoridp3} gives
\begin{equation}\label{eq.ratkryqoridpx}
r(A)u = V_m\,p(T_m) q_{m-1}^{-1}(T_m) (V_m,u)_{\Minr}
= V_m\,r(T_m) (V_m,u)_{\Minr}.
\end{equation}
With~$ K_m\,K_m^\Hast = I $ (see~\eqref{eq.Kmortho1}) the matrix $B_m$ in~\eqref{eq.defBm}
satisfies~$r(T_m)=K_m r(B_m) K_m^\Hast$,
and together with~$V_m K_m = U_m$~\eqref{eq.Kmortho1} we have
\begin{equation}\label{eq.ratkryqoridp4}
V_m\,r(T_m) (V_m,u)_{\Minr} = U_m\,r(B_m) (U_m,u)_{\Minr}.
\end{equation}
Combining~\eqref{eq.ratkryqoridpx} with~\eqref{eq.ratkryqoridp4}
results in~\eqref{eq.ratkryqoridr}.
\end{proof}

The following proposition is associated with identities of rational Gauss-Radau quadrature formulae, see also~\cite[§\,3.1.4.4]{Ga04}.
For more details on this relation see Section~\ref{sec.bounds3} below.
\begin{proposition}\label{prop.idforqorratKryinner}\label{prop.eqforqorratKry}
With $x=(U_m,u)_{\Minr}$,
\begin{equation}\label{eq.ratkryqoreqinnerprod}
(u,r(A)u)_{\Minr} = (x,r(B_m)x)_2~~~\text{for~\;$ r\in\Pi_{2m-2} /|q_{m-1}|^2$}.
\end{equation}
\end{proposition}
\begin{proof}
For rational functions~$r\in \Pi_{2m-2}/|q_{m-1}|^2$ we write $r=r_1\,r_2$,
where $r_1\in\Pi_{m-1}/\tbar{q}_{m-1}$ and $r_2\in\Pi_{m-1}/q_{m-1}$.
With this notation we write
\begin{equation}\label{eq.ratkryeqinnerprodp1x}
(u,r(A)u)_{\Minr} = (\tbar{r}_1(A)u,r_2(A)u)_{\Minr},
~~~\text{and}~~~
(x,r(B_m)x)_2 = (\tbar{r}_1(B_m) x,r_2(B_m) x)_2.
\end{equation}
For $r_1,\,r_2\in\Pi_{m-1}/q_{m-1}$ we apply Proposition~\ref{prop.idforqorratKry}
to conclude
\begin{equation}\label{eq.ratkryeqinnerprodp2x}
(\tbar{r}_1(A)u,r_2(A)u)_{\Minr}
= (U_m\,\tbar{r}_1(B_m)x,U_m\,r_2(B_m)x)_{\Minr}.
\end{equation}
Combining~\eqref{eq.ratkryeqinnerprodp1x}
with~\eqref{eq.ratkryeqinnerprodp2x} and making use of
$(U_m,U_m)_{\Minr}=I$ we conclude~\eqref{eq.ratkryqoreqinnerprod}.
\end{proof}

The definition of $B_m$ in~\eqref{eq.defBm} is of a theoretical nature.
We propose a setup in which $B_m$ can be computed efficiently.
Let $\RKry_{m-2}(A,u)$ be a rational Krylov subspace with
arbitrary poles~$s_1,\ldots,s_{m-3}\in\C\cup \{\pm\infty\}$.
The poles $s_1,\ldots,s_{m-3}$ define the denominator $q_{m-3}$ and we write $u_q=q_{m-3}^{-1}(A) u$.
We also recall the identities
$$
\RKry_{m-2}(A,u)=\Kry_{m-2}(A,u_q)~~~\text{and}~~
\Kry_m(A,u_q)=\Kry_{m-2}(A,u_q)\oplus\vspan\{A\,u,A^2u \}.
$$
We extend the rational Krylov subspace $\RKry_{m-2}(A,u)$
by two additional polynomial Krylov steps,
i.e.,
\begin{equation}\label{eq.xtendedKSinprop}
\Kry_m(A,u_q) = \RKry_{m-2}(A,u) \,\oplus\, \vspan\{A\,u,A^2u \},
~~~~\text{where}~~~u_q=q_{m-3}^{-1}(A) u.
\end{equation}

The Krylov subspace~$\Kry_m(A,u_q)$ can be referred to as an extended Krylov subspace,
and some of the following results are related to~\cite[Section 5]{DK98}.

\begin{proposition}\label{prop.extendedKrylovalgorithmrmk}
Let $m$ be fixed and $u_q=q_{m-3}^{-1}(A) u$ for a given denominator~$q_{m-3}$.
Let $U_{m-2}\in\C^{n\times m-2}$ be a given
\tMinr-orthonormal basis of~$\RKry_{m-2}(A,u)=\Kry_{m-2}(A,u_q)$,
and $A_{m-2}=(U_{m-2},A\,U_{m-2})_{\Minr}$.
Let~$\Kry_m(A,u_q)$ refer to the
extended Krylov subspace given in~\eqref{eq.xtendedKSinprop}.
Let $V_m=(v_1,\ldots,v_m)\in\C^{n\times m}$ be the
\tMinr-orthonormal basis of~$\Kry_m(A,u_q)$ provided by the Lanczos method.
Then the following statements hold true and provide a procedure
to compute the qor-Krylov representation of $B_m$
for a given $\xi\in\R$ and the basis $\widetilde{U}_m$ (given below)
of the extended Krylov subspace $\Kry_m(A,u_q)$.
\begin{enumerate}[label=(\roman*)]
\item\label{item.extendedKrylovalgorithmrmki1}
With~$\widetilde{U}_m=(U_{m-2},v_{m-1},v_m)\in\C^{n\times m}$ we have
an \tMinr-orthonormal basis of the extended Krylov subspace $\Kry_{m}(A,u_q)$,
i.e.,~$\vspan\{\widetilde{U}_m\}=\Kry_{m}(A,u_q)$ and $(\widetilde{U}_m,\widetilde{U}_m)_{\Minr} = I$.
Furthermore, $\widetilde{U}_m$ can be computed
without reference to~$u_q$.

\item
The Rayleigh quotient
$\widetilde{A}_m = (\widetilde{U}_m,A\,\widetilde{U}_m)_{\Minr}$
of the extended Krylov subspace is given by
\begin{align}
&\widetilde{A}_m=\begin{pmatrix}
\widetilde{A}_{m-1}&\beta_{m-1}e_{m-1}\\
\beta_{m-1}e_{m-1}^\Hast &a_m
\end{pmatrix}\in\C^{m\times m},
~~~\text{with}\label{eq.specialAmratqor}\\
&\widetilde{A}_{m-1}=\begin{pmatrix}
A_{m-2}&\widetilde{a}\\
\widetilde{a}^\Hast&a_{m-1}
\end{pmatrix}\in\C^{m-1\times m-1},
~~~
\text{and}~~~
\widetilde{a} = (U_{m-2},A\,v_{m-1})_{\Minr}\in\C^{m-2}.\notag
\end{align}
Furthermore, $a_m=(J_m)_{m,m}$, $a_{m-1}=(J_m)_{m-1,m-1}$ and $\beta_{m-1}=(J_m)_{m,m-1}$
for the Jacobi matrix $J_m$ of $\Kry_m(A,u_q)$.
The matrix entries~$a_m$, $a_{m-1}$, $\beta_{m-1}$, and $\widetilde{a}$
are computed in course of the orthogonalization procedure in~\ref{item.extendedKrylovalgorithmrmki1}.

\item
For the basis transformation $\widetilde{K}_m=(V_m,\widetilde{U}_m)_{\Minr}$ we have
\begin{equation}\label{eq.specialKm}
\widetilde{K}_m = \begin{pmatrix}
K_{m-2}&0\\
0&I_2
\end{pmatrix},
~~~\text{with}~~\;
I_2 = \begin{pmatrix}
1&0\\
0&1
\end{pmatrix},
\end{equation}
and $K_{m-2}=(V_{m-2},U_{m-2})_{\Minr}$.

\item\label{item.extendedKrylovalgorithmrmkilast}
Let $T_m$ be defined by~\eqref{eq.defJmxi} for~$\Kry_m(A,u_q)$.
Then, $\widetilde{B}_m=\widetilde{K}_m^\Hast\,T_m\,\widetilde{K}_m$ satisfies

\begin{subequations} 
\begin{equation}\label{eq.initemBmqorratKry}
\widetilde{B}_m = \begin{pmatrix}
\widetilde{A}_{m-1}&\beta_{m-1}\,e_{m-1}\\
\beta_{m-1}\,e_{m-1}^\Hast &\omega_m
\end{pmatrix},
\end{equation}
with
\begin{equation}\label{eq.initemomegamqorratKry}
\omega_m = \xi + \beta_{m-1}^2  e_{m-1}^\Hast ( \widetilde{A}_{m-1}  - \xi I)^{-1}e_{m-1}.
\end{equation}
\end{subequations}

\end{enumerate}
\end{proposition}
\begin{proof}
$~$ \newline \vspace{-5mm}
\begin{enumerate}[label=(\roman*)]
\item We have $\vspan\{U_{m-2}\}=\Kry_{m-2}(A,u_q)=\vspan\{V_{m-2}\}$,
and by adding $v_{m-1}$ and $v_m$ to the basis we have $\vspan\{\widetilde{U}_m\}=\Kry_{m}(A,u_q)$.
With $(U_{m-2},U_{m-2})_{\Minr}=I$ and $v_{m-1},v_m \perp_{\Minr} \Kry_{m-2}(A,u_q)$
this also implies $(\widetilde{U}_m,\widetilde{U}_m)_{\Minr} = I$.

The \tMinr-orthonormal basis $\widetilde{U}_m$ can be constructed without
referring to $u_q$ by the following procedure:
We construct $v_{m-1}$ by orthogonalizing $\widetilde{v} = A\,u$
w.r.t.~$U_{m-2}$ and normalizing.
In a similar manner we construct $v_m$ via $A v_{m-1}$.
To demonstrate that this procedure yields the correct results,
we argue as follows:
We recall $u_q = q_{m-3}^{-1}(A)u$, hence, $u = q_{m-3}(A)u_q$.
We introduce the notation $\widetilde{p}(\lambda)=\lambda\,q_{m-3}(\lambda)$,
where $\lambda\,q_{m-3}(\lambda)=\lambda^{m-2} +\widetilde{p}_{m-3}(\lambda)$
for a polynomial~$\widetilde{p}_{m-3}\in\Pi_{m-3}$.
Let $\widetilde{v} = A\,u$, then $\widetilde{v} =  \widetilde{p}(A)u_q$.
We recall $v_j=p_{j-1}(A)u_q$ for the orthonormal polynomials $p_0,\ldots,p_m$
provided by Proposition~\ref{prop.threetermorthopol}.
The polynomials $\widetilde{p}$ and $p_{m-2}$ both have a positive real-valued leading coefficient.
Hence, we obtain $p_{m-2}$ by orthogonalizing $\widetilde{p}$
w.r.t.~$p_0,\ldots,p_{m-3}$ and normalizing.
Analogously, we obtain $v_{m-1}$ by orthogonalizing
$\widetilde{v} = A\,u$ w.r.t.~$U_{m-2}$ and normalizing as stated above.

\item
In order to specify~$\widetilde{A}_m = (\widetilde{U}_m,A\,\widetilde{U}_m)_{\Minr}$
we recall~$\widetilde{U}_m=(U_{m-2},v_{m-1},v_m)$.
The upper left submatrix of~$\widetilde{A}_m$
is given by $A_{m-2} = (U_{m-2},A\,U_{m-2})_{\Minr}$.
In a similar manner we deduce $\widetilde{a}$, $a_m$, $a_{m-1}$ and $\beta_{m-1}$.
Here $a_m=(v_m,A\,v_m)_{\Minr}$ by the structure of $\widetilde{U}_m$
and with $J_m=(V_m,A\,V_m)_{\Minr}$ we also have $a_m= (J_m)_{m,m}$.
Analogously, $a_{m-1}$ and $\beta_{m-1}$ are equal to entries of $J_m$.
We introduce the notation $\widetilde{a} = (U_{m-2}, A\,v_{m-1})_{\Minr}\in\C^{m-2}$.
The entries $(U_{m-2},A\,v_m)_{\Minr}$ are zero due to
$A\,v_{m} \in\vspan\{v_{m-1},v_m,v_{m+1}\}$ being \tMinr-orthogonal to $U_{m-2}$.

\item
The basis transformation $\widetilde{K}_m=(V_m,\widetilde{U}_m)_{\Minr}$
for $\widetilde{U}_m=(U_{m-2},v_{m-1},v_m)$ and $V_m=(v_1,\ldots,v_m)$,
where $\widetilde{U}_m$ and $V_m$ are \tMinr-orthonormal bases,
indeed has the simple structure~\eqref{eq.specialKm}.

\item
We proceed with the matrix entry $\omega_m$ of $T_m$ in~\eqref{eq.defJmxi}.
Following Algorithm~\ref{alg.qorKrypol}, $\omega_m$~evaluates to
\begin{equation}\label{eq.omegaTm2}
\omega_m = \xi + \beta_{m-1}^2 e_{m-1}^\Hast(J_{m-1} - \xi I)^{-1}e_{m-1},
\end{equation}
where $\beta_{m-1}$ refers to $(J_m)_{m,m-1}$ which is equal
to $(\widetilde{A}_m)_{m,m-1}$, see~\eqref{eq.specialAmratqor}.
By the matrix structure of $\widetilde{K}_m$ (see~\eqref{eq.specialKm}) we have $\widetilde{K}_{m-1}e_{m-1}=e_{m-1}$,
thus,
\begin{equation}\label{eq.TminvbyAms0}
e_{m-1}^\Hast(J_{m-1} - \xi I)^{-1}e_{m-1}
= e_{m-1}^\Hast \widetilde{K}_{m-1}^\Hast ( J_{m-1}  - \xi I)^{-1}\widetilde{K}_{m-1} e_{m-1}
\end{equation}
Furthermore, $\widetilde{K}_{m-1}^\Hast \widetilde{K}_{m-1} = I$~\eqref{eq.Kmortho1}
together with $\widetilde{A}_{m-1} = \widetilde{K}_{m-1}^\Hast\,J_{m-1}\,\widetilde{K}_{m-1} $
yield
\begin{equation}\label{eq.TminvbyAm}
e_{m-1}^\Hast \widetilde{K}_{m-1}^\Hast ( J_{m-1}  - \xi I)^{-1}\widetilde{K}_{m-1} e_{m-1}
= e_{m-1}^\Hast ( \widetilde{A}_{m-1}  - \xi I)^{-1}e_{m-1}.
\end{equation}
Combining~\eqref{eq.omegaTm2} with~\eqref{eq.TminvbyAms0}
and~\eqref{eq.TminvbyAm} we conclude~\eqref{eq.initemomegamqorratKry}.

Compare $J_m$~\eqref{eq.defJm} with $T_m$~\eqref{eq.defJmxi} to observe
$$
T_m = J_m + (\omega_m-a_m) e_m e_m^\Hast.
$$
With $\widetilde{K}_m^\Hast\,e_{m}= e_{m}$ and
$\widetilde{A}_{m} = \widetilde{K}_{m}^\Hast\,J_{m}\,\widetilde{K}_{m}$
this implies
\begin{equation}\label{eq.TmtoBm}
\widetilde{B}_m =  \widetilde{K}_m^\Hast\,T_m\,\widetilde{K}_m
= \widetilde{A}_m + (\omega_m-a_m) e_m e_m^\Hast.
\end{equation}
With~\eqref{eq.specialAmratqor} and~\eqref{eq.TmtoBm} we conclude~\eqref{eq.initemBmqorratKry}.
\end{enumerate}
\end{proof}

\begin{remark}\label{rmk:BmforratKry}
The approach of Proposition~\ref{prop.extendedKrylovalgorithmrmk}
provides $B_m$ for an extended Krylov subspace
and can be slightly modified to fit for a fully rational Krylov subspace $\RKry_m(A,u)$.
Let $s_1,\ldots,s_{m-1}\in\C\cup\{\pm\infty\}$, where $s_{m-2},s_{m-1}\in\C$, and let $q_{m-1}$ be the respective denominator.
We recall
\begin{equation}\label{eq.exKryp2x1purec}
\RKry_m(A,u) = \Kry_m(A,u_q),~~~\text{where}~~u_q=q_{m-1}^{-1}(A)u.
\end{equation}
We introduce the modified initial vector $\widehat{u}$ and denominator $\widehat{q}_{m-3}(A)$ as
\begin{equation}\label{eq.exKryp2x1puhat}
\begin{aligned}
\widehat{u} &= (A-s_{m-2}I)^{-1}(A-s_{m-1}I)^{-1} u,~~~\text{and}~~\\
\widehat{q}_{m-3}(A) &=  (A-s_1 I) (A-s_2 I) \cdots (A-s_{m-3}I) .
\end{aligned}
\end{equation}
Let~$\RKry_{m-2}(A,\widehat{u})$ be the rational Krylov subspace according to the 
initial vector $\widehat{u}$ and poles $s_1,\ldots,s_{m-3}$.
Then
$$
\RKry_{m-2}(A,\widehat{u} \hskip 1pt) = \Kry_{m-2}(A, \widehat{q}_{m-3}^{\,-1}(A) \widehat{u}).
$$
Due to~\eqref{eq.exKryp2x1puhat},
this initial vector satisfies $ \widehat{q}_{m-3}^{\,-1}(A) \widehat{u} = u_q$
for $u_q$ given in~\eqref{eq.exKryp2x1purec}.
This implies
\begin{equation}\label{eq.exKryp2x1}
\RKry_{m-2}(A,\widehat{u} \hskip 1pt) = \Kry_{m-2}(A,u_q).
\end{equation}
To apply Proposition~\ref{prop.extendedKrylovalgorithmrmk}
for the rational Krylov subspace $\RKry_{m-1}(A,u)$ in~\eqref{eq.exKryp2x1purec},
we represent $\RKry_{m-1}(A,u)$
with poles $s_1,\ldots,s_{m-1}\in\C$
as an extended Krylov subspace of the form~\eqref{eq.xtendedKSinprop}.
Substituting $\widehat{u}$ for
the initial vector~$u$ in extended Krylov subspace in~\eqref{eq.xtendedKSinprop},
we have
\begin{equation}\label{eq.rewriteextorat01}
\RKry_{m-2}(A,\widehat{u}\hskip 1pt) \,\oplus\, \vspan\{A\,\widehat{u},A^2\widehat{u}\hskip 1pt \}.
\end{equation}
We proceed to show that this accumulated vector space
coincides with $\RKry_m(A,u)$.
Substituting
$$
\Kry_{m}(A,u_q)=\Kry_{m-2}(A,u_q) \oplus \vspan\{A^{m-2} u_q, A^{m-1} u_q \}
$$
for $\Kry_{m}(A,u_q)$ in~\eqref{eq.exKryp2x1purec}, we have 
\begin{equation}\label{eq.rewriteextorat02}
\RKry_m(A,u) = \Kry_{m-2}(A,u_q) \oplus \vspan\{A^{m-2} u_q, A^{m-1} u_q \}.
\end{equation}
Substituting $ \widehat{u} = \widehat{q}_{m-3}(A) u_q$,
we rewrite the right-hand term in~\eqref{eq.rewriteextorat01} to
$$
 \vspan\{A\,\widehat{u},A^2\widehat{u}\hskip 1pt \}
 = \vspan\{A\,\widehat{q}_{m-3}(A) u_q,A^2\widehat{q}_{m-3}(A) u_q\hskip 1pt \}.
$$
The matrix polynomials $A\,\widehat{q}_{m-3}(A)$
and $A^2\widehat{q}_{m-3}(A)$
correspond to polynomials of degree $m-2$ and $m-1$, respectively,
and this implies
$$
\Kry_{m-2}(A,u_q) \oplus \vspan\{A\,\widehat{u},A^2\widehat{u}\hskip 1pt \}
= \Kry_{m-2}(A,u_q) \oplus \vspan\{A^{m-2} u_q, A^{m-1} u_q \}.
$$
Combining~\eqref{eq.exKryp2x1} and~\eqref{eq.rewriteextorat02} with this identity, we conclude
$$
\RKry_m(A,u) = \RKry_{m-2}(A,\widehat{u}\hskip 1pt) \,\oplus\, \vspan\{A\,\widehat{u},A^2\widehat{u}\hskip 1pt \}.
$$

Thus, this rational Krylov subspace corresponds
to an extended Krylov subspace with initial vector $\widehat{u}=(A-s_{m-2}I)^{-1}(A-s_{m-1}I)^{-1} u$,
and the approach of of Proposition~\ref{prop.extendedKrylovalgorithmrmk} provides an algorithm
to compute a rational qor-Krylov representation $B_m$ of $\RKry_{m}(A,u)$ without accessing $q_{m-1}^{-1}(A) u$.
\end{remark}

Following Remark~\ref{rmk:BmforratKry}, the approach of Proposition~\ref{prop.extendedKrylovalgorithmrmk}
provides a procedure to compute the matrix $B_m$ for a rational Krylov subspace.
For the SaI Krylov subspace with a single pole $s\in\C$ of multiplicity $m-1$
and a fixed $\xi\in\R$ this is specified in Algorithm~\ref{alg.qorKryrat}.

\begin{figure}
\centering
\begin{algorithm}[H]
 \caption{An algorithm to compute the matrix $B_m$
for the SaI Krylov subspace with a single pole $s\in\C$ of multiplicity $m-1$
and a preassigned Ritz value $\xi\in\R$.
This algorithm follows Proposition~\ref{prop.extendedKrylovalgorithmrmk}
for a modified starting vector $\widehat{u}= X^2 u$ with $X=(A-s\,I)^{-1}$.}
\label{alg.qorKryrat}
\SetAlgoLined
$\widehat{u}= X^2 u$ with $X=(A-s\,I)^{-1}$\;
 run Algorithm~\ref{alg.SaI}
 to compute $U_{m-2}$ and $A_{m-2}$ for the SaI Krylov subspace $\Kry_{m-2}(X,\widehat{u} \hskip 1pt)$\;
 $\widetilde{v}=A\widehat{u}$\;
 orthogonalize $\widetilde{v}$ with $U_{m-2}$ and set $v_{m-1}=\widetilde{v}/\|\widetilde{v}\hskip 1pt\|_{\Minr}$\;
 $\widehat{v}=A\,v_{m-1}$\;
 for $j=1,\ldots,m-2$\;
 ~~~~  $y_j=(u_j,\widehat{v})_{\Minr}$\;
 ~~~~  $\widehat{v} \leftarrow \widehat{v} - y_j u_j$\;
 $a_{m-1} = (v_{m-1},\widehat{v}\,)_{\Minr}$ and $\widehat{v} \leftarrow \widehat{v} - a_{m-1} v_{m-1}$\;
 $\beta_{m-1}= \|\widehat{v}\hskip 1pt\|_{\Minr}$ and $v_m=\widehat{v}/\beta_{m-1}$\;
 $A_{m-1} = [\,A_{m-2}, y\,;~ y^\Hast, a_{m-1}\,] $\;
 $\omega_m=\xi + \beta_{m-1}^2 e_{m-1}^\Hast(A_{m-1} - \xi I)^{-1}e_{m-1}$\;
 $B_m  = [\,A_{m-1}, \beta_{m-1}e_{m-1}\,; ~\beta_{m-1}e_{m-1}^\Hast, \omega_m\,] $\;
 $U_m=(U_{m-2},v_{m-1},v_m)$\;
 $x=(U_m,u)_{\Minr}$\;
 return $x,U_m,B_m$\;
\end{algorithm}
\end{figure}

\paragraph{A rational qor-Krylov approximation
           to matrix functions~$ f(A)u $.}
$~$\newline
We refer to
\begin{equation}\label{eq.rqorfappr}
U_m f(B_m)x \approx f(A)u
\end{equation}
as a rational quasi-orthogonal residual (qor-)Krylov approximation.

\section[The CMS Theorem
for polynomial and some rational Krylov subspaces]
{The Separation Theorem of Chebyshev-Markov-Stieltjes (CMS Theorem)
for polynomial and some rational Krylov subspaces}
\label{sec.bounds3}
\medskip

The CMS Theorem
states that the accumulated quadrature weights of Gaussian quadrature formulae
are bounded by Riemann-Stieltjes integrals over the intervals
between the left integral limit and the quadrature nodes.
In Subsection~\ref{subsec.bounds31} we first reformulate previously stated identities
of the Krylov representation
(namely, Proposition~\ref{prop.inprodpol},~\ref{prop.ratkryeqinnerprod},~\ref{prop.eqforqorpolyKryinner} and~\ref{prop.eqforqorratKry})
as Gaussian quadrature formulae for the Riemann-Stieltjes integral associated
with the step function $\alpha_n$;
this allows us to present results in the following subsections
(which apply in the Krylov setting)
for a more general setting,~i.e.,~for Gaussian quadrature formulae.
We also recall some notation for Gaussian quadrature formulae
of Riemann-Stieltjes integrals, and we link classical notations
to the previously introduced setting.

In Subsection~\ref{subsec.bounds3polSepTh} we recapitulate
the CMS Theorem for the polynomial Krylov setting, and in
Subsection~\ref{subsec.SepTheoremRatKrylov}--\ref{subsec.STCMSxKry}
we introduce CMS type results for various rational
Krylov settings.

Throughout the
\ifthesis
present chapter,
\else
present work,
\fi
we consider integrals associated with a non-decreasing step function~$\alpha_n$ with~$n$ points of strict increase.
However, most of the results in the present section
hold true for integrals associated with non-decreasing continuous
functions $\alpha$ in a similar manner;
the case of $\alpha$ being a continuous
is not discussed in detail in the present work.

\subsection{Gaussian quadrature formulae and Krylov subspaces. Historical context}
\label{subsec.bounds31}

The integral associated with the step function $\alpha_n$
is to be understood as a Riemann-Stieltjes integral.
Gaussian quadrature formulae for Riemann-Stieltjes integrals are also referred to as
Gauss-Christoffel quadrature formulae in the literature, for previous remarks see also Subsection~\ref{subsec.introhistoric}.
For the Gauss-Christoffel quadrature formula which integrates polynomials of degree $\leq 2m-1$ exactly,
the quadrature nodes are given by the zeros of the associated orthogonal polynomial of degree $m$,
and the quadrature weights are given by so called Christoffel numbers.
Similar results hold for Gauss-Radau formulae for which 
the quadrature nodes and weights coincide with zeros of quasi-orthogonal polynomials and respective Christoffel numbers.
We briefly recapitulate the relation between Gaussian quadrature formulae
and the Jacobi matrix, which is also  mentioned in
Subsection~\ref{subsec.introhistoric};
for further details on Gaussian quadrature formulae we refer to~\cite{Ga81} and others.
Further below in the present subsection, we recall similar results for rational
Gaussian quadrature formulae.

\paragraph{The Christoffel numbers
and the eigendecomposition of the Jacobi matrix.}

For the orthonormal polynomials $p_0,\ldots,p_{m-1}$ associated with
the distribution $\dd\alpha_n$, see Proposition~\ref{prop.threetermorthopol},
we define
$$
\rho_{m-1}(\lambda) = 1 \Big/\sum_{k=0}^{m-1} p_k(\lambda)^2 \in\R.
$$
We recall that the Ritz values $\theta_1,\ldots,\theta_m\in\R$
correspond to the zeros of $p_m$.
The numbers $\rho_{m-1}(\theta_1),\ldots,\rho_{m-1}(\theta_m)$
are also referred to as {\em Christoffel numbers} in the literature.

We proceed to recall the relation between Christoffel numbers and
entries of eigenvectors of the Jacobi matrix which goes back to~\cite{Wi62,GW69}.
We introduce the denotation $c_1,\ldots,c_m\in\R$ for the spectral
coefficients of the vector $\beta_0\,e_1$ in
the eigenbasis of $J_m$,
which further correspond to the first components of the scaled eigenvectors:
Let $\widehat{q}_1,\ldots,\widehat{q}_m\in\R^m$ denote the
$\ell^2$-ortho\-normal eigenvectors of $J_m$,~i.e.,~$ J_m \widehat{q}_j = \theta_j \widehat{q}_j $
for the Ritz values $\theta_j$ and
$(\widehat{q}_j,\widehat{q}_k)_2=\delta_{jk}$, then
\begin{equation}\label{eq.cjintro}
c_j = \beta_0 \,(\widehat{q}_j,e_1)_2\in\R.
\end{equation}
The Christoffel numbers correspond to the first components
of the eigenvectors of the Jacobi matrix:
We recall the following results for the eigenvectors of $J_m$.
Following Section~\ref{sec.intro1},
the eigenvector for the eigenvalue $\theta_j$ is given by
\begin{equation}\label{eq.evJmbyp1m}
(p_0(\theta_j),\ldots,p_{m-1}(\theta_j))^\intercal \in\R^m.
\end{equation}
For the first component of the eigenvector we have $p_0=1/\beta_0$.
Thus, the first component of the $j$-th normalized eigenvector
scaled by $\beta_0$ and squared satisfies
\begin{equation}\label{eq.christoffelnumber}
c_j^2 = 1 \Big/\sum_{k=0}^{m-1} p_k(\theta_j)^2 \in\R,~~~j=1,\ldots,m,
\end{equation}
and for the Christoffel numbers we have the identity
\begin{equation}\label{eq.STintroandcj}
c_j^2=\rho_{m-1}(\theta_j),~~~j=1,\ldots,m.
\end{equation}
The Christoffel numbers are nonzero,\footnote{
       The result $c_j\neq 0$ is clarified in Appendix~\ref{appendix.A2},
       Proposition~\ref{prop.app3.cjneqzero}.
}~i.e.,~$c_j\neq 0$.
Although $c_j$ is real-valued, we also write $|c_j|^2$ in place of $c_j^2$.

\smallskip
Similar results hold for the spectrum of the qor-Krylov representation
$T_m$ introduced in Subsection~\ref{sec.polynomialqorKry}.
We reuse some notation associated with the spectrum of $J_m$ for $T_m$:
Corresponding to $T_m$ the denotations~$\theta_1,\ldots,\theta_m$
and~$c_1,\ldots,c_m$ refer to the eigenvalues of $T_m$ and
the spectral coefficients of $\beta_0\,e_1$ in the $\ell^2$-orthonormal
eigenbasis of $T_m$, respectively.
For the qor-Krylov representation $T_m$
we assume that the preassigned eigenvalue $\xi$ is given such that
the underlying quasi-orthogonal polynomial is well-defined,
and we assume that the eigenvalues of $T_m$ are included within
the integral limits of the respective Riemann-Stieltjes integral.
(See Proposition~\ref{prop.ritzqorpol} for some details on the location
of the eigenvalues of $T_m$.)
Following~\eqref{eq.qorpmatrixform}, the eigenvectors of $T_m$
conform to~\eqref{eq.evJmbyp1m} when $\theta_1,\ldots,\theta_m$
refer to the respective eigenvalues.
Similar to the case of the Jacobi matrix,
the representation~\eqref{eq.christoffelnumber}
and the identity~\eqref{eq.STintroandcj} also hold true for $T_m$.


\paragraph{A review on Gaussian quadrature formulae for the Riemann-Stieltjes integral.}
We proceed to reformulate Proposition~\ref{prop.inprodpol}
and~\ref{prop.ratkryeqinnerprod} as Gaussian quadrature formulae
for the Riemann-Stieltjes integral associated with the step function $\alpha_n$.
We recall that $\alpha_n$ is based on the eigenvalues of $A$ and the spectral coefficients of $u$.

For a complex-valued function $f\colon \R\to\C$,
where we consider polynomials or rational functions later on,
the following formulations are equivalent (see also~\eqref{eq.innerrseuclidean}),
\begin{subequations}\label{eq.innerprodrepresentation}
\begin{equation}\label{eq.innerprodrepresentationA}
\int_a^b f(\lambda)\dd \alpha_n(\lambda) = (u,f(A)u)_{\Minr}
= \sum_{j=1}^n f(\lambda_j) |w_j|^2.
\end{equation}
In a similar manner, the orthonormal eigendecomposition of $J_m$ yields
\begin{equation}\label{eq.innerprodrepresentationJm}
\beta_0^2 (e_1,f(J_m)e_1)_2  = \sum_{j=1}^m f(\theta_j) |c_j|^2.
\end{equation}
\end{subequations}
Identity~\eqref{eq.innerprodrepresentationJm} also holds true for $T_m$
if $\theta_j$ and $c_j$ refer to the spectrum of $T_m$.

The Ritz values $\theta_j$ and Christoffel numbers $\rho_{m-1}(\theta_j)$
provide quadrature nodes and weights, respectively, for the Gaussian quadrature
formulae which are also referred to as Gauss-Christoffel quadrature formulae
in the literature, see also~\cite{Ga81}.
We recapitulate classical results on Gaussian quadrature formulae
using the notation $|c_j|^2$ for the Christoffel numbers,
see~\eqref{eq.STintroandcj}.

\begin{remark}[Gaussian quadrature property, e.g.,~Subsection~6.2~\cite{GM10}]\label{rmk.gaussquad}
The Ritz values $\theta_1,\ldots,\theta_m$
and the spectral coefficients~$c_1,\ldots,c_m$~w.r.t.\ $J_m$
constitute a Gaussian quadrature formula for the Riemann-Stieltjes integral~\eqref{eq.innerrsintintro},
\begin{equation}\label{eq.gaussquadpolid}
\int_a^b p(\lambda)\dd \alpha_n(\lambda) = \sum_{j=1}^m p(\theta_j) |c_j|^2 ,~~~p\in \Pi_{2m-1}.
\end{equation}
Here, the Ritz values and the spectral coefficients
represent the quadrature nodes and quadrature weights, respectively.
On the basis of results of the present work,
identity~\eqref{eq.gaussquadpolid} can be verified via
the identities for the inner product in~\eqref{eq.defalphaninnerprod}
and~\eqref{eq.innerprodrepresentation},
$$
\int_a^b p(\lambda)\dd \alpha_n(\lambda)  = (u,p(A)u)_{\Minr} =
\beta_0^2(e_1,p(J_m)e_1)_2  = \sum_{j=1}^m p(\theta_j) |c_j|^2 ,~~~p\in \Pi_{2m-1}.
$$
Analogously, the qor-Krylov representation $T_m$ provides the following quadrature formula.
Let~$\theta_1,\ldots,\theta_m$ and~$c_1,\ldots,c_m$ be the eigenvalues and spectral coefficients of $T_m$,
then the identities~\eqref{eq.quasiKrylovidentityinner} for $p\in \Pi_{2m-2}$ together with~\eqref{eq.innerprodrepresentation} imply
\begin{equation}\label{eq.gaussquadpolid2}
\int_a^b p(\lambda)\dd \alpha_n(\lambda) = \sum_{j=1}^m p(\theta_j) |c_j|^2 ,~~~p\in \Pi_{2m-2}.
\end{equation}
When $\xi=a$ (thus, $\theta_1=a$) or $\xi=b$ (thus, $\theta_m=b$) is preassigned
this is also referred to as a Gauss-Radau quadrature formula.
\end{remark}

In view of Remark~\eqref{rmk.gaussquad} we summarize results for the Jacobi matrix $J_m$ and the qor-Krylov representation $T_m$.
For these results we write out the Riemann-Stieltjes integral~\eqref{eq.innerrsintintro}
in terms of its sum representation.
\begin{corollary}\label{cor.identityJmTm}
Let $\theta_1,\ldots,\theta_m$ and $c_1,\ldots,c_m$ denote the eigenvalues and spectral coefficients, respectively,
of either $J_m$ or $T_m$,
where the spectral coefficients $c_j$ refer to the vector $\beta_0\,e_1$.
Then,
\begin{equation}\label{eq.summarypolgaussquad}
\int_a^b p(\lambda)\dd \alpha_n(\lambda)
= \sum_{j=1}^m p(\theta_j) |c_j|^2,~~~p\in\Pi_{2m-2}.
\end{equation}
\end{corollary}

\paragraph{Rational Gaussian quadrature formulae and rational Krylov subspaces.}

~ For  
rational Krylov subspaces~$\RKry_m(A,u)$ we recall the definition
of the Rayleigh quotient~$A_m=(U_m, AU_m)_{\Minr}$,
where $U_m$ is an orthonormal basis of~$\RKry_m(A,u)$.
Furthermore, the vector $x=(U_m,u)_{\Minr}$
and the rational qor-Krylov representation $B_m$
(introduced in Subsection~\ref{sec.ratqor} via~\eqref{eq.defBm})
implicitly depend on $U_m$.
In the sequel we consider $U_m$ to be fixed,
and we assume that $B_m$ is well-defined.
For the latter we refer to the conditions concerning the definition
of $T_m$ in Section~\ref{sec.2qorpol}.
We proceed to reuse the denotation~$\theta_1,\ldots,\theta_m$
for the eigenvalues of~$A_m$ (`rational' Ritz values),
and $c_1,\ldots,c_m$ for the spectral coefficients of~$x$
in the orthonormal eigenbasis of~$A_m$:
Let $\widehat{q}_j\in\C^m$
denote the $\ell^2$-orthonormal eigenvectors
of $A_m$,~i.e.,~$A_m \widehat{q}_j = \theta_j \widehat{q}_j$ and
$(\widehat{q}_j ,\widehat{q}_k)_2=\delta_{jk}$, then
\begin{equation}\label{eq.cjforAm}
c_j=(\widehat{q}_j,x)_2\in\C,~~~j=1,\ldots,m.
\end{equation}
We remark that the coefficients $|c_j|$ are independent of
the explicit choice of the orthonormal basis $U_m$,
this is clarified in Proposition~\ref{prop.cjindependentofUm},
Appendix~\ref{appendix.A2}.
For a function $f\colon\R\to\C$, the eigendecomposition of $A_m$
yields
\begin{equation}\label{eq.innerprodrepresentationAm}
(x,f(A_m)x)_2  = \sum_{j=1}^m f(\theta_j) |c_j|^2.
\end{equation}

In the context of the rational qor-Krylov representation $B_m$ the denotation
$\theta_1,\ldots,\theta_m$ and $c_1,\ldots,c_m$ is reused accordingly,
and an identity similar to~\eqref{eq.innerprodrepresentationAm}
holds true for $B_m$ when $\theta_j$ and $c_j$ refer to the spectrum of $B_m$.

\begin{remark}\label{rmk.ratgaussquad}
Similar to Remark~\ref{rmk.gaussquad}, the identity
in~\eqref{eq.ratkryeqinnerprod}
corresponds to the following rational Gaussian quadrature formula.
Let $\theta_1,\ldots,\theta_m$ and $c_1,\ldots,c_m$ refer to
the spectrum of $A_m$, then
\begin{equation}\label{eq.gaussquadratid}
\int_a^b r(\lambda)\,\dd \alpha_n(\lambda)
= \sum_{j=1}^m r(\theta_j) |c_j|^2,~~~r\in\Pi_{2m-1}/|q_{m-1}|^2.
\end{equation}
To demonstrate~\eqref{eq.gaussquadratid} we recall the identities for the inner product
in~\eqref{eq.ratkryeqinnerprod},~\eqref{eq.innerprodrepresentationA},
and~\eqref{eq.innerprodrepresentationAm},
$$
\int_a^b r(\lambda)\dd \alpha_n(\lambda) = (u,r(A)u)_{\Minr}
= (x,r(A_m)x)_2 = \sum_{j=1}^m |c_j|^2 r(\theta_j),~~~r\in\Pi_{2m-1}/|q_{m-1}|^2.
$$
The rational qor-Krylov representation $B_m$
provides the following quadrature formula via Proposition~\ref{prop.eqforqorratKry},
$$
\int_a^b r(\lambda)\dd \alpha_n(\lambda)
= \sum_{j=1}^m r(\theta_j) |c_j|^2,~~~r\in\Pi_{2m-2}/|q_{m-1}|^2.
$$
When the preassigned eigenvalue of $B_m$ is set to one of the
integral limits, i.e.,~$\theta_1=a$ or $\theta_m=b$,
then this formula is also referred to as
rational Gauss-Radau quadrature formula.
\end{remark}

We summarize the statements of Remark~\ref{rmk.ratgaussquad}
concerning $A_m$ and $B_m$.
\begin{corollary}\label{cor.identityAmBm}
Let $\theta_1,\ldots,\theta_m$ and $c_1,\ldots,c_m$ denote the eigenvalues and spectral coefficients, respectively,
of either $A_m$ or $B_m$,
where the spectral coefficients refer to the vector~$x$.
Then,
\begin{equation}\label{eq.summaryratgaussquad}
\int_a^b r(\lambda)\,\dd \alpha_n(\lambda)
= \sum_{j=1}^m r(\theta_j) |c_j|^2 ,~~~r\in\Pi_{2m-2}/|q_{m-1}|^2.
\end{equation}
\end{corollary}

\subsection{The CMS Theorem for the polynomial case}
\label{subsec.bounds3polSepTh}
The CMS Theorem dates back to works of Chebyshev, Markov
and Stieltjes in the 19th century 
and also goes by the name Chebyshev-Markov-Stieltjes inequalities.
For further historical and technical remarks
we refer to~\cite[Section~3.41]{Sze85}
(including an extensive survey of this
theorem),~\cite[Theorem~2.54]{Ak65},~\cite[Section~4]{vA93},~\cite[Section~3]{LS13},~\cite{Chi78}
and others.

The Riemann-Stieltjes integral
associated with $\alpha_n$~\eqref{eq.measureAv}
over a subset of $(a,b)$ can be understood as a measure of such a subset.
Namely, with $\alpha_n(a)=0$ we consider $\alpha_n(\theta)$
to be the associated measure of the interval $(a,\theta]$ for $\theta\in(a,b)$.
To simplify the notation in the sequel,
we let $\mu_n(R)$ denote the measure of a subset $R$ of $(a,b)$ associated with $\alpha_n$.
More precisely, we first define
\begin{subequations}\label{eq.defmun0}
\begin{equation}\label{eq.defJR0x}
J(R)=\big\{j:~\lambda_j\in R \big\}\subset \{1,\ldots,n\},~~~~\text{for a set $R\subset (a,b)$}.
\end{equation}
The sum of the spectral coefficients $w_j$ over the index set $J(R)$
corresponds to the measure of the set $R$ associated with $\alpha_n$, and we define
\begin{equation}\label{eq.defmun}
\mu_n(R) = \sum_{j \in J(R)} |w_j|^2.
\end{equation}
\end{subequations}
Thus, we have $\mu_n((a,\theta]) = \alpha_n(\theta)$ for $\theta\in(a,b)$.
Furthermore, we proceed to use the notation $\mu_n$ and $\alpha_n$
for the measure of an interval $(a,\theta]$ in an equivalent manner.
Similarly, we use the notation $\alpha_n(\theta-)$
for the measure of the open interval~$(a,\theta)$,~i.e.,
$$
\alpha_n(\theta-) := \lim_{\varepsilon\to 0+} \alpha_n(\theta-\varepsilon) =\mu_n((a,\theta)).
$$


\smallskip
We proceed to recall the CMS Theorem.
This theorem is based on the Gaussian quadrature properties~\eqref{eq.summarypolgaussquad}
as in Corollary~\ref{cor.identityJmTm},
and thus, the following results hold true when $\theta_j$ and $c_j$ refer to
the spectrum of the Jacobi matrix $J_m$
or the qor-Krylov representation $T_m$.

\begin{theorem}[Separation Theorem of Chebyshev-Markov-Stieltjes,
see also Section 3.41 in~\cite{Sze85}]\label{thm.STCMS}
\label{Kryproponweights}
Let~$\theta_1,\ldots,\theta_m\in(a,b)$ and~$c_1,\ldots,c_m\in\C$ satisfy the Gaussian quadrature property~\eqref{eq.summarypolgaussquad},
then
\begin{equation}\label{eq.cjwjres1}
\alpha_n(\theta_k)<|c_1|^2+\ldots+|c_k|^2<\alpha_n(\theta_{k+1}-),~~~~k=1,\dots,m-1.
\end{equation}
\end{theorem}

We point out that for~$k=m$ the bounds in~\eqref{eq.cjwjres1} can be replaced by the following identity.
The Gaussian quadrature property~\eqref{eq.summarypolgaussquad} for $p=1$ implies
\begin{equation}\label{eq.coeftrivial1tom}
\sum_{j=1}^m |c_j|^2 = \alpha_n(b).
\end{equation}
(This also results directly from $\|u\|_{\Minr}=\beta_0 \|e_1\|_2$.)

To recall a classical proof of the CMS Theorem we introduce the following polynomials.
\begin{proposition}[Eq.~(3.411.1) in~\cite{Sze85}, part of Theorem (2.5.4)
in~\cite{Ak65} and others\footnote{
     A classical proof of Proposition~\ref{prop.thatpol1}
     is recapitulated in Appendix~\ref{appendix.A4}.
}]\label{prop.thatpol1}
Let~$\theta_1<\ldots<\theta_m\in\R$ and let~$k$ be fixed
with~$1\leq k<m$.
Then there exist polynomials~$p_{\{+,k\}}$ and~$p_{\{-,k\}}\in \Pi_{2m-2}$
which satisfy\footnote{
     In the sequel statements concerning $p_{\{\pm,k\}}$
     apply to $p_{\{+,k\}}$ and $p_{\{-,k\}}$ individually.}
\begin{equation}\label{eq.ppmequal}
p_{\{\pm,k\}}(\theta_j)=\left\{
\begin{array}{ll}
1,&j=1,\ldots,k,\\
0,&j=k+1,\ldots,m,
\end{array}
\right.
\end{equation}
together with
\begin{equation}\label{eq.ppminequal}
p_{\{+,k\}}(\lambda) \geq \left\{
\begin{array}{ll}
1,&\lambda \leq \theta_k,\\
0,& \lambda > \theta_k,
\end{array}
\right.
~~~\text{and}~~~
p_{\{-,k\}}(\lambda) \leq \left\{
\begin{array}{ll}
1,& \lambda < \theta_{k+1},\\
0,& \lambda \geq \theta_{k+1}.
\end{array}
\right.
\end{equation}
Additionally, the inequalities in~\eqref{eq.ppminequal} are strict inequalities for $\lambda\notin\{\theta_1,\ldots,\theta_m\}$.
\end{proposition}

\begin{figure}
\centering
\begin{tabular}{cc}
\begin{overpic}
[width=0.5\textwidth]{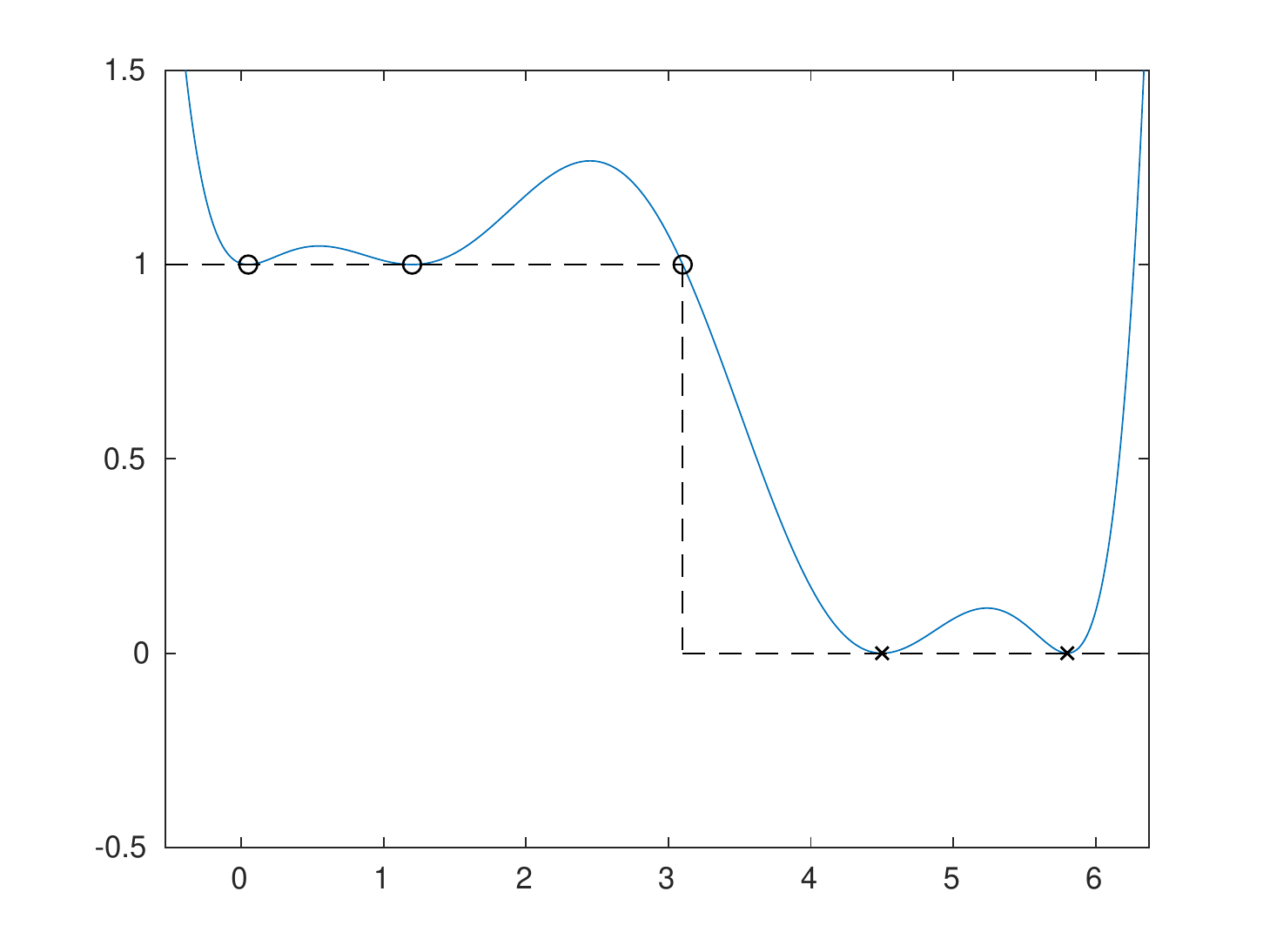}
\put(50,1){\small~$\lambda$}
\put(0,30){\rotatebox{90}{\small~$p_{\{+,k\}}(\lambda)$}}
\put(52.3,56.5){\small~$p_{\{+,k\}}(\theta_k)$}
\end{overpic}
&
\hspace{-0.5cm}
\begin{overpic}
[width=0.5\textwidth]{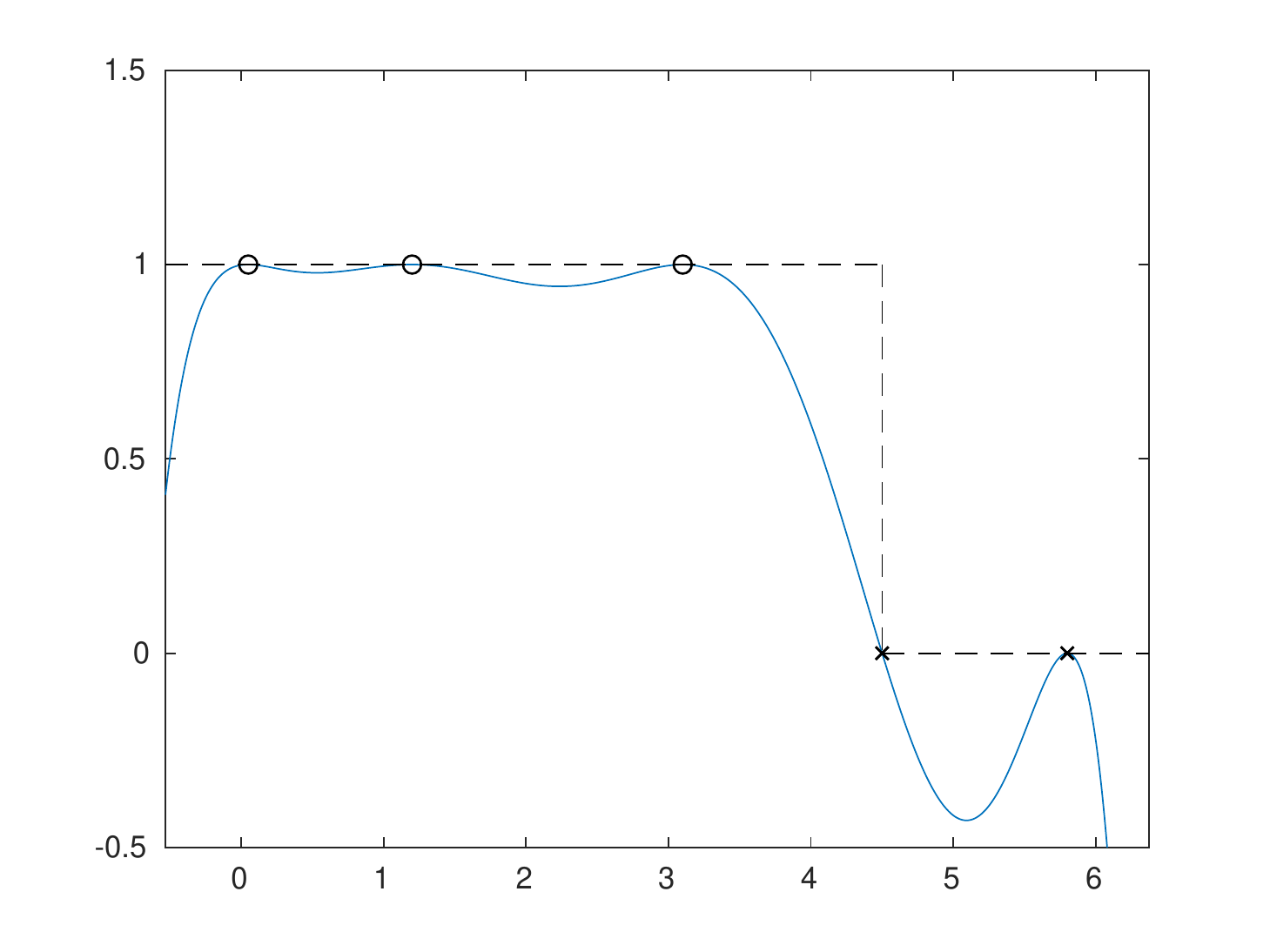}
\put(50,1){\small~$\lambda$}
\put(1,30){\rotatebox{90}{\small~$p_{\{-,k\}}(\lambda)$}}
\put(52,57){\small~$p_{\{-,k\}}(\theta_k)$}
\end{overpic}
\end{tabular}
\caption{This figure illustrates the polynomials~$p_{\{+,k\}}$ (left) and~$p_{\{-,k\}}$ (right) given in Proposition~\ref{prop.thatpol1}
for given nodes $\theta_1,\ldots,\theta_m$ with $m=5$.
The identities~\eqref{eq.ppmequal} are illustrated for~$\theta_1,\ldots,\theta_k$ ('$\circ$')
and~$\theta_{k+1},\ldots,\theta_m$ ('$\times$') with~$k=3$,
and the dashed line illustrates the bounds~\eqref{eq.ppminequal}.
}
\label{fig:pplusminusleft}
\end{figure}

The polynomials of Proposition~\ref{prop.thatpol1}
are illustrated in \figref{fig:pplusminusleft}
for a numerical example.

\smallskip
With Proposition~\ref{prop.thatpol1} we proceed to prove Theorem~\ref{Kryproponweights}.
\begin{proof}[Proof of \,\textbf{Theorem~\ref{Kryproponweights}}]
Let~$p_{\{\pm,k\}}\in \Pi_{2m-2}$ be given according to Proposition~\ref{prop.thatpol1}
for the eigenvalues~$\theta_1<\ldots<\theta_m$ and~$k=1,\ldots,m-1$.
The polynomials~$p_{\{\pm,k\}}$ satisfy~\eqref{eq.ppmequal},
and this implies
\begin{equation}\label{eq.prc1.lb1}
\sum_{j=1}^m p_{\{\pm,k\}}(\theta_j) |c_j|^2
= \sum_{j=1}^k |c_j|^2.
\end{equation}
On the other hand, identity~\eqref{eq.summarypolgaussquad} yields
$$
\sum_{j=1}^m p_{\{\pm,k\}}(\theta_j) |c_j|^2
= \int_a^b p_{\{\pm,k\}}(\lambda)\,\dd \alpha_n(\lambda).
$$
Evaluating the Riemann-Stieltjes integral in this identity, we arrive at
\begin{equation}\label{eq.prc1.lb2}
\sum_{j=1}^m p_{\{\pm,k\}}(\theta_j) |c_j|^2
= \sum_{j=1}^n p_{\{\pm,k\}}(\lambda_j) |w_j|^2.
\end{equation}

Let the index set $J((a,\theta_k])\subset\{1,\ldots,n\}$ be given as in~\eqref{eq.defJR0x}.
Then, the inequalities for $p_{\{+,k\}}$ in~\eqref{eq.ppminequal} imply
\begin{equation}\label{eq.prc1.lb3}
\sum_{j=1}^n p_{\{+,k\}}(\lambda_j) |w_j|^2
\, > \sum_{j\in J((a,\theta_k])} |w_j|^2
= \alpha_n(\theta_k).
\end{equation}
This inequality is strict due to the interlacing property of
the eigenvalues $\lambda_j$ and $\theta_j$, see Proposition~\ref{prop.ritzvaluespol}
and~\ref{prop.ritzqorpol}.
Combining~\eqref{eq.prc1.lb1},~\eqref{eq.prc1.lb2} and~\eqref{eq.prc1.lb3}
yields the lower bound in~\eqref{eq.cjwjres1}.

\smallskip
Similarly to~\eqref{eq.prc1.lb3},
the inequalities for~$p_{\{-,k\}}$ in~\eqref{eq.ppminequal} imply
\begin{equation}\label{eq.prc1.lb4}
\sum_{j=1}^n p_{\{-,k\}}(\lambda_j) |w_j|^2
\, < \sum_{ j\in J((a,\theta_{k+1}))}|w_j|^2
= \alpha_n(\theta_{k+1}-).
\end{equation}
Combining~\eqref{eq.prc1.lb1},~\eqref{eq.prc1.lb2} and~\eqref{eq.prc1.lb4}
yields the upper bound in~\eqref{eq.cjwjres1}.
\end{proof}

The inequalities~\eqref{eq.cjwjres1} in Theorem~\ref{thm.STCMS} yield
the following bounds on the measure of the intervals located between Ritz values.
In the following, we use the notation $\mu_n$ for the measure as in~\eqref{eq.defmun}.
\begin{corollary}\label{cor.piecewiseboundspol}
In the setting of Theorem~\ref{thm.STCMS}, the following inequalities hold true.
\begin{itemize}
\item For indices $j,k$ with $1<j<k<m$,
\begin{subequations}\label{eq.septhmpolcorboundsin}
\begin{equation}\label{eq.septhmpolcorboundsina}
\mu_n([\theta_j,\theta_k])
< |c_j|^2+|c_{j+1}|^2+\ldots+|c_k|^2
< \mu_n((\theta_{j-1},\theta_{k+1}))
\end{equation}

\item Furthermore, the 
accumulated spectral coefficients satisfy
\begin{equation}\label{eq.septhmpolcorboundsinb}
\mu_n([\theta_j,b))
< |c_j|^2+\ldots+|c_m|^2
< \mu_n((\theta_{j-1},b)),
~~~j=2,\ldots,m.
\end{equation}
\end{subequations}
\end{itemize}
\end{corollary}

\begin{proof}
Applying~\eqref{eq.cjwjres1} twice
(once we substitute $j-1$ for the index $k$ therein)
and subtracting, we observe
\begin{equation*}
\alpha_n(\theta_k) - \alpha_n(\theta_j-)
< |c_j|^2+|c_{j+1}|^2+\ldots+|c_k|^2
< \alpha_n(\theta_{k+1}-) - \alpha_n(\theta_{j-1}),
\end{equation*}
this shows~\eqref{eq.septhmpolcorboundsina}.
Subtracting~\eqref{eq.cjwjres1} for the index $j-1$
from~\eqref{eq.coeftrivial1tom}, we arrive~at
\begin{equation}
\alpha_n(b) - \alpha_n(\theta_j-) 
< |c_j|^2+\ldots+|c_m|^2
< \alpha_n(b) - \alpha_n(\theta_{j-1}),
\end{equation}
which entails~\eqref{eq.septhmpolcorboundsinb}.
\end{proof}


We proceed to specify the intertwining property of the distributions $\dd\alpha_n$
and $\dd\alpha_m$ which already appeared in the introduction
of Subsection~\ref{subsec.introhistoric}:
Similarly to~$\alpha_n$ in~\eqref{eq.measureAv},
we introduce the step function
\begin{equation}\label{eq.alpham}
\alpha_m(\lambda)=\left\{
\begin{array}{ll}
0,&\lambda<\theta_1,\\
\sum_{j=1}^\ell |c_j|^2, & \theta_\ell \leq \lambda < \theta_{\ell+1},~~~\ell=1,\ldots,m-1,\\
\sum_{j=1}^m |c_j|^2, & \theta_m \leq \lambda.
\end{array}
\right.
\end{equation}
For $f:\R\to\C$
the Riemann-Stieltjes integral associated with $\alpha_m$ reads
$$
\int_a^b f(\lambda)\,\dd\alpha_m(\lambda) = \sum_{j=1}^m |c_j|^2 f(\theta_j).
$$
Thus, the quadrature property in Corollary~\ref{cor.identityJmTm} coincides with the identity
\begin{equation}\label{eq.polmoments}
\int_a^b \lambda^j \,\dd\alpha_n(\lambda) = \int_a^b \lambda^j \,\dd\alpha_m(\lambda),~~~j=0,\ldots,2m-2.
\end{equation}
The integral terms in~\eqref{eq.polmoments} correspond to the moments of
the distributions $\dd\alpha_n$ and~$\dd\alpha_m$,
and thus, the Gaussian quadrature property in Corollary~\ref{cor.identityJmTm}
coincides with $\dd\alpha_n$ and $\dd\alpha_m$
having matching moments up to order $2m-2$.
We define the auxiliary function
\begin{equation} \label{eq:Fxalpha}
F(\lambda)=\alpha_n(\lambda)-\alpha_m(\lambda),
\end{equation}
and remark the following properties of $F$.
The step functions $\alpha_n(\lambda)$ and $\alpha_m(\lambda)$ are both increasing in $\lambda$,
whereat the step function $\alpha_m(\lambda)$ has
exactly $m$ points of increase at $\lambda=\theta_1,\ldots,\theta_m$.
Thus, the function $F(\lambda)$ is
increasing for $\lambda\in(\theta_k,\theta_{k+1})$, $k=1,\ldots,m-1$,
and away from the boundaries $\lambda<\theta_1$ and $\lambda>\theta_m$.
Furthermore, Theorem~\ref{Kryproponweights} yields
$$
\alpha_n(\theta_k)- (|c_1|^2+\ldots+|c_k|^2)
< 0
< \alpha_n(\theta_{k+1}-) - (|c_1|^2+\ldots+|c_k|^2),~~~~k=1,\dots,m-1.
$$
The accumulated coefficients $c_k$ correspond to the step function $\alpha_m$~\eqref{eq.alpham},
namely,
\begin{equation}\label{eq.alpmthetakpm}
|c_1|^2+\ldots+|c_k|^2 = \alpha_m(\theta_k) = \alpha_m(\theta_{k+1}-),
\end{equation}
and we observe the inequalities
\begin{equation}\label{eq.Flowerbound}
 F(\theta_k) < 0 < F(\theta_{k+1}-)~~~\text{for}~~k=1,\ldots,m-1.
\end{equation}
More precisely, the inequalities~\eqref{eq.Flowerbound} are
equivalent to the assertion of the
CMS Theorem (Theorem~\ref{Kryproponweights}).

To clarify the intertwining property of $\dd\alpha_n$ and $\dd\alpha_m$
in this context:
The CMS Theorem relies on quadrature properties
which correspond to~\eqref{eq.polmoments},~i.e.,~$\dd\alpha_n$ and $\dd\alpha_m$
having matching moments,
and the result of the CMS Theorem corresponds
to~\eqref{eq.Flowerbound}, which can be understood as
an intertwining property of $\dd\alpha_n$ and $\dd\alpha_m$.

Besides these remarks, the function $F$ is further used
in the following subsection to rewrite CMS type results
for rational cases, and in Section~\ref{sec.experiments4} below
where we verify results of the present section
for numerical examples.

\begin{remark}\label{rmk.boundsforweightederrornorm}
In the
\ifthesis
present chapter,
\else
present work,
\fi
the measure $\alpha_n$ is introduced
based on eigenvalues $\lambda_j$ of $A$
and the spectral coefficients $w_j$
of the initial vector $u$ in the eigenbasis of $A$ as in~\eqref{eq.measureAv}.
Thus, the bounds given by
the CMS Theorem
reveal bounds for the accumulated spectral coefficients $w_j$.
To simplify the notation we proceed with the setting of the
Jacobi matrix $J_m$,~i.e.,~the eigenvalues $\theta_j$ and
spectral coefficients $c_j$ refer to the spectrum of the Jacobi matrix.
In a similar manner such results also hold for the
qor-Krylov representation $T_m$ as specified below.
The bounds on $\alpha_n$ provided by the CMS Theorem
are computable,~i.e.,~$\theta_j$ and $c_j$ are available via an eigendecomposition
of the Jacobi matrix which can be computed using the Lanczos method.

\smallskip
We proceed in the setting of the Jacobi matrix.
For its eigenvalues $\theta_j$ we define
the index~$\ell=\ell(k)$ for~$k=1,\ldots,m$, such that
\begin{equation}\label{eq.deflk}
 \lambda_{\ell(k)} \leq \theta_k < \lambda_{\ell(k)+1}.
\end{equation}
The positioning of the eigenvalues,
which is specified in Proposition~\ref{prop.ritzvaluespol},
implies~$\ell(k)<\ell(k+1)$ for~$k=1,\ldots,m-1$
and~$1\leq \ell(k)< n$ for~$k=1,\ldots,m$.

With $\ell(k)$ defined in~\eqref{eq.deflk} we have the representations
\begin{equation}\label{eq.alphantktowjlk}
\alpha_n(\theta_k) = \sum_{j=1}^{\ell(k)} |w_j|^2 ,~~~~\text{and}
~~~\alpha_n(\theta_{k+1}) = \sum_{j=1}^{\ell(k+1)} |w_j|^2,~~~k=1,\ldots,m-1.
\end{equation}
Note that $\alpha_n(\theta_{k+1}-) \leq \alpha_n(\theta_{k+1})$;
to keep the notation simple, the case
$ \alpha_n(\theta_{k+1}-) <\alpha_n(\theta_{k+1}) $ is not treated separately here.
For the remainder of the present remark we assume 
$$
\lambda_{\ell(k)}\neq \theta_k,~~~k=1,\ldots,m.
$$

Thus, with~\eqref{eq.alphantktowjlk} Theorem~\ref{Kryproponweights} reads
\begin{equation}\label{eq.cjwjres1ell}
\sum_{j=1}^{\ell(k)} |w_j|^2 < \sum_{j=1}^k |c_j|^2 < \sum_{j=1}^{\ell(k+1)} |w_j|^2,~~~~k=1,\dots,m-1.
\end{equation}

Furthermore, for a set of eigenvalues of $A$ located between 
two Ritz values $\theta_j$ and $\theta_k$ with $j<k$ we recall
$$
\lambda_{\ell(j)} < \theta_{j} < \lambda_{\ell(j)+1} < \ldots < \lambda_{\ell(k)} < \theta_k
,~~~k=2,\ldots,m,
$$
and with~\eqref{eq.alphantktowjlk}, the sum of spectral coefficients $w_j$ associated
with these eigenvalues corresponds to
\begin{equation}\label{eq.pkboundwjpiecewise}
\sum_{\iota=\ell(j)+1}^{\ell(k)} |w_\iota|^2
= \alpha_n(\theta_k) - \alpha_n(\theta_{j}), ~~~j<k.
\end{equation}
Furthermore, combining this identity with~\eqref{eq.cjwjres1ell} or~\eqref{eq.septhmpolcorboundsin},
we obtain computable bounds on accumulated spectral coefficients of $u$.
E.g., for $1<j<k<m$ the inequality~\eqref{eq.septhmpolcorboundsina} yields
$$
|c_{j+1}|^2+\ldots+|c_{k-1}|^2
< \sum_{\iota=\ell(j)+1}^{\ell(k)} |w_\iota|^2
< |c_j|^2+\ldots+|c_k|^2,
$$
where the lower bound is trivial in the case $k=j+1$.

We remark that the results of the present subsection can be generalized
to the setting of the qor-Krylov representation $T_m$.
For the qor-Krylov representation,
the cases $\theta_1<\lambda_1$ and $\lambda_n<\theta_m$
have to be considered explicitly in the notation,~namely,~the indices
$\ell(1)$ and $\ell(m)$ have to be adapted accordingly for these cases.
\end{remark}

\begin{remark}\label{rmk:STwithsj}
In the present work the measure $\alpha_n$
is based on the spectrum of $A$ and has $n$ points of strict increase.
Thus, the identity of~\cite[eq. (3.41.3)]{Sze85} which relies on a
continuous measure does not hold true in the present case,~i.e.,
$$
\text{in general, we do not find any point $y_k\in\R$ such that}
~~ \alpha_n(y_k) = \sum_{j=1}^k |c_j|^2.
$$
Nevertheless, the inequalities in~\eqref{eq.cjwjres1ell}
imply that there exist indices $\nu_k$
with $\ell(k)<\nu_k\leq\ell(k+1)$
and numbers $\xi_k\in(0,1]$ for $k=1,\ldots,m-1$ such that
$$
\sum_{j=1}^{\nu_{k}-1} |w_j|^2 + \xi_{k} |w_{\nu_k}|^2 = \sum_{j=1}^k |c_j|^2,
$$
This can give further theoretical insight on the estimates provided
in Remark~\ref{rmk.boundsforweightederrornorm}.
Nevertheless, the indices $\nu_j$ and scaling factors $\xi_j$
are not computable in general.

The indices $\nu_k$ satisfy $\lambda_{\nu_k}\in(\theta_k,\theta_{k+1}]$, thus,
$$
\lambda_1 < \theta_1 < \lambda_{\nu_1} \leq \theta_2 < \lambda_{\nu_2}
< \ldots \leq \theta_{m-1} < \lambda_{\nu_{m-1}} \leq \theta_m < \lambda_n.
$$
For each spectral coefficient $c_k$, this implies
$$
\begin{aligned}
&|c_1|^2 =  \sum_{j=1}^{\nu_{1}-1} |w_j|^2 + \xi_{1} |w_{\nu_1}|^2,\\
&|c_k|^2 =  (1-\xi_{k-1}) |w_{\nu_{k-1}}|^2 + \sum_{j=\nu_{k-1}+1}^{\nu_{k}-1} |w_j|^2 + \xi_{k} |w_{\nu_k}|^2,
~~~k=2,\ldots,m-1,~~~\text{and}\\
&|c_m|^2 =  (1-\xi_{m-1}) |w_{\nu_{m-1}}|^2 + \sum_{j=\nu_{m-1}+1}^{n} |w_j|^2.
\end{aligned}
$$
\end{remark}

%
\subsection{The rational case with a single pole \texorpdfstring{$s\in\R$}{s in R} of higher multiplicity}
\label{subsec.SepTheoremRatKrylov}

In the present subsection we consider CMS type results
for the setting of a rational Krylov subspace $\RKry_m(A,u)$
with a single pole $s\in\R$, thus, we have the denominator
$q_{m-1}(\lambda)=(\lambda-s)^{m-1}$.
This subspace corresponds to a SaI Krylov subspace; for previous
remarks see also Subsection~\ref{subsec.remarksSaI}.
Following Subsection~\ref{subsec.bounds31},
the eigenvalues $\theta_1,\ldots,\theta_m\in (a,b)$ 
and spectral coefficients $c_1,\ldots,c_m\in\C$ of the respective
Rayleigh quotient $A_m$ or qor-representation $B_m$ satisfy the
quadrature property~\eqref{eq.summaryratgaussquad} in Corollary~\ref{cor.identityAmBm}.
To provide results in a more general setting, the results in the remainder
of the subsection are based on the quadrature property~\eqref{eq.summaryratgaussquad};
we provide results for a class of rational Gaussian quadrature formulae
which fit to the respective SaI Krylov setting.

Although the rational Krylov subspace corresponds to the polynomial Krylov subspace
$\Kry_m(A,u_q)$ with starting vector $u_q=q_{m-1}^{-1}(A)u$,
results of the previous subsection do not yield bounds associated with $\alpha_n$,
this is specified in the following remark.
\begin{remark}
The rational Krylov subspace~$\RKry_m(A,u)$
with the respective denominator $q_{m-1}$
is identical to~$\Kry_m(A,u_q)$ with $u_q=q_{m-1}^{-1}(A)u$.
This polynomial Krylov subspace is associated with the
step function $\widehat{\alpha}_n$ given in~\eqref{eq.measureAvrat}.
Let $J_m$ and $V_m$ denote the Jacobi matrix and
the \tMinr-orthonormal eigenbasis of~$\Kry_m(A,u_q)$
constructed by the Lanczos method.
In the setting of~$\Kry_m(A,u_q)$, Theorem~\ref{Kryproponweights}
yields bounds based on spectral coefficients of the vector
$x_q = (V_m, u_q)_{\Minr}$ in the eigenbasis of $J_m$
and the step function $\widehat{\alpha}_n$.
This does not entail
bounds based on spectral coefficients of $x = (V_m, u)_{\Minr}$
and the step function $\alpha_n$ in general.
\end{remark}

\begin{figure}
\centering
\begin{overpic}
[width=0.5\textwidth]{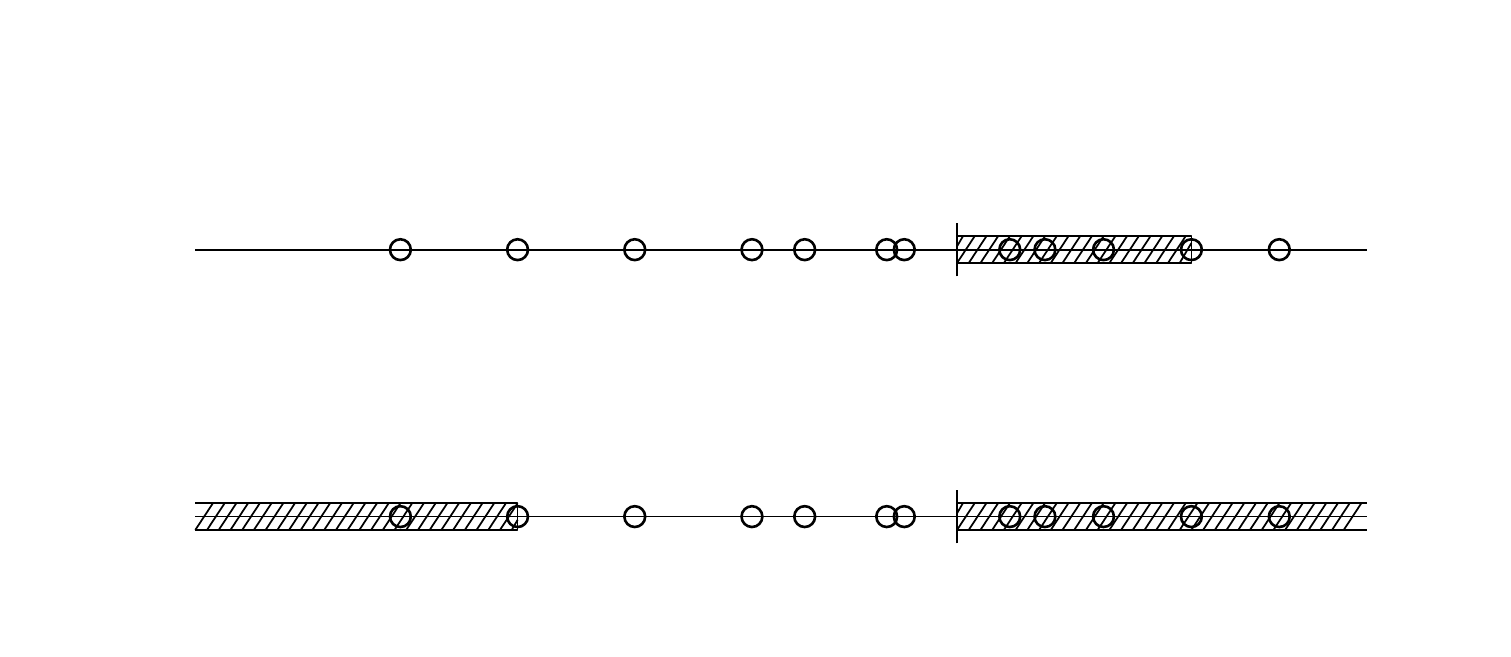}
\put(-30,32){\small \textbf{a)}}
\put(-20,32){\small {the case $\theta_k>s$:}}
\put(62.5,30){\small$ s $}   
\put(78,30){\small$\theta_k$}
\put(95,26){\small$R_k=(s,\theta_k]$}
\put(-30,14){\small \textbf{b)}}
\put(-20,14){\small {the case $\theta_k<s$:}}
\put(62.5,12){\small$ s $}   
\put(33,12){\small$\theta_k$}
\put(95,8){\small$R_k=(a,\theta_k]\cup (s,b)$}
\put(80,4){\vector(1,0){10}}
\put(82,1){\small$ +\infty $}
\put(23,4){\vector(-1,0){10}}
\put(12,1){\small$ -\infty $}
\end{overpic}
\caption{In Figure~\textbf{a)} and~\textbf{b)} we illustrate the
set $R_k\subset\R$ given in~\eqref{eq.thatr2defRk}
for a given sequence of nodes $\theta_1,\ldots,\theta_m$ ('$\circ$'),
and a given pole $s$ which satisfies $\theta_1<s<\theta_m$.
In Figure~a) we choose the index $k$ such that $\theta_k>s$,
and in Figure~b) we consider~$\theta_k<s$.
In each figure the set $R_k$ is highlighted by a dashed area.
}
\label{fig:thatrproRk0}
\end{figure}

To simplify the notation in the sequel,
we first define the indices $k_1$ and $k_m$ such that
\begin{subequations}\label{eq.thatr2defk1m}
\begin{equation}\label{eq.thatr2defk1ma}
\theta_{k_m} < s < \theta_{k_1},~~~k_1=k_m+1,~~
\text{in case of $s\in(\theta_1,\theta_m)$},
\end{equation}
and otherwise,
\begin{equation}\label{eq.thatr2defk1mb}
k_1=1~~~\text{and}~~k_m=m,~~~\text{in case of $s<\theta_1$ or $s>\theta_m$}.
\end{equation}
\end{subequations}
Furthermore, we define the sets~$I_k\subset\{1,\ldots,m\}$
and~$R_k\subset \R$ for~$k=1,\ldots,m$ by
\begin{equation}\label{eq.thatr2defRk}
\begin{aligned}
&I_k=\left\{
\begin{array}{ll}
\{k_1,\ldots,k \}, & ~~ k_1\leq k \leq m,\\
\{1,\ldots,k,k_1,\ldots,m \}, & ~~ 1 \leq k < k_1,
\end{array}\right.~~~\text{and}\\
&R_k=\left\{
\begin{array}{ll}
(s,\theta_k],  & ~~ \theta_k>s, \\
(a,\theta_k] \cup (s,b ), & ~~ \theta_k<s.
\end{array}\right.
\end{aligned}
\end{equation}
The set $R_k$ is illustrated in \figref{fig:thatrproRk0}.

Let $\mu_n(R_k)$ and $\mu_n(R_{k}^{\setopen})$ as in~\eqref{eq.defmun}
denote the measure of the sets $R_k$ and~\footnote{
          In the sequel, we let $R^{\setopen}$ denote the interior of a set $R$.
}$R_{k}^{\setopen}$, respectively.
Thus, we have
\begin{subequations}\label{eq.STSaibyalphan0}
\begin{equation}\label{eq.STSaibyalphan1}
\mu_n(R_k)
=\left\{
\begin{array}{ccll}
\mu_n((s,\theta_k]) &=&\alpha_n(\theta_k) - \alpha_n(s),  & ~~  \theta_k>s, \\
\mu_n((a,\theta_k] \cup (s,b )) &=&\alpha_n(\theta_k) + \alpha_n(b) - \alpha_n(s), & ~~ \theta_k<s,
\end{array}\right.
\end{equation}
and
\begin{equation}\label{eq.STSaibyalphan2}
\mu_n(R_{k}^{\setopen})
=\left\{
\begin{array}{ccll}
\mu_n((s,\theta_k)) &=& \alpha_n(\theta_k-) - \alpha_n(s),  & ~~  \theta_k>s, \\
\mu_n((a,\theta_k) \cup (s,b )) &=&\alpha_n(\theta_k-) + \alpha_n(b) - \alpha_n(s), & ~~ \theta_k<s.
\end{array}\right.
\end{equation}
\end{subequations}


In the following theorem we provide a CMS type result
for a class of rational Gaussian quadrature formulae
which applies to the setting of SaI Krylov subspaces with a shift $s\in\R$.

\begin{theorem}[A separation theorem
for rational Gaussian quadrature formulae with a single single pole $s\in\R$ of higher multiplicity]
\label{thm.saixxx}
Let~$\theta_1,\ldots,\theta_m\in(a,b)$ and~$c_1,\ldots,c_m\in\C$
satisfy the rational Gaussian quadrature properties~\eqref{eq.summaryratgaussquad}
for $q_{m-1}(\lambda)=(\lambda-s)^{m-1}$ with $s\in\R$.
Let the index $k_m$ be defined as in~\eqref{eq.thatr2defk1m}.
Let the sets $I_k\subset\{1,\ldots,m\}$
and~$R_k\subset \R$ for~$k=1,\ldots,m$
be defined as in~\eqref{eq.thatr2defRk},
and let $\mu_n$ be defined as in~\eqref{eq.defmun} (analogously,~\eqref{eq.STSaibyalphan0}).
Additionally, define $R_{m+1}:=R_1$.
Then,
\begin{equation}\label{eq.STratbounds}
\mu_n(R_k)
< \sum_{j\in I_k} |c_j|^2
< \mu_n(R_{k+1}^{\setopen}),
~~~~k\in\{1,\ldots,m\}\setminus\{k_m\}.
\end{equation}
\end{theorem}

\smallskip

The case $k=k_m$ is not discussed in Theorem~\ref{thm.saixxx}.
In this case we have $I_k = \{1,\ldots,m\}$
and the bounds~\eqref{eq.STratbounds} can be replaced
by the identity
\begin{equation}\label{eq.normconservedratKry}
\sum_{j=1}^m |c_j|^2 = \alpha_n(b).
\end{equation}
This identity corresponds to the identity~\eqref{eq.summaryratgaussquad}
for $r=1$ (or directly
results from $\|u\|_{\Minr}=\|x\|_2$).

To prove Theorem~\ref{thm.saixxx}, we first introduce rational functions
which constitute bounds on a Heaviside type step function, similar to the
polynomials given in Proposition~\ref{prop.thatpol1}.
\begin{proposition}\label{prop.thatr2}
Let~$\theta_1<\ldots<\theta_m$ be a given sequence
and let $s\in\R$ be a given pole
which is distinct to~$\theta_1,\ldots,\theta_m$.
We make use of the denotations $k_1$ and $k_m$ introduced
in~\eqref{eq.thatr2defk1m}.
Furthermore, let the sets $I_k\subset\{1,\ldots,m\}$
and~$R_k\subset \R$ for~$k=1,\ldots,m$
be defined as in~\eqref{eq.thatr2defRk}, and we define $R_{m+1}:=R_1$.

For~$k\in\{1,\ldots,m\}\setminus\{k_m\}$ and $q_{m-1}(\lambda)=(\lambda-s)^{m-1}$
there exist rational functions~$r_{\{+,k\}}$ and~$r_{\{-,k\}} \in \Pi_{2m-2}/|q_{m-1}|^2$
which satisfy\footnote{
   Analogously to $p_{\{\pm,k\}}$, the denotation $r_{\{\pm,k\}}$ refers to
   $r_{\{+,k\}}$ and $r_{\{-,k\}}$ individually.
}
\begin{equation}\label{eq.thatr2eq}
r_{\{\pm,k\}}(\theta_j)=\left\{
\begin{array}{ll}
1,& j\in I_k,\\
0,&\text{otherwise}.
\end{array}
\right.
\end{equation}
Furthermore, we have
\begin{equation}\label{eq.thatr2ineq}
r_{\{+,k\}}(\lambda) \geq \left\{
\begin{array}{ll}
1,&\lambda \in R_k,\\
0,&\lambda \in\R_s \setminus R_k,
\end{array}
\right.
~~~~\text{and}~~~
r_{\{-,k\}}(\lambda) \leq \left\{
\begin{array}{ll}
1,& \lambda \in R_{k+1}^{\setopen}, \\
0,& \lambda \in \R_s\setminus R_{k+1}^{\setopen},
\end{array}
\right.
\end{equation}
where $\R_s=(a,b)\setminus\{s\}$.
The inequalities in~\eqref{eq.thatr2ineq}
are strict for $\lambda\notin\{\theta_1,\ldots,\theta_m \}$.
Without loss of generality, we assume $(a,b)=\R$ in the present proposition.
\end{proposition}
\begin{proof}
See Appendix~\ref{appendix.A4}.
\end{proof}

\begin{figure}
\centering
\begin{tabular}{cc}
\begin{overpic}
[width=0.5\textwidth]{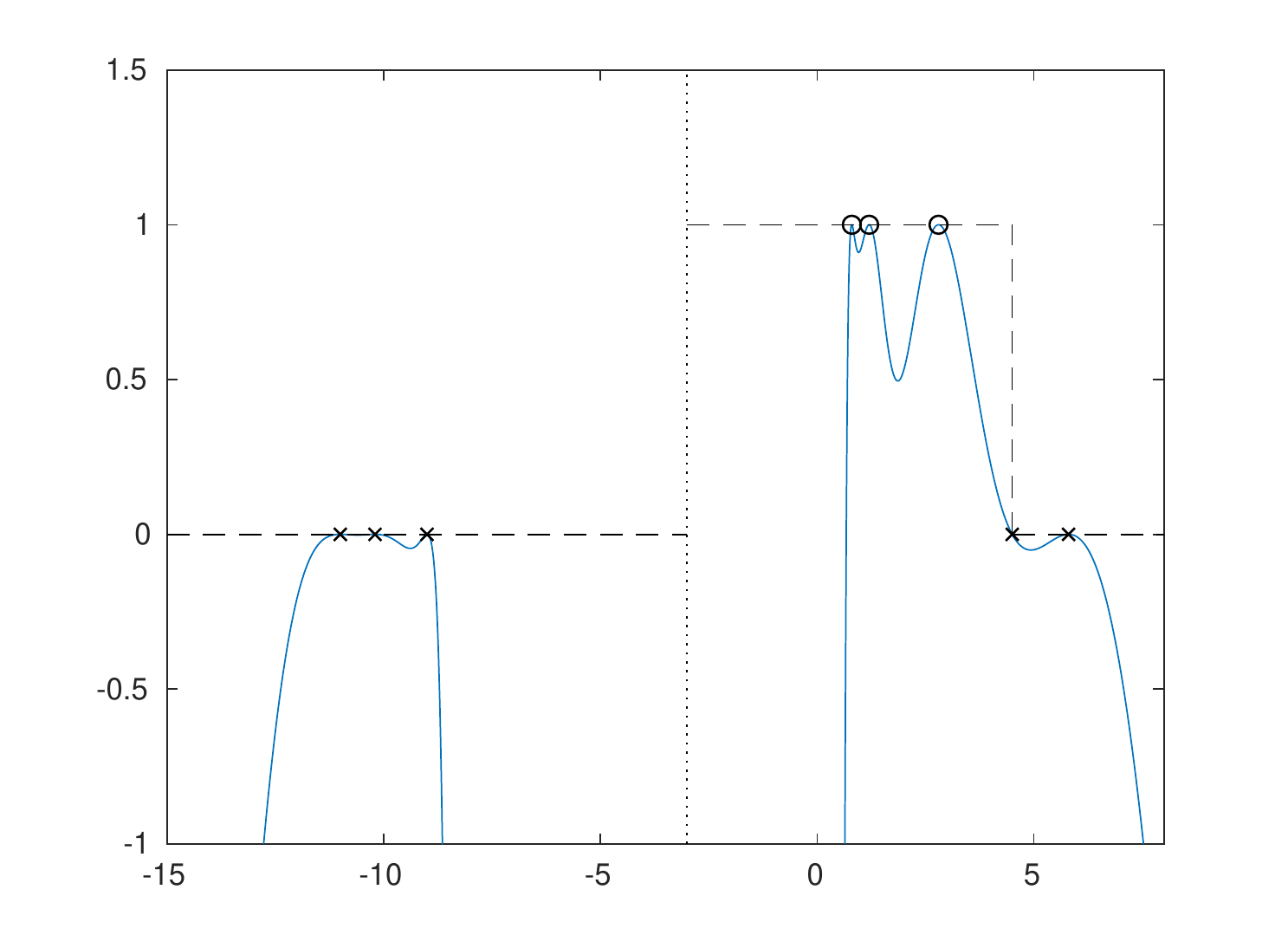}
\put(0,27){\rotatebox{90}{\small$r_{\{-,k\}}(\lambda)$}}
\put(15,63){\small$k=6$}
\put(-4,65){ \textbf{a)}}
\end{overpic}
&
\begin{overpic}
[width=0.5\textwidth]{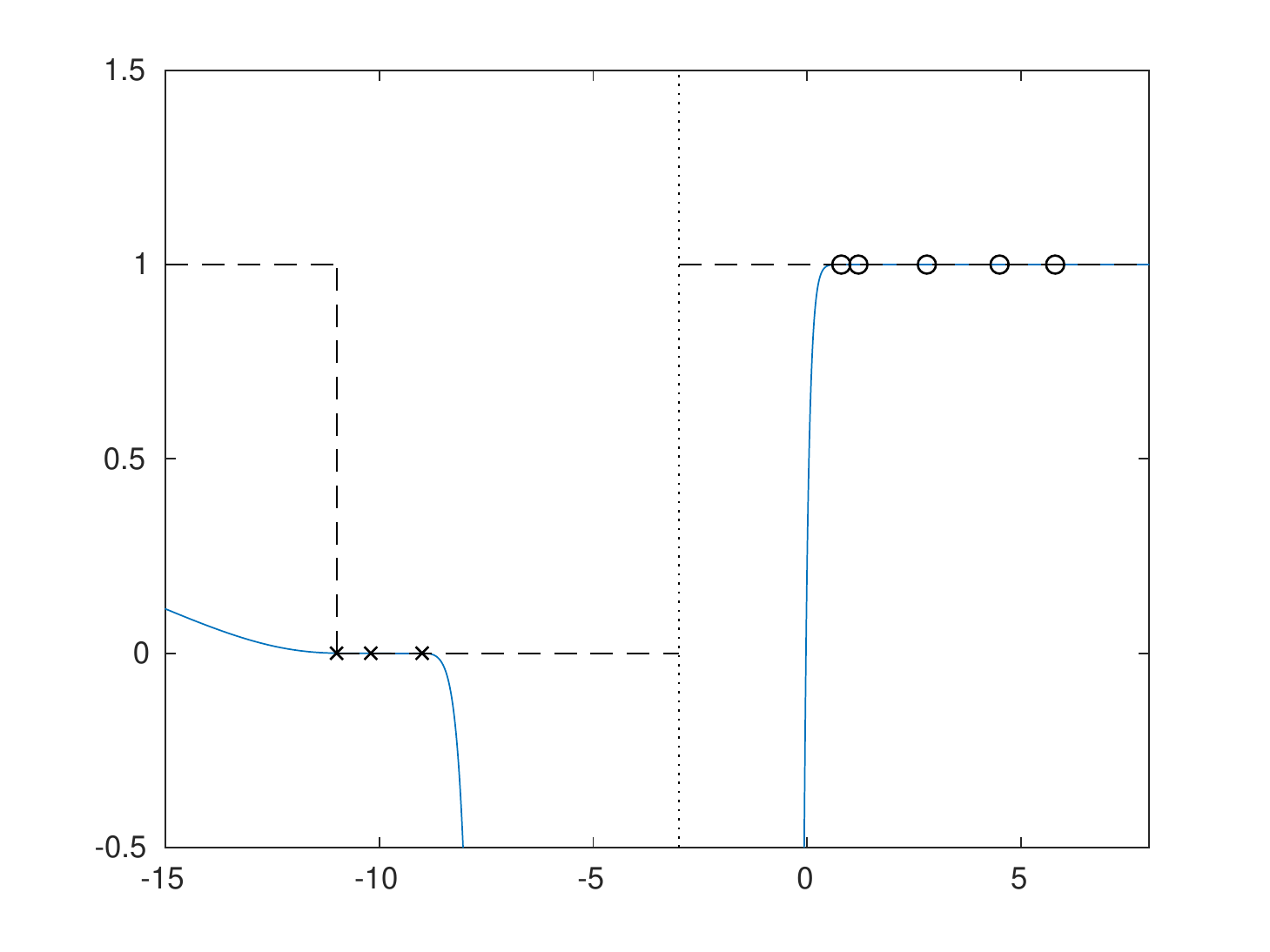}
\put(15,63){\small$k=8$}
\put(50,-1){\small$\lambda$}
\put(-4,65){ \textbf{b)}}
\end{overpic}
\\
\begin{overpic}
[width=0.5\textwidth]{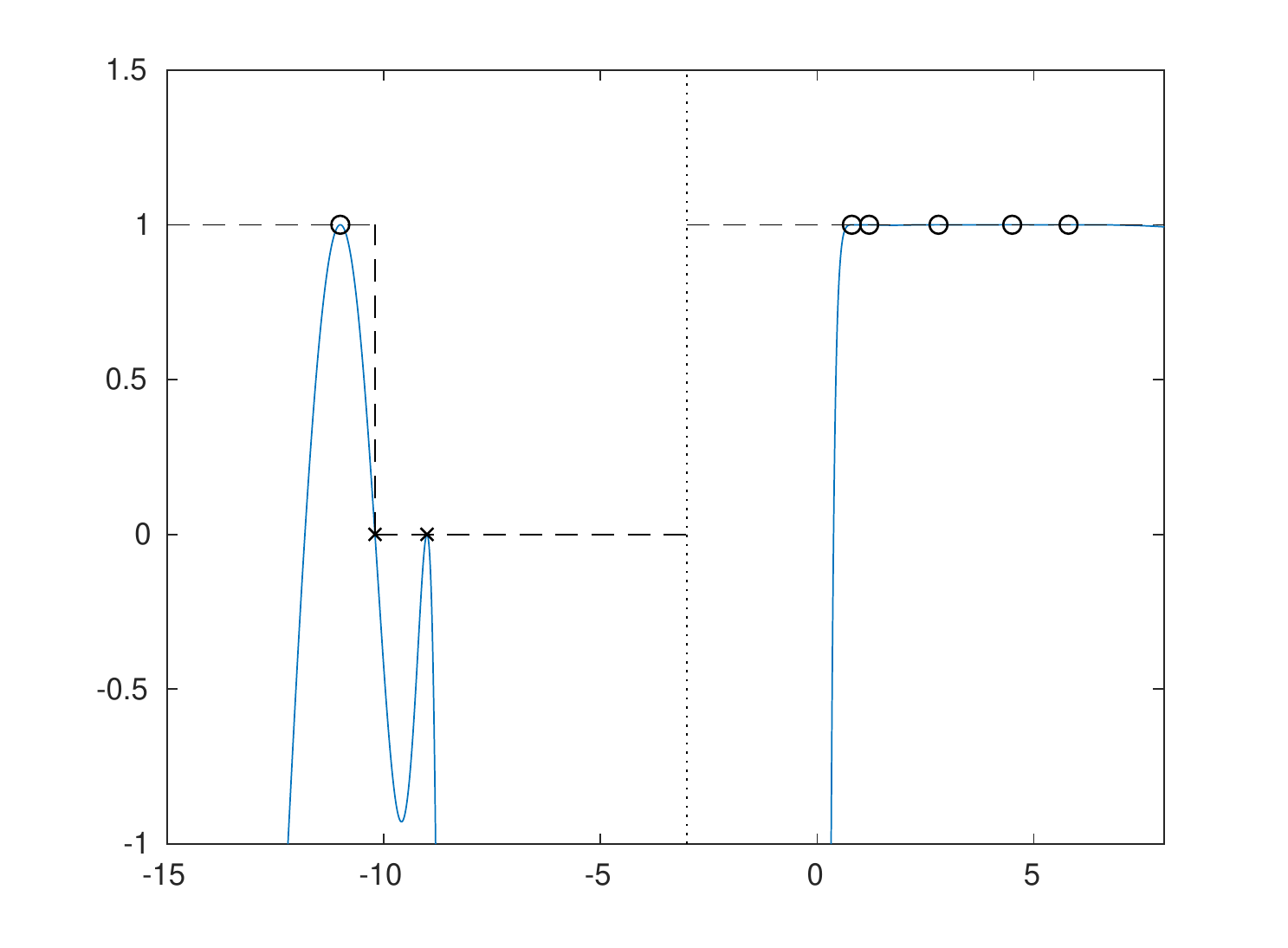}
\put(0,27){\rotatebox{90}{\small$r_{\{-,k\}}(\lambda)$}}
\put(15,63){\small$k=1$}
\put(50,-1){\small$\lambda$}
\put(-4,65){ \textbf{c)}}
\end{overpic}
&
\end{tabular}
\caption{In Figure~\textbf{a)--c)} we show~$r_{\{-,k\}}$ for a given
sequence of nodes $\theta_1,\ldots,\theta_m$ with $m=8$,
and a given pole $s$ which is located between the
nodes,~i.e.,~$\theta_1<s<\theta_m$.
These figures show results for different choices
of~$k\in\{1,\ldots,m\}\setminus\{k_m\}$ where $k_m=3$
(following~\eqref{eq.thatr2defk1mb}).
In each figure the symbols ('$\circ$') and ('$\times$')
mark~$r_{\{-,k\}}(\theta_j)$
for $j\in I_k$ and $j\notin I_k$, respectively.
The dashed lines illustrate the upper bounds given in~\eqref{eq.thatr2ineq}.
Figure~b) shows the special case $k=m$ for which the
upper bound~\eqref{eq.thatr2ineq} relies on $R_{m+1}=R_1$.
For additional illustrations considering~$r_{\{\pm,k\}}$ we refer to
\figref{fig:thatrpro1} and~\ref{fig:thatrpro2} in Appendix~\ref{appendix.A4}.
}
\label{fig:lemrpk}
\end{figure}

Rational functions $r_{\{-,k\}}$
as introduced in Proposition~\ref{prop.thatr2}
are illustrated in \figref{fig:lemrpk} for a numerical example.

We proceed with the proof of Theorem~\ref{thm.saixxx}.
\begin{proof}[Proof of \,\textbf{Theorem~\ref{thm.saixxx}}]
Let $k_1$ and $k_m$ be given in~\eqref{eq.thatr2defk1m},
and let $k\in\{1,\ldots,m\}\setminus\{k_m\}$ be fixed.
For the eigenvalues $\theta_1,\ldots,\theta_m$ we let $r_{\{\pm,k\}}$
denote the rational functions given in~Proposition~\ref{prop.thatr2}.
We proceed to prove the lower bound in~\eqref{eq.STratbounds}.
The identities~\eqref{eq.thatr2eq} imply
\begin{subequations}\label{eq.psai1.1c}
\begin{equation}
\sum_{j\in I_k} |c_j|^2
= \sum_{j=1}^m  r_{\{\pm,k\}}(\theta_j) |c_j|^2 .
\end{equation}
Furthermore, the quadrature property~\eqref{eq.summaryratgaussquad} implies
$$
\int_a^b r_{\{\pm,k\}}(\lambda)\,\dd \alpha_n(\lambda) =\sum_{j=1}^m r_{\{\pm,k\}}(\theta_j) |c_j|^2.
$$
Rewriting this Riemann-Stieltjes integral as in~\eqref{eq.innerprodrepresentationA},
we arrive at
\begin{equation}
\sum_{j=1}^n r_{\{\pm,k\}}(\lambda_j) |w_j|^2 =\sum_{j=1}^m r_{\{\pm,k\}}(\theta_j) |c_j|^2 .
\end{equation}
\end{subequations}
The inequalities in~\eqref{eq.thatr2ineq} for $r_{\{+,k\}}$ entail
\begin{equation}\label{eq.psai1.plus}
\sum_{j=1}^n r_{\{+,k\}}(\lambda_j) |w_j|^2
> \sum_{j\in J(R_k)} |w_j|^2.
\end{equation}
This inequality is strict due to the interlacing property of
the eigenvalues $\lambda_j$ and $\theta_j$, see Proposition~\ref{prop.ritzvaluespol}
and~\ref{prop.ritzqorpol}.
Combine~\eqref{eq.psai1.1c} and~\eqref{eq.psai1.plus} to conclude
with the lower bound in~\eqref{eq.STratbounds}.

In a similar manner, $r_{\{-,k\}}$ reveals the upper bound in~\eqref{eq.STratbounds};
the inequalities in~\eqref{eq.thatr2ineq} for $r_{\{-,k\}}$ yield
\begin{equation}\label{eq.psai1.minus}
\sum_{j=1}^n r_{\{-,k\}}(\lambda_j) |w_j|^2
< \sum_{j\in J(R_{k+1}^{\setopen})} |w_j|^2.
\end{equation}
Indeed, the identities~\eqref{eq.psai1.1c} together with~\eqref{eq.psai1.minus} conclude
the upper bound in~\eqref{eq.STratbounds}.
\end{proof}


We reformulate the result of Theorem~\ref{thm.saixxx}
in the following proposition.

\begin{proposition}\label{prop.SaIRip}
In the setting of Theorem~\ref{thm.saixxx}, the following inequality holds true,
\begin{subequations}\label{eq.SaIxxxv2}
\begin{equation}\label{eq.SaIxxxv2a}
\alpha_n(\theta_k) \leq |c_1|^2 + \ldots + |c_k|^2 + \gamma \leq \alpha_n(\theta_{k+1}-),~~~k=1,\ldots,m-1,
\end{equation}
with
\begin{equation}
\gamma = \alpha_n(s) - \alpha_m(s),
\end{equation}
\end{subequations}
where $\alpha_m$ is given in~\eqref{eq.alpham}.
The inequalities in~\eqref{eq.SaIxxxv2a} are strict for $k\neq k_m$.
Additionally, the case $k=m\neq k_m$ in~\eqref{eq.STratbounds} corresponds to
\begin{equation}\label{eq.SaIxxxv2xc}
\alpha_n(\theta_m) \leq |c_1|^2 + \ldots + |c_m|^2 + \gamma ,~~~~\text{and}~~~ \gamma \leq \alpha_n(\theta_1-).
\end{equation}
\end{proposition}

\begin{proof}[Proof of \,\textbf{Proposition~\ref{prop.SaIRip}}]
We first prove~\eqref{eq.SaIxxxv2}.
Here, we consider different cases for the index $k=1,\ldots,m-1$.
\begin{itemize} 
\item $k<k_1$ with $k_1$ as in~\eqref{eq.thatr2defk1m};
this case only occurs for $k_1 > 1$ which follows from $s\in(\theta_1,\theta_m)$.
This case also implies $k_m=k_1-1$,
\begin{equation}\label{eq.SaIxxxv2p31}
\alpha_m(s)  = |c_1|^2 + \ldots + |c_{k_m}|^2,
~~~\text{and}~~
\alpha_m(b)-\alpha_m(s) = |c_{k_1}|^2 + \ldots + |c_m|^2,
\end{equation}
and with $I_k$ as in~\eqref{eq.thatr2defRk},
\begin{equation}\label{eq.SaIxxxv2p11}
\sum_{j\in I_k} |c_j|^2 = |c_1|^2 + \ldots + |c_k|^2 + \alpha_m(b)-\alpha_m(s).
\end{equation}

-- Additionally, let $k<k_1-1=k_m$.
This case implies $\theta_k,\theta_{k+1}<s$ and as given in~\eqref{eq.STSaibyalphan0}, 
\begin{equation}\label{eq.SaIxxxv2p12}
\mu_n(R_k) = \alpha_n(\theta_k) + \alpha_n(b) - \alpha_n(s),~~~\text{and}~~
\mu_n(R_{k+1}^{\setopen}) = \alpha_n(\theta_{k+1}-) + \alpha_n(b) - \alpha_n(s).
\end{equation}

Substituting~\eqref{eq.SaIxxxv2p11} and~\eqref{eq.SaIxxxv2p12}
in the inequalities~\eqref{eq.STratbounds},
subtracting $\alpha_n(b)$ ($=\alpha_m(b)$) and adding $\alpha_n(s)$,
we conclude~\eqref{eq.SaIxxxv2} for $k<k_1-1$.

-- Let $k=k_1-1$. Thus, $k = k_m\neq m$, and for this case the inequalities~\eqref{eq.STratbounds} do not apply;
we show~\eqref{eq.SaIxxxv2} in a direct manner:

For $k=k_m$ we have $|c_1|^2 + \ldots + |c_k|^2 = \alpha_m(s)$ as in~\eqref{eq.SaIxxxv2p31}.
With this identity, the enclosed term in~\eqref{eq.SaIxxxv2a} simplifies to
\begin{equation}\label{eq.SaIxxxv2p32u}
|c_1|^2 + \ldots + |c_{k}|^2 + \gamma = \alpha_n(s).
\end{equation}
Due to $\alpha_n$ being an increasing function and $\theta_{k}<s<\theta_{k+1}$ for $k=k_1-1$ we have
\begin{equation}\label{eq.SaIxxxv2p32}
\alpha_n(\theta_{k}) \leq \alpha_n(s) \leq \alpha_n(\theta_{k+1}-),~~~\text{for $k=k_1-1$}.
\end{equation}
Combining~\eqref{eq.SaIxxxv2p32u} and~\eqref{eq.SaIxxxv2p32}, we conclude~\eqref{eq.SaIxxxv2} for the case $k = k_1-1 $.

\item $k\geq k_1$. For this case we further distinguish between $s<\theta_m$ and $s>\theta_m$.

-- Let $s<\theta_m$ (this includes the case $s<\theta_1$).
With $I_k$ as in~\eqref{eq.thatr2defRk}, we have
\begin{equation}\label{eq.SaIxxxv2p21x}
\sum_{j\in I_k} |c_j|^2 = |c_1|^2 + \ldots + |c_k|^2 - \alpha_m(s)
\end{equation}
Furthermore, this case implies $\theta_k,\theta_{k+1}>s$ (we recall $k < m$), and as in~\eqref{eq.STSaibyalphan0},
\begin{subequations}\label{eq.SaIxxxv2p2}
\begin{equation}\label{eq.SaIxxxv2p22}
\mu_n(R_k)  = \alpha_n(\theta_k) - \alpha_n(s),
\end{equation}
and
\begin{equation}
\mu_n(R_{k+1}^{\setopen}) = \alpha_n(\theta_{k+1}-) - \alpha_n(s).
\end{equation}
\end{subequations}
Substituting~\eqref{eq.SaIxxxv2p2} and~\eqref{eq.SaIxxxv2p21x} in the inequalities~\eqref{eq.STratbounds},
we conclude~\eqref{eq.SaIxxxv2} for $k\geq k_1$ and $s<\theta_m$.

-- Otherwise, for $k\geq k_1$ and $s>\theta_m$ our notation simplifies to $k_1=1$ and
\begin{equation}\label{eq.SaIxxxv2p21}
\sum_{j\in I_k} |c_j|^2= |c_1|^2 + \ldots + |c_k|^2.
\end{equation}
The case $s>\theta_m$ implies $\alpha_m(b) = \alpha_m(s)$, and due to $\alpha_m(b)=\alpha_n(b)$, we have
\begin{equation}\label{eq.SaIxxxv2p21x2}
\alpha_n(b) = \alpha_m(s).
\end{equation}
Furthermore, we have $\theta_k,\theta_{k+1}<s$, and $\mu_n(R_k)$ and $\mu_n(R_{k+1}^{\setopen})$
correspond to~\eqref{eq.SaIxxxv2p12} further above.
Making use of~\eqref{eq.SaIxxxv2p21x2} in~\eqref{eq.SaIxxxv2p12} and substituting $\gamma$, we simplify
\begin{equation}\label{eq.SaIxxxv2p21x3} 
\mu_n(R_k) = \alpha_n(\theta_k) - \gamma,~~~\text{and}~~
\mu_n(R_{k+1}^{\setopen}) = \alpha_n(\theta_{k+1}-) - \gamma.
\end{equation}
Substituting~\eqref{eq.SaIxxxv2p21} and~\eqref{eq.SaIxxxv2p21x3}
in the inequalities~\eqref{eq.STratbounds},
we conclude~\eqref{eq.SaIxxxv2} for $k\geq k_1$ and $s>\theta_m$.

\end{itemize}

We proceed with the proof of~\eqref{eq.SaIxxxv2xc}.
The case $k=m\neq k_m$ only occurs for $s\in(\theta_1,\theta_m)$.
Thus with $s<\theta_m$, $\mu_n(R_m)$ corresponds to~\eqref{eq.SaIxxxv2p22}.
Substituting $\mu_n(R_m)$ as in~\eqref{eq.SaIxxxv2p22}
and the sum over $I_m$ as in~\eqref{eq.SaIxxxv2p21x}
in the lower bound in~\eqref{eq.STratbounds},
we conclude
the inequality on the left-hand side of~\eqref{eq.SaIxxxv2xc}.

To prove the inequality on the right-hand side of~\eqref{eq.SaIxxxv2xc},
we first recall $\theta_1<s$, and as in~\eqref{eq.SaIxxxv2p12}
\begin{equation}\label{eq.SaIxxxv2p4}
\mu_n(R_{1}^{\setopen}) = \alpha_n(\theta_1-) + \alpha_n(b) - \alpha_n(s).
\end{equation}
Substituting~\eqref{eq.SaIxxxv2p21x} and~\eqref{eq.SaIxxxv2p4}
in the upper bound in~\eqref{eq.STratbounds}
(for the case $k=m\neq k_m$ with $\mu_n(R_{m+1}^{\setopen})=\mu_n(R_{1}^{\setopen})$ due to convention),
we arrive at
$$
|c_1|^2 + \ldots + |c_m|^2 - \alpha_m(s) < \alpha_n(\theta_1-) + \alpha_n(b) - \alpha_n(s)
$$
On the left-hand side we can further simplify $|c_1|^2 + \ldots + |c_m|^2=\alpha_n(b)$ 
and subtract this term, which entails the inequality on the right-hand side of~\eqref{eq.SaIxxxv2xc}.
\end{proof}


\begin{remark}\label{rmk.SaIboundssnotinspec}
For the case $\alpha_m(s)=\alpha_n(s)$ the constant $\gamma$
in Proposition~\ref{prop.SaIRip} is zero, and the inequalities~\eqref{eq.SaIxxxv2a}
coincide with the inequalities given by Theorem~\ref{thm.STCMS},
i.e.,~the CMS Theorem for polynomial Gaussian quadrature formulae.
Furthermore, for this case the inequalities
given in Corollary~\ref{cor.piecewiseboundspol} hold true.
Here, we highlight the case $s\notin(\lambda_1,\lambda_n)$
for the Gaussian quadrature formulae without preassigned nodes
(this implies $\theta_j\in(\lambda_1,\lambda_n)$);
a prominent case for which $\alpha_m(s)=\alpha_n(s)$ holds true a~priori.
\end{remark}

\begin{remark}\label{rmk.CMSforSaIrepres}
Following Remark~\ref{rmk.quadSaIKry},
the SaI Krylov representation $X_m^{-1}+sI$ provides a Gaussian quadrature formula.
For the case of a real shift $s\in\R$, the respective quadrature nodes are located
on the real axis, and at least one eigenvalue $\lambda_j$
is located between each neighboring pair of quadrature nodes.
Thus, the result of Theorem~\ref{thm.saixxx} and its corollaries hold true
in this setting.
However, results concerning the SaI Krylov representation $X_m^{-1}+sI$
are not discussed in further detail in the present work.

CMS type results for the SaI Krylov representation
are also given in~\cite{ZTK19}. In the present work we include
the case of a shift $s$ being located inside the convex hull of the spectrum of $A$,
which extends some results of~\cite{ZTK19}.
\end{remark}

We proceed to specify the results of Proposition~\ref{prop.SaIRip}
for the pole $s$ being located in the convex hull of the
rational Ritz values,~i.e.,~$s\in(\theta_1,\theta_m)$.
This case implies $k_m\neq m$,
and substituting $\alpha_n(b)=|c_1|^2 + \ldots + |c_m|^2$ in~\eqref{eq.SaIxxxv2xc},
we observe
\begin{equation}\label{eq.boundsongam}
\alpha_n(\theta_m) - \alpha_n(b) \leq \gamma \leq \alpha_n(\theta_1-).
\end{equation}
With these inequalities, we further specify the results of Proposition~\ref{prop.SaIRip}:
The following corollary states some bounds on piecewise accumulated quadrature weights,
similar to Corollary~\ref{cor.piecewiseboundspol} in the previous subsection
for the polynomial case.
\begin{corollary}\label{cor.piecewiseboundsrat}
Additionally to the setting of Theorem~\ref{thm.saixxx}, we assume $s\in(\theta_1,\theta_m)$.
Then Proposition~\ref{prop.SaIRip} yields the following inequalities.
\begin{itemize}
\item The accumulated quadrature weights satisfy
\begin{subequations}
\begin{equation}\label{eq.boundscorsaixa}
\mu_n([\theta_1,\theta_k]) \leq |c_1|^2 + \ldots + |c_k|^2 
 \leq  \mu_n((a,\theta_{k+1})\cup (\theta_m,b)) ,~~~k=1,\ldots,m-1.
\end{equation}

\item For indices $j,k$ with $1<j<k<m$, the following piecewise accumulated quadrature weights
satisfy
\begin{equation}\label{eq.boundscorsaixb}
\mu_n([\theta_j,\theta_k]) \leq |c_j|^2 + \ldots + |c_k|^2 
\leq \mu_n((\theta_{j-1},\theta_{k+1})).
\end{equation}

\item Furthermore, the 
accumulated quadrature weights satisfy
\begin{equation}\label{eq.boundscorsaixc}
\mu_n([\theta_j,\theta_m]) \leq |c_j|^2 + \ldots + |c_m|^2 
\leq \mu_n( (a,\theta_1) \cup (\theta_{j-1},b) )  ,~~~j=2,\ldots,m.
\end{equation}
\end{subequations}
\end{itemize}
\end{corollary}
\begin{proof}
The inequalities~\eqref{eq.SaIxxxv2a} in Proposition~\ref{prop.SaIRip} yield
$$
\alpha_n(\theta_k) - \gamma \leq |c_1|^2 + \ldots + |c_k|^2 \leq \alpha_n(\theta_{k+1}-) - \gamma.
$$
Substituting~\eqref{eq.boundsongam} for $\gamma$, we arrive at
$$
\alpha_n(\theta_k) - \alpha_n(\theta_1-) \leq |c_1|^2 + \ldots + |c_k|^2 
 \leq \alpha_n(\theta_{k+1}-) + \alpha_n(b) - \alpha_n(\theta_m).
$$
This implies~\eqref{eq.boundscorsaixa}.

To prove the inequalities in~\eqref{eq.boundscorsaixc}, we first remark
$$
|c_j|^2 + \ldots + |c_k|^2 = |c_1|^2 + \ldots + |c_k|^2 + \gamma - \big(|c_1|^2 + \ldots + |c_{j-1}|^2 + \gamma\big).
$$
Applying~\eqref{eq.SaIxxxv2a} twice (once we substitute $j-1$ for the index $k$ therein)
we observe
$$
\alpha_n(\theta_k) - \alpha_n(\theta_j-) \leq |c_j|^2 + \ldots + |c_k|^2
\leq \alpha_n(\theta_{k+1}-) - \alpha_n(\theta_j),
$$
which implies~\eqref{eq.boundscorsaixb}.

To show~\eqref{eq.boundscorsaixc}, apply~\eqref{eq.boundscorsaixa}
for the index $j-1$ and subtract the result from
$|c_1|^2 + \ldots + |c_m|^2 = \mu_n((a,b))$.
\end{proof}

\begin{remark}\label{rmk.CMSqorsplit}
For the case $s\in(\theta_1,\theta_m)$ as in Corollary~\ref{cor.piecewiseboundsrat}, bounds on quadrature weights
related to the leftmost or rightmost quadrature nodes potentially
depend on the measure of an interval including the opposite integral limit.
This relation can be avoided by preassigning one of the quadrature nodes at the integral limit,
using a rational Gauss-Radau formula associated with the spectrum of a rational qor-Krylov representation $B_m$
in the Krylov setting.
\begin{itemize}
\item For a preassigned node $\xi<\lambda_1$, we have $\theta_1=\xi$ and $\alpha_n(\theta_1)=0$.
Thus, the inequalities in~\eqref{eq.boundscorsaixa} correspond to
\begin{subequations}
\begin{equation}
\mu_n([a,\theta_k]) \leq |c_1|^2 + \ldots + |c_k|^2 
 \leq  \mu_n((a,\theta_{k+1})\cup (\theta_m,b)),
\end{equation}
and the inequalities in~\eqref{eq.boundscorsaixc} correspond to
\begin{equation}
\mu_n([\theta_j,\theta_m]) \leq |c_j|^2 + \ldots + |c_m|^2 
\leq \mu_n( (\theta_{j-1},b) ).
\end{equation}

\item For a preassigned node $\xi>\lambda_n$, we have $\theta_m=\xi$ and $\alpha_n(\theta_m)=\alpha_n(b)$.
Thus, the inequalities in~\eqref{eq.boundscorsaixa} correspond to
\begin{equation}
\mu_n([\theta_1,\theta_k]) \leq |c_1|^2 + \ldots + |c_k|^2 
 \leq  \mu_n((a,\theta_{k+1})),
\end{equation}
and the inequalities in~\eqref{eq.boundscorsaixc} correspond to
\begin{equation}
\mu_n([\theta_j,b]) \leq |c_j|^2 + \ldots + |c_m|^2 
\leq \mu_n( (a,\theta_1) \cup (\theta_{j-1},b) ).
\end{equation}
\end{subequations}
\end{itemize}
\end{remark}

We proceed to introduce a step function $F_s$ which changes its sign
at each rational Ritz value according to Proposition~\ref{prop.SaIRip};
with the step function $F$ given in~\eqref{eq:Fxalpha} we introduce
\begin{equation}\label{eq.verifySaI2Fs}
F_s(\lambda)=F(\lambda)-F(s).
\end{equation}
Here, $ F(s) = \gamma $ with $\gamma$ as in Proposition~\ref{prop.SaIRip}.
As previously stated in~\eqref{eq.alpmthetakpm}, we have
$$
\alpha_m(\theta_k) = \alpha_m(\theta_{k+1}-) = |c_1|^2 + \ldots + |c_k|^2.
$$
Then the inequalities~\eqref{eq.SaIxxxv2a} correspond to
\begin{subequations}\label{eq.Fsineqall}
\begin{equation}
F_s(\theta_k) \leq 0 \leq F_s(\theta_{k+1}-),~~~k=1,\ldots,m-1,
\end{equation}
whereat these inequalities are strict for $k\neq k_m$.
Furthermore, the inequalities~\eqref{eq.SaIxxxv2xc} correspond to
\begin{equation}
F_s(\theta_m) < 0,~~~\text{and}~~0 < F_s(\theta_{1}-),
\end{equation}
\end{subequations}
for $k_m\neq m$.

\begin{remark}\label{rmk.FsandThmsaix}
In~\eqref{eq.Fsineqall}, the special case $k=k_m$ holds true due to the identity~\eqref{eq.normconservedratKry};
the case $k\in\{1,\ldots,m\}\setminus\{k_m\}$ corresponds to the result of Theorem~\ref{thm.saixxx}.
Namely,~the result of Theorem~\ref{thm.saixxx} conforms to the following inequality in an equivalent manner,
$$
F_s(\theta_k) \leq 0 \leq F_s(\theta_{k+1}-)~~~\text{for}~~k\in\{1,\ldots,m\}\setminus\{k_m\},~\text{and with}~\theta_{m+1}=\theta_1.
$$
\end{remark}


\subsection{The rational case with a single pole \texorpdfstring{$s\in\C\setminus\R$}{s in C minus R} of higher multiplicity}
\label{subsec.CMSSaIcomplex}

In the present subsection, we consider rational Gaussian quadrature formulae which
satisfy the quadrature property~\eqref{eq.summaryratgaussquad}
with $q_{m-1}(\lambda) = (\lambda-s)^{m-1}$ for $s\in\C\setminus\R$.
To specify, these quadrature formulae are exact for rational functions
with denominator $|q_{m-1}(\lambda)|^2 = ((\lambda-\real s )^2 + (\imag s)^2)^{m-1}$ where $\imag s\neq 0$,
i.e., rational functions with complex-conjugate poles of higher multiplicity.
Considering Krylov subspaces, these quadrature formulae are related to
SaI Krylov subspaces with a complex shift $s\in\C\setminus\R$.

As a main result of the present subsection,
the following Proposition yields upper bounds on
the measure of
the intervals between neighboring quadrature nodes,
and the measure at the boundary of the spectrum.
\begin{proposition}\label{prop.SaIC1}
Let $c_1,\ldots,c_m$ and $\theta_1<\ldots<\theta_m$ satisfy
the quadrature property~\eqref{eq.summaryratgaussquad}
with $q_{m-1}(\lambda) = (\lambda-s)^{m-1}$ for $s\in\C\setminus\R$.
Then, with $\mu_n$ given in~\eqref{eq.defmun0}
\begin{subequations}
\begin{equation}\label{eq.SaIC1bk}
\mu_n([\theta_k,\theta_{k+1}]) \leq |c_k|^2+|c_{k+1}|^2,~~~~k=1,\ldots,m-1,
\end{equation}
and
\begin{equation}\label{eq.SaIC1b1m}
\mu_n((a,\theta_1]) + \mu_n([\theta_m,b)) \leq |c_1|^2+|c_m|^2.
\end{equation}
\end{subequations}
\end{proposition}

Before proving Proposition~\ref{prop.SaIC1},
we proceed with some auxiliary results.
The results of the previous subsection do not apply for the case $s\in\C\setminus\R$.
However, the present class of rational functions
can be related to polynomials on the unit circle $\T$ and vice versa.
To specify this relation, we recall the Cayley transform as in~\eqref{eq.cayleyscalarSaI},
$$
\tau(\lambda)=(\lambda-\tbar{s})(\lambda-s)^{-1},
~~~~ \tau\colon \R \to \T\setminus\{1\}.
$$
For a complex polynomial $p\in\Pi_{m-1}$ we 
consider $ p(\tau(\lambda)) $ 
as a function of $\lambda$; normalizing shows
$$
p(\tau(\lambda)) = g(\lambda)/q_{m-1}(\lambda),~~~\text{for some $g\in \Pi_{m-1}$},
$$
For $\lambda\in\R$ we conclude
\begin{equation}\label{eq.ponTtorat}
|p(\tau(\lambda))|^2 = \tbar{g}(\lambda) g(\lambda)/|q_{m-1}(\lambda)|^2,~~~
\text{where}~~\tbar{g}g\in\Pi_{2m-2}.
\end{equation}
In the following corollary we introduce rational majorants on a
Heaviside type step function, based on interpolating polynomials
on the unit circle given in~\cite[Lemma~4]{Go02}.
\begin{corollary}[A corollary of Lemma~4 in~\cite{Go02}]\label{cor.thatrcomplex}
Let $\theta_1,\ldots,\theta_m$ be a given sequence of nodes,
and let $q_{m-1}(\lambda)=(\lambda-s)^{m-1}$
for a given pole $s\in\C\setminus\R$.

\begin{enumerate}[label=(\roman*)]
\item Let $k\in\{1,\ldots,m-1\}$ be fixed.
There exists a rational function $r_k\in\Pi_{2m-2}/|q_{m-1}|^2$
with
\begin{subequations}\label{eq.thatrcomplexkprop}
\begin{equation}\label{eq.thatrcomplexkequal}
r_k(\theta_j)=\left\{
\begin{array}{ll}
1,&~~~ j\in\{k,k+1\},\\
0,&~~~ \text{otherwise},
\end{array}\right.
\end{equation}
and
\begin{equation}\label{eq.thatrcomplexkbound}
r_k(\lambda) \geq \left\{
\begin{array}{ll}
1,&~~~ \lambda \in [\theta_k,\theta_{k+1}],\\
0,&~~~ \lambda \in (-\infty,\theta_k)\cup(\theta_{k+1},\infty).
\end{array}\right.
\end{equation}
\end{subequations}

\item Additionally, there exists a function $r_m\in\Pi_{2m-2}/|q_{m-1}|^2$ with
\begin{subequations}\label{eq.thatrcomplex1mprop}
\begin{equation}\label{eq.thatrcomplex1mequal}
r_m(\theta_j)=\left\{
\begin{array}{ll}
1,&~~~ j\in\{1,m\},\\
0,&~~~ \text{otherwise},
\end{array}\right.
\end{equation}
and
\begin{equation}\label{eq.thatrcomplex1mbound}
r_m(\lambda) \geq \left\{
\begin{array}{ll}
1,&~~~ \lambda \in (-\infty,\theta_1]\cup[\theta_m,\infty),\\
0,&~~~ \lambda \in (\theta_1,\theta_m).
\end{array}\right.
\end{equation}
\end{subequations}

\end{enumerate}
\end{corollary}
\begin{proof}
The Cayley transform $\tau$ as in~\eqref{eq.cayleyscalarSaI} reads
$$
\tau(\lambda)=(\lambda-\tbar{s})(\lambda-s)^{-1},
~~~~ \tau\colon \R \to \T\setminus\{1\}.
$$
Simplifying this fraction yields
$$
\tau(\lambda)= 1 + \frac{s-\tbar{s}}{\lambda-s} =  1 + \frac{2\ii \,\imag s}{\lambda-s}.
$$
Here, $\tau\colon \R \to \T\setminus\{1\}$ is a continuous and bijective function,
and with the previous representation, we observe
$$
\tau(-\infty) = \left\{
\begin{array}{ll}
1 + \ii 0+,& ~~~ \imag s < 0,\\
1 + \ii 0-,& ~~~ \imag s > 0.
\end{array}
\right.
$$
Thus, $\tau$ maps $\R$ to $\T\setminus\{1\}$ in counter-clockwise and clockwise order
for $\imag s<0$ and $\imag s>0$, respectively.

We proceed to define distinct points $\zeta_1,\ldots,\zeta_m\in\T$ by
\begin{equation}\label{eq.tharcplxdefzeta}
\zeta_j := \left\{
\begin{array}{ll}
\tau(\theta_j),& ~~~ \imag s < 0,\\
\tau(\theta_{m-j+1}),& ~~~ \imag s > 0,
\end{array}
\right. 
~~~~ j=1,\ldots,m.
\end{equation}
The points $\zeta_1,\ldots,\zeta_m\in\T$ are distinct, in counter-clockwise order,
and the point $1$ is located between $\zeta_1$ and $\zeta_m$ on the unit circle.
For the remainder of the proof we assume the case $\imag s < 0$ to simplify the notation.
Thus, we consider $\zeta_j = \tau(\theta_j)$.

Additionally to~\eqref{eq.tharcplxdefzeta}, we define $\zeta_{m+1}:=\zeta_1$.
Let $p_k\in\Pi_{m-1}$ denote the complex polynomial given by~\cite[Lemma~4]{Go02} for the points $\zeta_j$
and a fixed index $k\in\{1,\ldots,m\}$. Here, we also normalize $p_k$ at $\zeta_{k}$.
Thus, $p_k$ satisfies $|p_k(\zeta_{k})|=|p_k(\zeta_{k+1})|=1$.
Following~\eqref{eq.ponTtorat}, the function $r_k(\lambda):=|p_k(\tau(\lambda))|^2$ conforms to a rational function
$r_k\in\Pi_{2m-2}/|q_{m-1}|^2$.

We proceed to show~\eqref{eq.thatrcomplexkprop} for the rational function $r_k$;
let $k\in\{1,\ldots,m-1\}$:
\begin{itemize}
\item With $r_k(\theta_j) = p_k(\zeta_j)$ the identities $|p_k(\zeta_{k})|=|p_k(\zeta_{k+1})|=1$ yield $r_k(\theta_k)=r_k(\theta_{k+1})=1$,
and the identity $p_k(\zeta_{j})=0$ for $j\neq k,k+1$ yields $r_k(\theta_j)=0$ for $j\neq k,k+1$,
which shows~\eqref{eq.thatrcomplexkequal}.
\item Due to $\tau$ being a continuous and bijective function,
the points $\zeta$ located on the unit circle between $\zeta_k$ and $\zeta_{k+1}$
(including $\zeta_k$ and $\zeta_{k+1}$)
are identical to the set $\{\zeta=\tau(\lambda)~|~\lambda \in [\theta_k,\theta_{k+1}]\}$.
As a result of~\cite[Lemma~4]{Go02}, the polynomial~$p_k$ satisfies $|p_k(\zeta)|\geq 1 $
for $\zeta$ in this set of
points,~i.e.,~$|p(\tau(\lambda))|\geq 1 $ for $\lambda \in [\theta_k,\theta_{k+1}]$.
Thus, we have $r_k(\lambda)\geq 1 $ for $\lambda \in [\theta_k,\theta_{k+1}]$;
furthermore, $r_k$ is positive for $\lambda\in\R$ due to $r_k(\lambda)=|p_k(\xi(\lambda))|^2$,
which implies~\eqref{eq.thatrcomplexkbound}.
\end{itemize}

We proceed to sketch the proof of~\eqref{eq.thatrcomplex1mprop} which corresponds to the case $k=m$.
The polynomial $p_m$ satisfies $|p_m(\zeta_m)|=|p_m(\zeta_1)|=1$.
Furthermore, the points $\zeta$ located between $\zeta_1$ and $\zeta_m$ correspond to the set
$ \{\zeta=\tau(\lambda)~|~\lambda \in (-\infty,\theta_1) \cup (\theta_m,\infty) \} \cup\{1\}\subset \T$.
Similar to previous arguments, this shows~\eqref{eq.thatrcomplex1mprop}.

Considering the definition of $\zeta_j$ in~\eqref{eq.tharcplxdefzeta},
similar arguments hold for the case~$\imag s>0$.
\end{proof}

We proceed with the proof of Proposition~\ref{prop.SaIC1}.
\begin{proof}[Proof of \,\textbf{Proposition~\ref{prop.SaIC1}}]
Let $k\in\{1,\ldots,m-1\}$ be fixed, we prove~\eqref{eq.SaIC1bk}.
For the nodes $\theta_1,\ldots,\theta_m$ and $k$ given,
we let $r_k\in\Pi_{2m-2}/|q_{m-1}|^2$ denote the
rational function given in Corollary~\ref{cor.thatrcomplex}
which satisfies~\eqref{eq.thatrcomplexkprop}.
Due to~\eqref{eq.thatrcomplexkequal} we have
\begin{subequations}
\begin{equation}\label{eq.proofSaIC1ka}
\sum_{j=1}^m |c_j|^2 r_k(\theta_j) = |c_k|^2 + |c_{k+1}|^2.
\end{equation}
The quadrature property~\eqref{eq.summaryratgaussquad} implies
\begin{equation}
\int_a^b r_k(\lambda) \,\dd \alpha_n(\lambda) = \sum_{j=1}^m |c_j|^2 r_k(\theta_j),
\end{equation}
and the inequality~\eqref{eq.thatrcomplexkbound} yields
\begin{equation}\label{eq.proofSaIC1kc}
\int_a^b r_k(\lambda) \,\dd \alpha_n(\lambda) \geq \mu_n([\theta_k,\theta_{k+1}]).
\end{equation}
\end{subequations}
Combining~\eqref{eq.proofSaIC1ka}--\eqref{eq.proofSaIC1kc}, we conclude~\eqref{eq.SaIC1bk}.

Analogously, making use of the rational function $r_m\in\Pi_{2m-2}/|q_{m-1}|^2$
(which satisfies the properties~\eqref{eq.thatrcomplex1mprop})
in combination with the quadrature property~\eqref{eq.summaryratgaussquad},
we conclude~\eqref{eq.SaIC1b1m}.

\end{proof}

\subsection{Results for an extended Krylov subspace}\label{subsec.STCMSxKry}

In the present subsection, we consider an extended Krylov subspace.
Namely, the Krylov subspace of so called Laurent polynomials which also appears in~\cite{DK98}
and corresponds to a rational Krylov subspace.
Here, we also include a shift $s <\lambda_1$.
This yields a rational Krylov subspace with denominator $q(\lambda) = (\lambda-s)^{\varrho-1}$ for $m=2\varrho-1$,~i.e.,
\begin{align}\label{eq.extendedKrysubspace}
\RKry_{2\varrho-1}(A,u) &= \vspan\{ (A-s)^{-\varrho+1}u,\ldots,(A-s)^{-1}u,u,Au,\ldots, A^{\varrho-1}u\}\\
&= \Kry_{2\varrho-1}(A, (A-s)^{-\varrho+1}u).\notag
\end{align}
Similar to previous sections, $U_m$ denotes an \tMinr-orthonormal basis of the Krylov subspace
and $A_m=(U_m,A\,U_m)_{\Minr}$ denotes the associated Rayleigh quotient.
As previously, we let $x=(U_m,u)_{\Minr}$.
An extended Lanczos recurrence to compute $U_m$ and $A_m$ in an efficient
manner is given in~\cite[Section~5]{DK98} and summarized in Algorithm~\ref{alg.extendedLanczos}.

\begin{figure}
\centering
\begin{algorithm}[H]
\caption{A summary of the extended Lanczos recurrence in~\cite[Section~5]{DK98};
an algorithm to compute the \tMinr-orthogonal basis $U_m$ and 
the Rayleigh quotient $A_m=(U_m,A\,U_m)_{\Minr}$ of the extended Krylov subspace
given in~\eqref{eq.extendedKrysubspace}. Here, $m=2\varrho-1$.}
\label{alg.extendedLanczos}
\SetAlgoLined
 run Algorithm~\ref{alg.SaI}
 to compute $\beta_0=\|u\|_{\Minr}$, $U^{\text{SaI}}_{\varrho}$ and
$A^{\text{SaI}}_{\varrho}=(U^{\text{SaI}}_{\varrho},A\,U^{\text{SaI}}_{\varrho})_{\Minr}$
for the SaI Krylov subspace $\Kry_{\varrho}(X,u)$
 with $X=(A-sI)^{-1}$\;
 $\widetilde{v}=Au$\;
 orthogonalize $\widetilde{v}$ with $U^{\text{SaI}}_{\varrho}$ and set $U_{\varrho}=U_{\varrho}^{\text{SaI}}$
and $u_{\varrho+1}=\widetilde{v}/\|\widetilde{v}\|_{\Minr}$\;
 $\widehat{v}=A\,u_{\varrho+1}$\;
 for $j=1,\ldots,\varrho$\;
 ~~~~  $y_j=(u_j,\widehat{v})_{\Minr}$\;
 ~~~~  $\widehat{v} \leftarrow \widehat{v} - y_j u_j$\;
 $a_{1} = (u_{\varrho+1},\widehat{v}\,)_{\Minr}$ and $\widehat{v} \leftarrow \widehat{v} - a_{1} u_{k+1}$\;
 $\beta_{1}= \|\widehat{v}\|_{\Minr}$ and $u_{\varrho+2}=\widehat{v}/\beta_{1}$\;
 consider $u_{\varrho+1},u_{\varrho+2}$, and $a_1$ and $\beta_1$ to be the result of two initial Lanczos steps,
 and continue the Lanczos procedure to compute $u_{\varrho+3},\ldots,u_{2\varrho-1}$ and 
 the Jacobi matrix $J_{\varrho-1}$ (using a total of $\varrho-1$ Lanczos steps)\;
 $A_m = [\,A^{\text{SaI}}_{\varrho}, y\,e_1^{\Hast};~ e_1\,y^\Hast, J_{\varrho-1}\,] $,~~~where $y\,e_1^{\Hast}\in\C^{\varrho\times \varrho-1}$\;
 $x=\beta_0\,e_1$\;
 return $x,U_m,A_m$\;
\end{algorithm}
\end{figure}

Let $\theta_1,\ldots,\theta_m\in(a,b)$ and $c_1,\ldots,c_m\in\C$ denote
the eigenvalues and spectral coefficients of the Rayleigh quotient $A_m$.
Following Proposition~\ref{eq.ratkryeqinnerprod} and Corollary~\ref{cor.identityAmBm},
these eigenvalues and spectral coefficients
satisfy the identity~\eqref{eq.summaryratgaussquad}
 for $r\in \Pi_{2m-1}/q^2$
with $q^2(\lambda) = (\lambda-s)^{2\varrho-2} = (\lambda-s)^{m-1}$,~i.e.,
\begin{equation}\label{eq.idinrprodextendedkry}
\int_a^b r(\lambda)\,\dd \alpha_n(\lambda)
= \sum_{j=1}^m r(\theta_j) |c_j|^2,~~~~ r\in \Pi_{2m-1}/(\lambda-s)^{m-1}.
\end{equation}
The CMS type results given in~\cite{Li98} apply to
rational quadrature formulae which satisfy~\eqref{eq.idinrprodextendedkry}.
i.e.,~a Gaussian quadrature formulae for so called Laurent polynomials.

We proceed to recapitulate results given in~\cite{Li98}
for the setting of the extended Krylov subspace~\eqref{eq.extendedKrysubspace}.
To this end, we first recall the following rational functions introduced~\cite{Li98}
which yield majorants and minorants on a Heaviside step function
similar to the polynomials in Proposition~\ref{prop.thatpol1}.
\begin{proposition}[Theorem~4 and~5  in~\cite{Li98}]\label{prop.thatLaurentpol}
Let~$\nu_1<\ldots<\nu_m\in\R$ with $\nu_1>0$ and let~$k$ be fixed
with~$1\leq k<m$.
Then there exist rational functions~$\widehat{r}_{\{+,k\}}$ and~$\widehat{r}_{\{-,k\}}\in \Pi_{2m-2}/\lambda^{m-1}$
which satisfy
\begin{equation}\label{eq.lpmequal}
\widehat{r}_{\{\pm,k\}}(\nu_j)=\left\{
\begin{array}{ll}
1,&j=1,\ldots,k,\\
0,&j=k+1,\ldots,m,
\end{array}
\right.
\end{equation}
together with
\begin{equation}\label{eq.lpminequal}
\widehat{r}_{\{+,k\}}(\lambda) \geq \left\{
\begin{array}{ll}
1,&\lambda \leq \nu_k,\\
0,& \lambda > \nu_k,
\end{array}
\right.
~~~\text{and}~~~
\widehat{r}_{\{-,k\}}(\lambda) \leq \left\{
\begin{array}{ll}
1,& \lambda < \nu_{k+1},\\
0,& \lambda \geq \nu_{k+1}.
\end{array}
\right.
\end{equation}
Additionally, the inequalities in~\eqref{eq.ppminequal} are strict inequalities for $\lambda\notin\{\nu_1,\ldots,\nu_m\}$.
\end{proposition}
In the proof of Proposition~\ref{prop.STCMSex} below, we apply these results for the shifted case with $s\leq a<\lambda_1$.

We proceed to recapitulate~\cite[eq.~(4) in Theorem~1]{Li98}.
For the following proposition, we recall that $\lambda_1<\theta_1$ holds true when $\theta_j$
refers to the eigenvalues of the Rayleigh quotient $A_m$, thus, for a pole $s<\lambda_1$
the condition $s<\theta_1$ is satisfied.
\begin{proposition}[Eq.~(4) in Theorem~1 in~\cite{Li98}]\label{prop.STCMSex}
~ Let $\theta_1,\ldots,\theta_m\in(a,b)$ and $c_1,\ldots,c_m\in\C$ satisfy~\eqref{eq.idinrprodextendedkry}
for $r\in \Pi_{2m-2}/(\lambda-s)^{m-1}$ and a pole $s < \lambda_1,\theta_1$. Then,
\begin{equation}\label{eq.cjwjresex1}
\alpha_n(\theta_k)<|c_1|^2+\ldots+|c_k|^2<\alpha_n(\theta_{k+1}-),~~~~k=1,\dots,m-1.
\end{equation}
\end{proposition}
\begin{proof}
The proof of this proposition is similar to the proof of Theorem~\ref{thm.STCMS},
and is also provided in~\cite{Li98}.
We proceed with a sketch of the proof.

We first introduce $ \nu_j = \theta_j-s $ for $j=1,\ldots,m$.
The nodes $\nu_j$ are positive due to $s < \theta_1$
and for a fixed $k=1,\dots,m-1$
we let $\widehat{r}_{\pm,k}\in\Pi_{2m-2}/\lambda^{m-1}$ denote the rational functions
given in Proposition~\ref{prop.thatLaurentpol}.
Based on these rational functions,
we consider the rational functions $r_{\pm,k}(\lambda) = \widehat{r}_{\pm,k}(\lambda-s)$
in the class $\Pi_{2m-2}/(\lambda-s)^{m-1}$; and
based on properties of $\widehat{r}_{\pm,k}$ given in Proposition~\ref{prop.thatLaurentpol}
the functions $r_{\pm,k}$
satisfy
\begin{subequations}
\begin{equation}\label{eq.lpsmequal}
r_{\{\pm,k\}}(\theta_j)=\left\{
\begin{array}{ll}
1,&j=1,\ldots,k,\\
0,&j=k+1,\ldots,m,
\end{array}
\right.
\end{equation}
together with
\begin{equation}\label{eq.lpsminequal}
r_{\{+,k\}}(\lambda) \geq \left\{
\begin{array}{ll}
1,&\lambda \leq \theta_k,\\
0,& \lambda > \theta_k,
\end{array}
\right.
~~~\text{and}~~~
r_{\{-,k\}}(\lambda) \leq \left\{
\begin{array}{ll}
1,& \lambda < \theta_{k+1},\\
0,& \lambda \geq \theta_{k+1}.
\end{array}
\right.
\end{equation}
\end{subequations}

The identity~\eqref{eq.lpsmequal} implies
\begin{subequations}\label{eq.proofofSTCMSextended123}
\begin{equation}\label{eq.proofofSTCMSextended1}
\sum_{j=1}^m |c_j|^2 r_{\{\pm,k\}}(\theta_j) = |c_1|^2+\ldots+|c_k|^2,
\end{equation}
Analogously, the inequalities in~\eqref{eq.lpsminequal} imply
\begin{equation}\label{eq.proofofSTCMSextended2}
\alpha_n(\theta_k)
= \sum_{\{j:\lambda_j \leq \theta_{k}\}} |w_j|^2
< \sum_{j=1}^n |w_j|^2 r_{\{+,k\}}(\lambda_j),
\end{equation}
and
\begin{equation}\label{eq.proofofSTCMSextended3}
\alpha_n(\theta_{k+1}-)
= \sum_{\{j:\lambda_j < \theta_{k+1}\}} |w_j|^2
> \sum_{j=1}^n |w_j|^2 r_{\{-,k\}}(\lambda_j).
\end{equation}
The right-hand sides of~\eqref{eq.proofofSTCMSextended2}
and~\eqref{eq.proofofSTCMSextended3} can be understood as a Riemann-Stieltjes
integral~\eqref{eq.innerprodrepresentationA},
for which the quadrature property~\eqref{eq.idinrprodextendedkry} yields
\begin{equation}
\sum_{j=1}^n |w_j|^2 r_{\{\pm,k\}}(\lambda_j)
= \sum_{j=1}^m |c_j|^2 r_{\{\pm,k\}}(\theta_j),
\end{equation}
\end{subequations}
Combining the identities and inequalities in~\eqref{eq.proofofSTCMSextended123},
we conclude~\eqref{eq.cjwjresex1}; 
for further details we also refer to the proof of Theorem~\ref{thm.STCMS}.
\end{proof}
As previously discussed in Remark~\ref{rmk:STwithsj};
in the Krylov setting the measure $\alpha_n$ is not continuous
and a property as in~\cite[eq.~(3) in Theorem~1]{Li98} does not hold in general.

The spectrum of the Rayleigh quotient $A_m$
for the extended Krylov subspace given in~\eqref{eq.extendedKrysubspace}
defines a measure $\alpha_m$, as in~\eqref{eq.alpham}.
The result of Proposition~\ref{prop.STCMSex}
can be understood as an intertwining property
of the distributions $\dd\alpha_n$ and $\dd\alpha_m$,
similar as in the polynomial case in Subsection~\ref{subsec.bounds31}.
This property is illustrated for a numerical example
in Section~\ref{sec.experiments4} below.

Considering a rational qor-Krylov setting,
previous results for the qor-Krylov representation $B_m$
also apply to the extended Krylov subspace
(which does correspond to the rational Krylov subspace
$\RKry_{2\varrho-1}(A,u)$ with denominator $q(\lambda)=(\lambda-s)^{\varrho-1}$ as previously mentioned).
Thus, Proposition~\ref{prop.STCMSex} holds true for the qor-Krylov representation $B_m$
(assuming $s\leq a<\theta_1$).
However, these results are not specified here.




\section{Numerical illustrations}\label{sec.experiments4}

In the present section we verify the results of
Theorem~\ref{Kryproponweights} (Subsection~\ref{subsec.bounds3polSepTh})
and~\ref{thm.saixxx} (Subsection~\ref{subsec.SepTheoremRatKrylov}),
and Proposition~\ref{prop.SaIC1} (Subsection~\ref{subsec.CMSSaIcomplex})
and~\ref{prop.STCMSex} (Subsection~\ref{subsec.STCMSxKry})
by numerical experiments.

For the present numerical examples,
the notation $\theta_j$ and $c_j$ refers to the quadrature nodes and
weights, respectively, satisfying different polynomial and rational
Gaussian quadrature formulae which originate from polynomial Krylov subspaces $\Kry_m(A,u)$ 
and rational Krylov subspaces $\RKry_m(A,u)$ with different choices of poles.
Here, the matrix~$A\in\R^{n\times n}$ corresponds to
the finite-difference discretization
of the negative 1D Laplace operator with~$n=1200$,
and $u$ is a random starting vector which is normalized.
The \tMinr-inner product corresponds to the Euclidian inner product.

For the polynomial case, the quadrature nodes and weights are based
on the spectrum of the Jacobi matrix $J_m$ which is computed
using the Lanczos method.
Considering the rational case, we show results for SaI Krylov subspaces
with real and complex shifts. For the case of a real shift, we consider the
Rayleigh quotient $A_m$ as in Algorithm~\ref{alg.SaI}
and the rational qor-Krylov representation $B_m$ as in Algorithm~\ref{alg.qorKryrat}.
For a complex shift we show an example using the Rayleigh quotient $A_m$.
Furthermore, we show results for an extended Krylov subspace,
for which the Rayleigh quotient $A_m$ is computed using Algorithm~\ref{alg.extendedLanczos}.

\begin{figure}
\centering
\begin{overpic}
[width=0.8\textwidth]{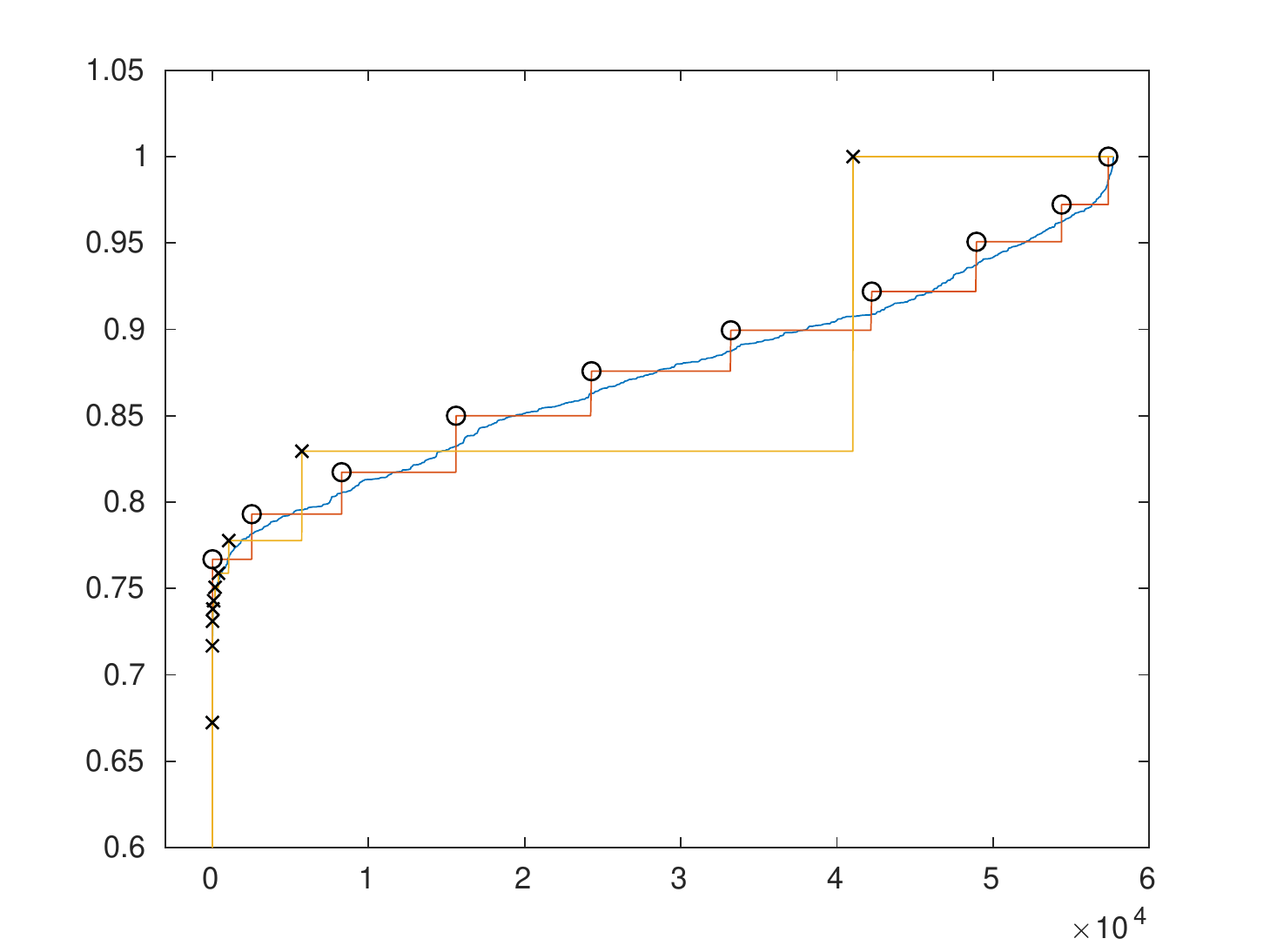}
\put(-1,38){$\alpha(\lambda)$}
\put(50,1){$\lambda$}
\put(31,8.8){\includegraphics[width=0.3\textwidth]{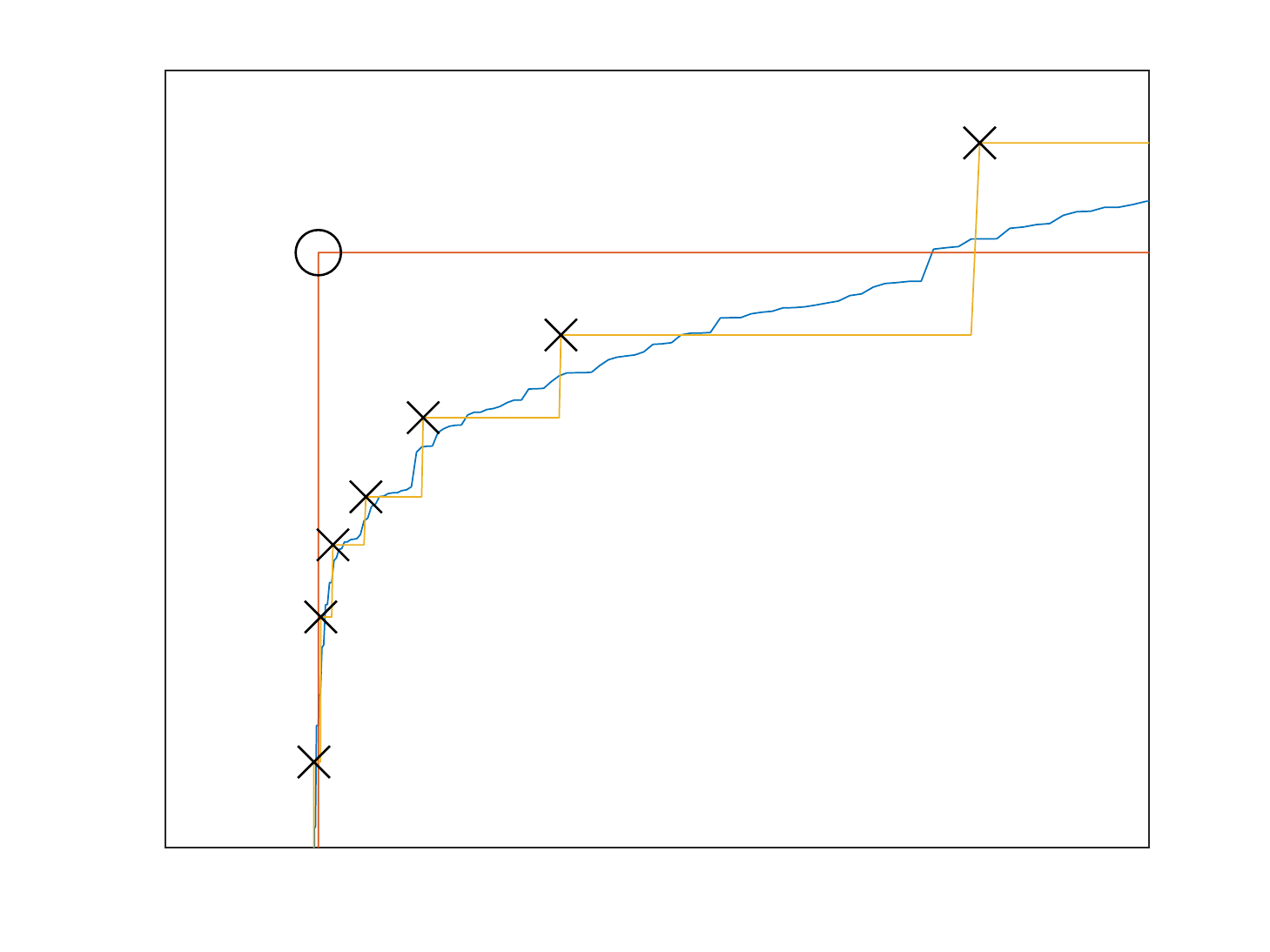}}
\put(33,25){\vector(-3,1){12}}
\end{overpic}
\caption{The matrix~$A\in\R^{n\times n}$ is given by
the finite-difference discretization
of the negative 1D Laplace operator with~$n=1200$,
and $u$ is a random starting vector which is normalized.
The continuous line without additional symbols illustrates the
step function~$\alpha_n$ associated with the eigenvalues
and spectral coefficients of~$u$ in the eigenbasis of~$A$.
The symbols ('$\circ$') mark~$\alpha_m(\theta_j)$
where $\theta_j$ are the Ritz values of the polynomial
Krylov subspace~$\Kry_m(A,u)$ with~$m=10$,
and~$\alpha_m$ is the respective step function given in~\eqref{eq.alpham}.
Similarly, the symbols ('$\times$') mark~$\alpha_m(\theta_j)$
where $\theta_j$ refer to the rational Krylov subspace~$\RKry_m(A,u)$
with~$m=10$ and a single pole~$s=-10^2$ of multiplicity $m-1$.
}
\label{fig.stepfcts}
\end{figure}

\begin{figure}
\centering
\begin{tabular}{cc}
\hspace{0.5cm}
\begin{overpic}
[width=0.5\textwidth]{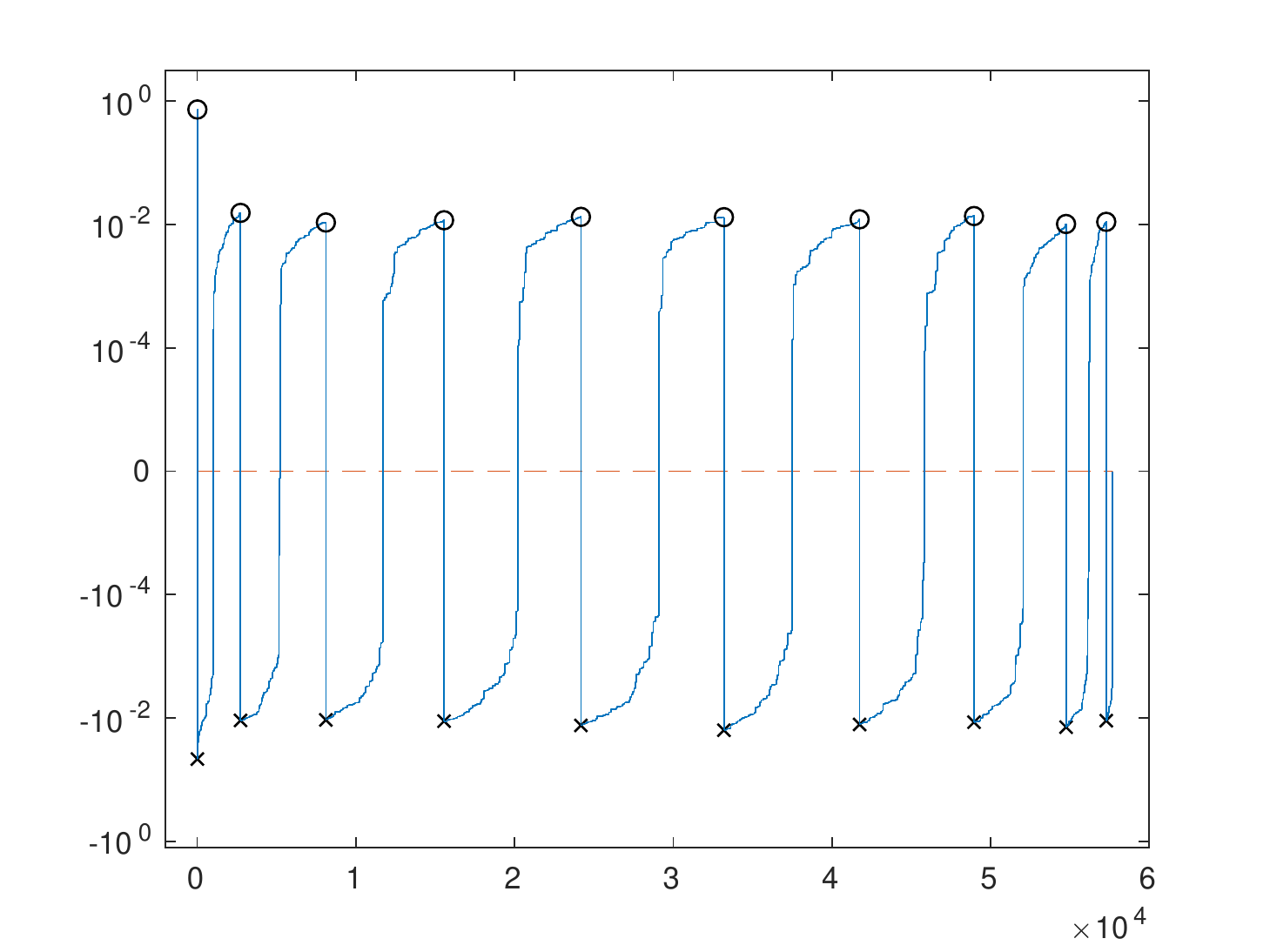}
\put(-6,38){~\small$F(\lambda)$}
\put(-3,65){ \textbf{a)}}
\put(50,-1){\small$\lambda$}
\end{overpic}
&
\hspace{-0.9cm}
\begin{overpic}
[width=0.5\textwidth]{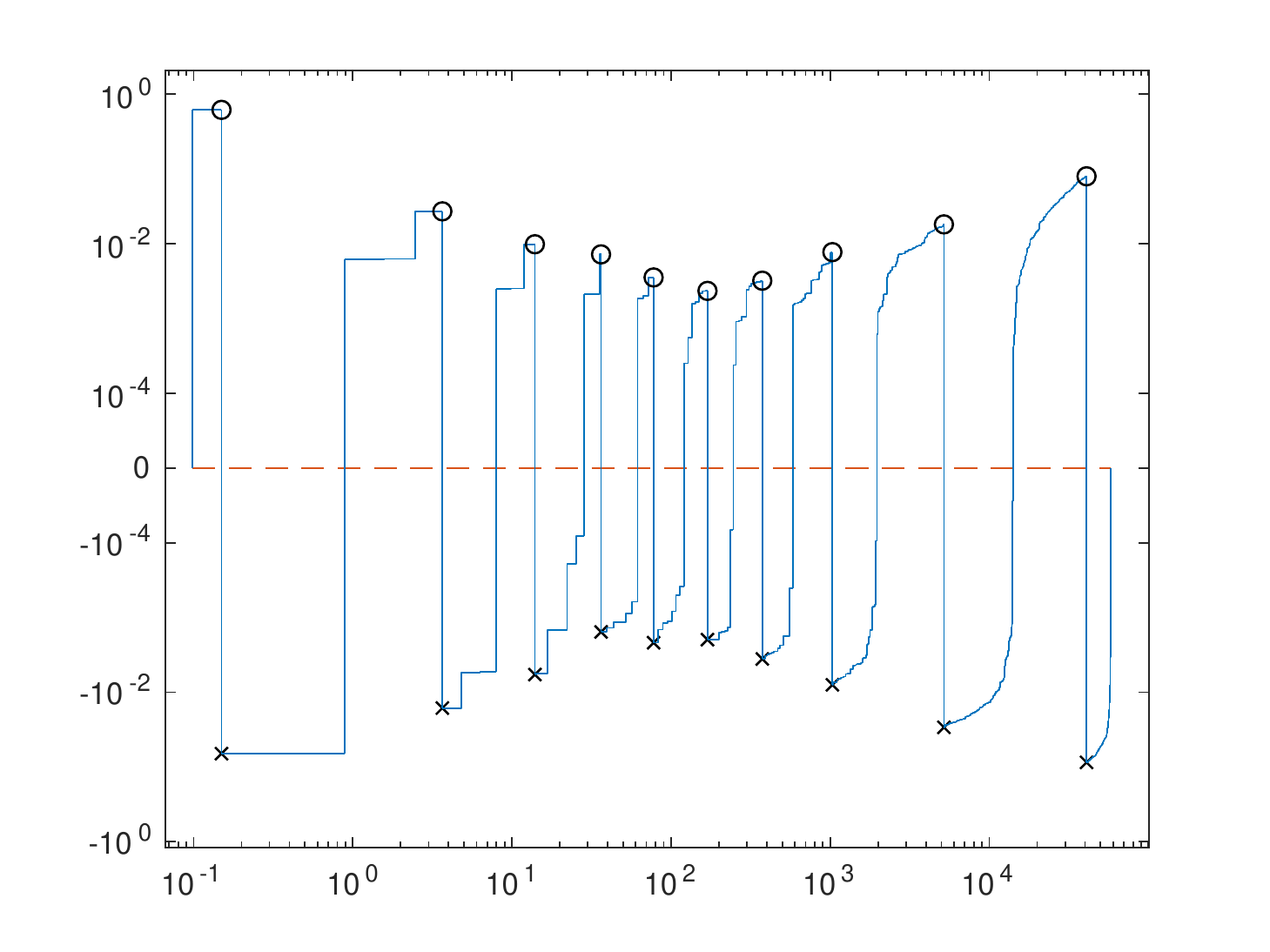}
\put(0,65){ \textbf{b)}}
\put(50,-1){\small$\lambda$}
\end{overpic}
\end{tabular}
\caption{In Figure~\textbf{a)} and~\textbf{b)} the matrix~$A\in\R^{n\times n}$ is
given by the finite-difference discretization of the negative 1D Laplace operator
with~$n=1200$.
The starting vector~$u$ is chosen at random and is normalized.
In these figures
we show the function $F=\alpha_n-\alpha_m$ where $\alpha_m$
originates from different settings as stated below.
The symbols~('$\circ$') and~('$\times$')
mark $F(\theta_k-)$ and $F(\theta_k)$, respectively.
\newline\noindent--
Figure~\textbf{a)}~shows $F$ with $\alpha_m$ given
by spectral weights and Ritz values of the Jacobi matrix $J_m$
for the polynomial Krylov subspace~$\Kry_m(A,u)$ with~$m=10$.
The $y$-axis is scaled logarithmically in positive and negative direction,
namely, with range $(-10^0,-10^{-6})\cup(10^{-6},10^0)$.
\newline\noindent--
Figure~\textbf{b)} shows $F$ where~$\alpha_m$
refers to the spectrum of the Rayleigh quotient $A_m$
for the rational Krylov subspace~$\RKry_m(A,u)$ with~$m=10$
and a single pole~$s=-10^2$ of multiplicity~$m-1$, thus, $s<\lambda_1$.
Similar to Figure~a) the $y$-axis is scaled logarithmically and
covers $(-10^0,-10^{-5})\cup(10^{-5},10^0)$.
Additionally, the $x$-axis is scaled logarithmically in a classical sense.
}
\label{fig:Fplots1}
\end{figure}

\medskip
The step functions $\alpha_n$~\eqref{eq.measureAv} and $\alpha_m$~\eqref{eq.alpham}
are illustrated for numerical examples in \figref{fig.stepfcts}.
The step function $\alpha_m$ is shown for the polynomial Krylov subspace
and a SaI Krylov subspace with a shift $s\in\R$ located outside of the convex hull of the matrix spectrum,
namely, $s<\lambda_1$.
In both cases the distributions $\dd\alpha_n$ and $\dd\alpha_m$ satisfy an intertwining property.
To provide a clear illustration of the results of the previous section we
also show the function $F(\lambda)$ given in~\eqref{eq:Fxalpha}
for a numerical example concerning the polynomial case in \figref{fig:Fplots1} a).
In this figure, we observe that $F(\lambda)$ changes its sign at the Ritz values,
and following~\eqref{eq.Flowerbound}, this verifies the result of Theorem~\ref{Kryproponweights}.
For the SaI Krylov subspace with a shift $s<\lambda_1$
we have $F(s)=0$ which implies $F_s(\lambda)=F(\lambda)$
for the function $F_s(\lambda)$ as given in~\eqref{eq.verifySaI2Fs}.
Considering this example, the function $F=F_s$ is illustrated in \figref{fig:Fplots1} b),
and following Remark~\ref{rmk.FsandThmsaix}, the change of the sign of $F_s$ at rational Ritz values
verifies Theorem~\ref{thm.saixxx}.

\begin{figure}
\centering
\begin{tabular}{cc}
\hspace{0.5cm}
\begin{overpic}
[width=0.5\textwidth]{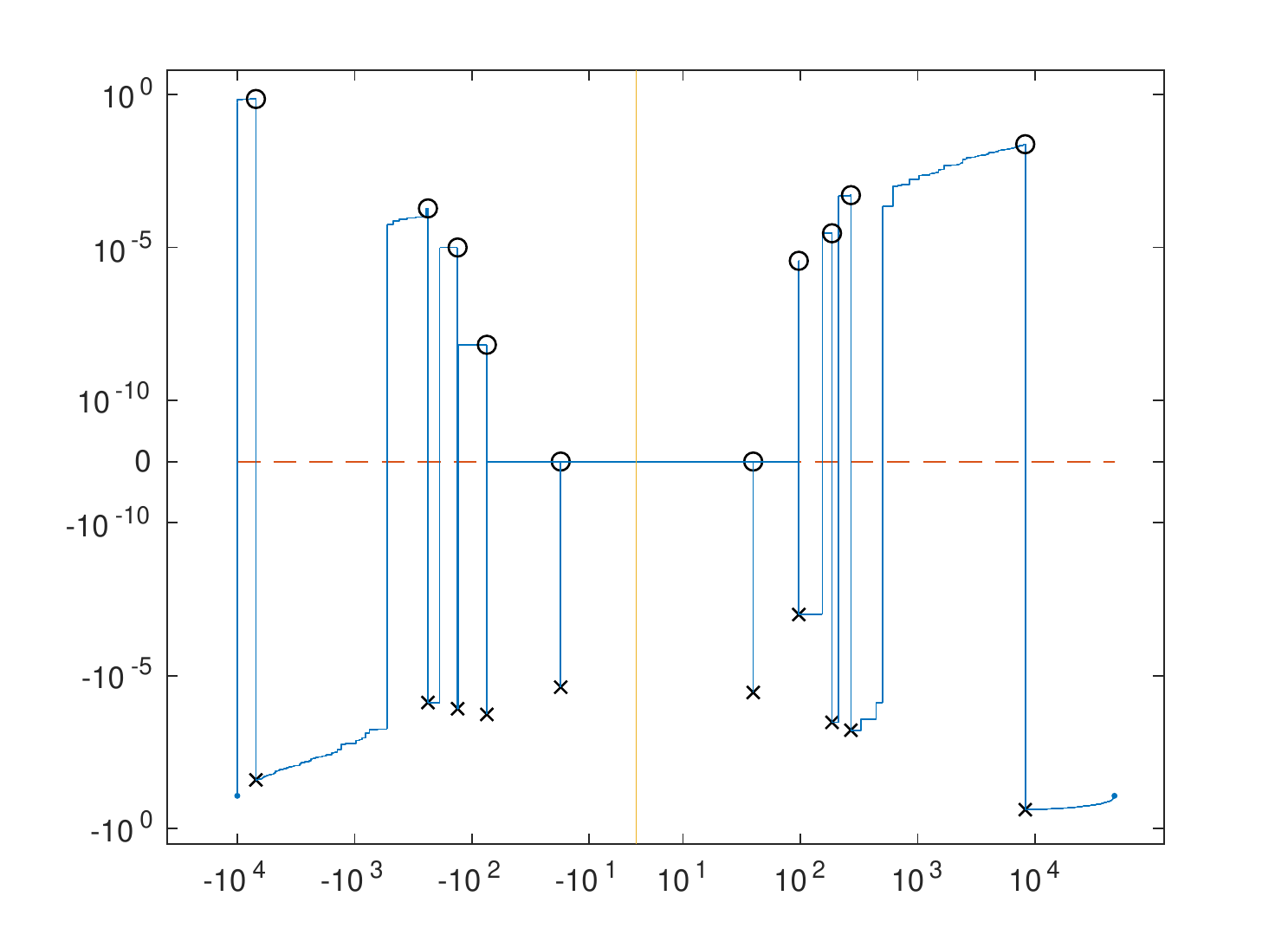}
\put(-7,36){\small$F_s(\lambda)$}
\put(51,63){\small$s=10^4$}
\put(44,-1){\small$\lambda-s$}
\put(0,65){ \textbf{a)}}
\end{overpic}
\hspace{-0.9cm}
&
\begin{overpic}
[width=0.5\textwidth]{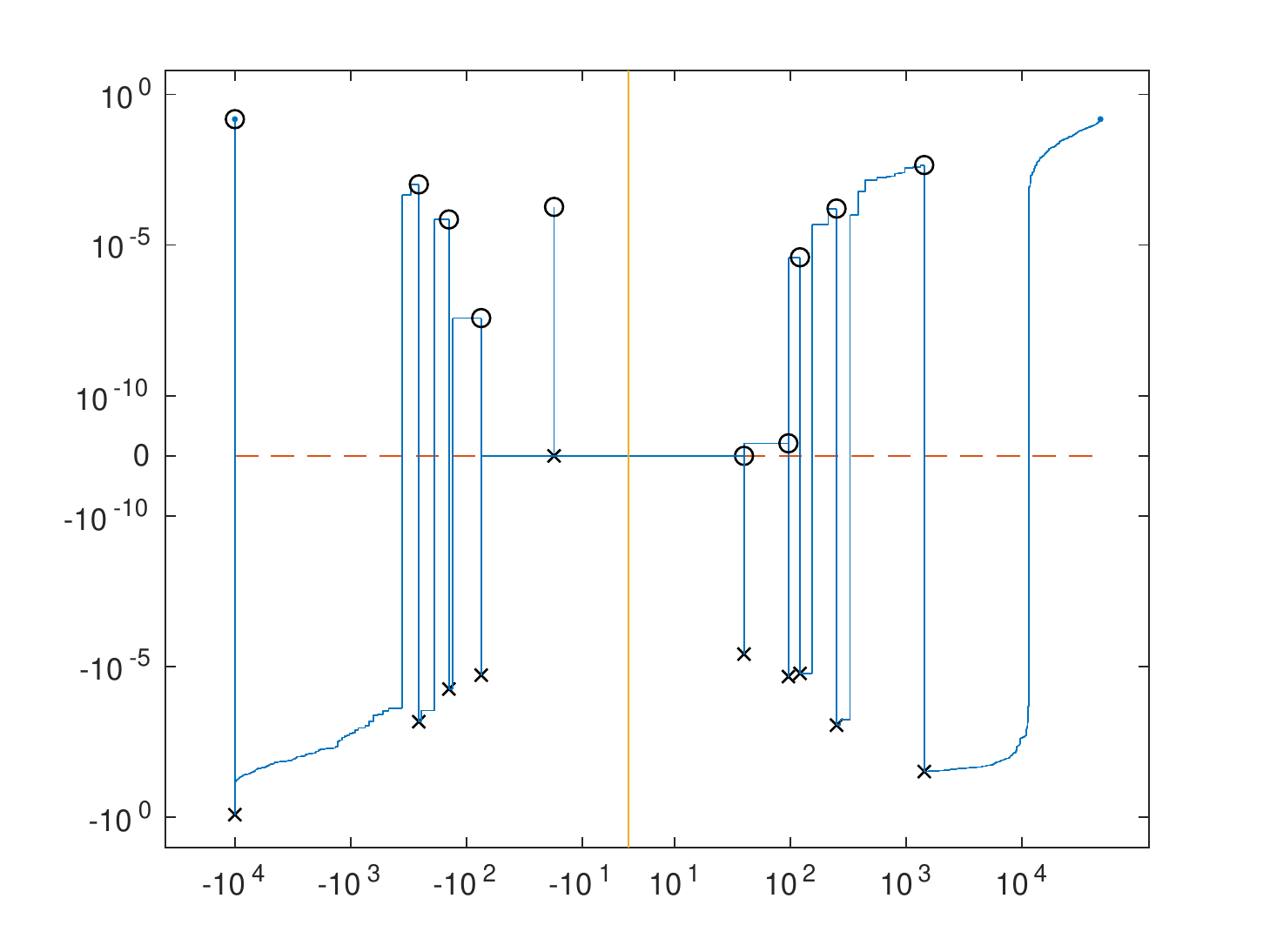}
\put(51,64){\small$s=10^4$}
\put(44,-1){\small$\lambda-s$}
\put(0,65){ \textbf{b)}}
\end{overpic}
\end{tabular}
\caption{In Figure~\textbf{a)} and~\textbf{b)} the matrix~$A\in\R^{n\times n}$ is
given by the finite-difference discretization of the negative 1D Laplace operator
with~$n=1200$. 
The starting vector~$u$ is chosen at random and is normalized.
These figures show $F_s(\lambda)=F(\lambda)-F(s)$
for different settings, and the symbols~('$\circ$')
and~('$\times$') mark $F_s(\theta_k-)$ and $F_s(\theta_k)$, respectively.
Similar to \figref{fig:Fplots1} the
$y$-axis is scaled logarithmically and
covers~$(-10^0,-10^{-12})\cup(10^{-12},10^0)$.
With $\lambda$ being the argument of the function $F_s(\lambda)$
as illustrated in the $y$-axis, the $x$-axis shows~$\lambda-s$,~i.e.~the
distance from the argument $\lambda$ to the pole $s=10^4$.
Furthermore, the $x$-axis is scaled logarithmically with a range of
approximately~$(-10^4,-10^{1})\cup(10^{1},10^4)$.
\newline\noindent--
Figure \textbf{a)} shows $F_s(\lambda)$
corresponding to the spectrum of $A_m$,
where $A_m$ is the Rayleigh quotient in the
rational Krylov subspace~$\RKry_m(A,u)$ with~$m=10$
and a single pole~$s=10^4$ of multiplicity~$m-1$.
Here, the pole $s$ is enclosed by the eigenvalues of~$A_m$.
\newline\noindent--
Figure \textbf{b)} shows $F_s(\lambda)$ where $\alpha_m$ corresponds
to the spectrum of $B_m$, which is the rational qor-Krylov representation
for which the eigenvalue~$\theta_1=-10$ is preassigned.
For the underlying rational Krylov subspace~$\RKry_m(A,u)$ we have~$m=10$
and a single pole~$s=10^4$ of multiplicity~$m-1$.
The pole $s$ is enclosed by the eigenvalues of~$B_m$.
}
\label{fig:Fplots2}
\end{figure}

\medskip
The case of a SaI Krylov subspace with a shift 
$s\in\R$ such that $\theta_1<s<\theta_m$ is illustrated in \figref{fig:Fplots2}.
As for the previous example, the change of the sign of $F_s$ at rational Ritz values
verifies Theorem~\ref{thm.saixxx} as stated in Remark~\ref{rmk.FsandThmsaix}.
Here, \figref{fig:Fplots2}~a) illustrates~$F_s$ for the Rayleigh quotient $A_m$ and
\figref{fig:Fplots2}~b) illustrates~$F_s$ for a rational qor-Krylov
representation with a preassigned eigenvalue $\xi\in\R$;
this verifies the result of Theorem~\ref{thm.saixxx} for these cases.

\begin{figure}
\centering
\begin{overpic}
[width=0.8\textwidth]{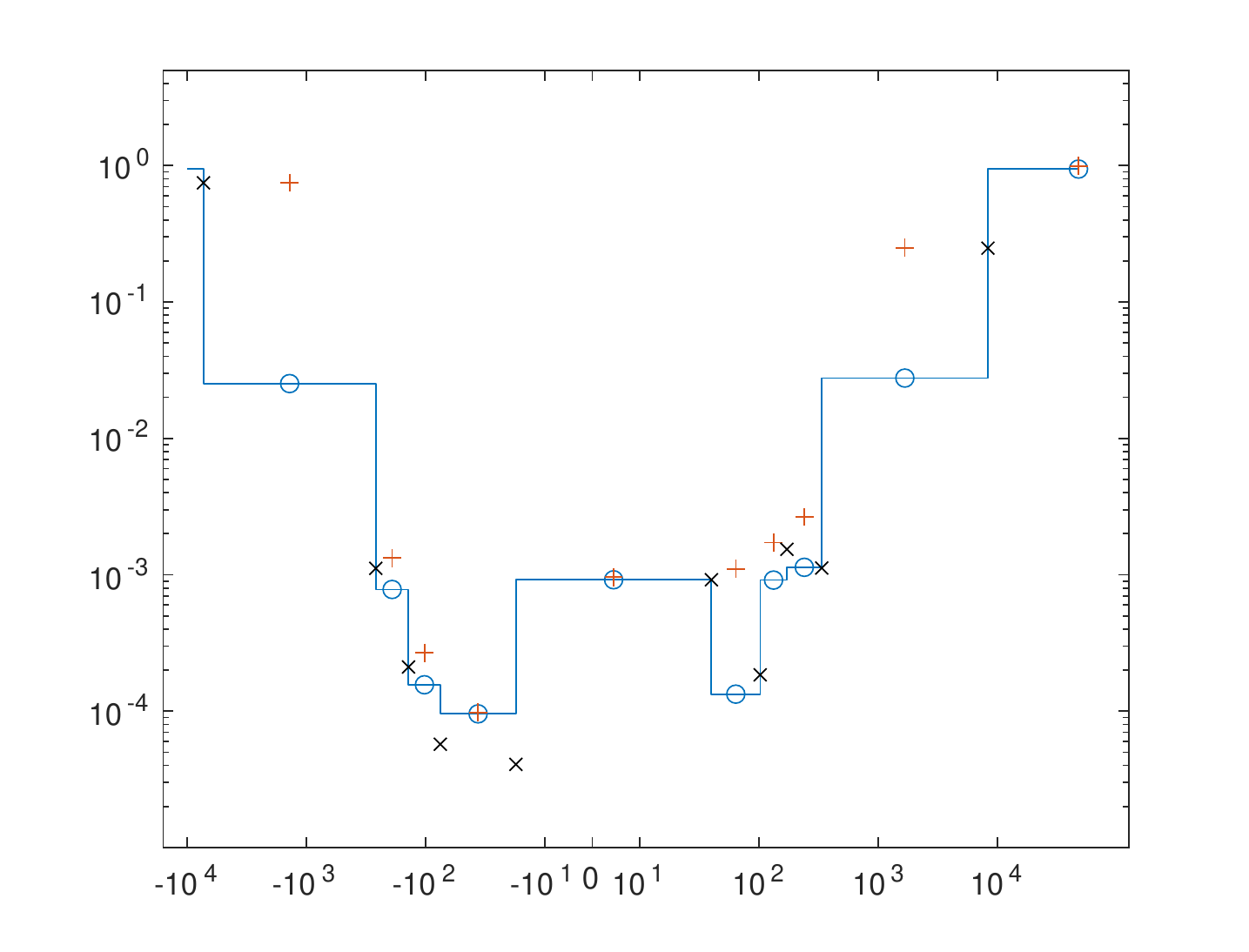}
\put(50,-1){\small$\lambda-\real s$}
\end{overpic}
\caption{The matrix~$A\in\R^{n\times n}$ is
given by the finite-difference discretization of the negative 1D Laplace operator
with~$n=1200$.
The starting vector~$u$ is chosen at random and is normalized.
For the present figure we consider the SaI Krylov subspace $\RKry_m(A,u)$ with $m=10$
and a complex shift~$s= 10^4 -10^2\ii$.
In the present caption, $c_j$ and $\theta_k$ refer to the entries of eigenvectors and the eigenvalues of
the respective Rayleigh quotient $A_m$.
The symbols~('$\times$') mark~$|c_j|^2$ over~$\theta_j$.
The symbols~('$+$') show $|c_k|^2+|c_{k+1}|^2$ over the midpoint of the interval $[\theta_k,\theta_{k+1}]$ for $k=1,\ldots,m-1$.
Furthermore, the symbol~('$+$') located at the right boundary of the spectrum shows $|c_m|^2+|c_1|^2$.
The line marked by~('$\circ$') shows the measure $\mu_n([\theta_k,\theta_{k+1}])$
over each interval $[\theta_k,\theta_{k+1}]$ for $k=1,\ldots,m-1$,
and the measure $\mu_n((-\infty,\theta_1]\cup[\theta_m,\infty))$ at the boundary.
The $y$-axis is scaled logarithmically in a classical sense,
and the $x$-axis shows~$\lambda-\real s$,~i.e.~the
distance from the argument $\lambda$ to the real part of the shift,~i.e.,~$\real s=10^4$.
Furthermore, the $x$-axis is scaled logarithmically with a range of
approximately~$(-10^4,-10^{1})\cup(10^{1},10^4)$.
}
\label{fig:Fplots3complexSaI}
\end{figure}

\medskip
In \figref{fig:Fplots3complexSaI} we consider a SaI Krylov subspace
with a complex shift $s\in\C\setminus\R$.
For this example, we illustrate $|c_k|^2+|c_{k+1}|^2$ for $k=1,\ldots,m-1$
and $|c_m|^2+|c_1|^2$, which yield upper bounds on $\mu_n([\theta_k,\theta_{k+1}])$ for $k=1,\ldots,m-1$
and $\mu_n((-\infty,\theta_1]\cup[\theta_m,\infty))$, respectively.
This verifies the result of Proposition~\ref{prop.SaIC1}.

\begin{figure}
\centering
\begin{overpic}
[width=0.8\textwidth]{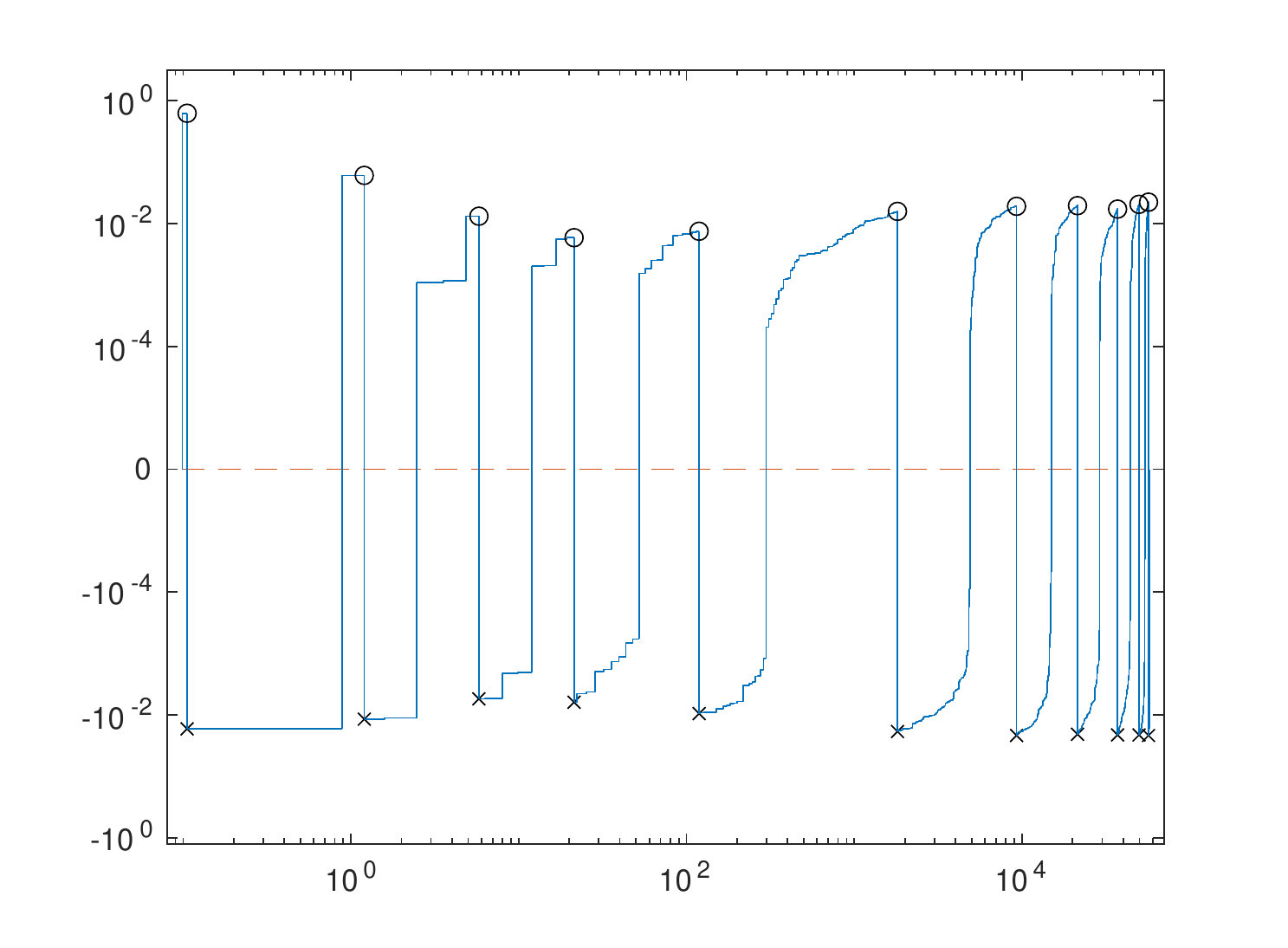}
\put(-6,38){~\small$F(\lambda)$}
\put(50,-1){\small$\lambda$}
\end{overpic}
\caption{The matrix~$A\in\R^{n\times n}$ is
given by the finite-difference discretization of the negative 1D Laplace operator
with~$n=1200$.
The starting vector~$u$ is chosen at random and is normalized.
In this figure
we show the function $F=\alpha_n-\alpha_m$ where $\alpha_m$
corresponds to the spectrum of the Rayleigh quotient $A_m$
given by the extended Krylov subspace~\eqref{eq.extendedKrysubspace}
with $m=11$ (thus, $\varrho=6$) and the shift~$s=-10<\lambda_1$.
The symbols~('$\circ$') and~('$\times$')
mark $F(\theta_k-)$ and $F(\theta_k)$, respectively,
where $\theta_k$ refers to the eigenvalues of $A_m$.
The $y$-axis is scaled logarithmically in positive and negative direction,
namely, with range $(-10^0,-10^{-6})\cup(10^{-6},10^0)$.
Additionally, the $x$-axis is scaled logarithmically in a classical sense.
}
\label{fig:Fplots3exKry}
\end{figure}

\medskip
For the extended Krylov subspace as in Subsection~\ref{subsec.STCMSxKry},
Proposition~\ref{prop.STCMSex} yields an intertwining property
for the distributions $\dd\alpha_n$ and $\dd\alpha_m$ as in the polynomial case;
the changing sign of $F$ as illustrated in \figref{fig:Fplots3exKry}
verifies the result of Proposition~\ref{prop.STCMSex}.


\ifthesis
\section*{Appendix}
\begin{subappendices}
\else
\appendix
\fi

\section{Some properties of Krylov subspaces}\label{appendix.A2}

\begin{proposition}\label{rmk.wlamrank}
Let $(q_1,\ldots,q_n)\in\C^{n\times n}$ be an orthogonal eigenbasis
of the matrix~$A\in\C^{n\times n}$.
Here, orthogonal is to be understood w.r.t.~a given
positive definite inner product.
Let $\lambda_1,\ldots,\lambda_n\in\C$ be the corresponding eigenvalues, and
$w_j=(q_j,u)_{\Minr}\in\C$ be the spectral coefficients of a given vector $u$. Then
\begin{equation}\label{eq.polAfullrankX}
\rank\{u,Au,\ldots,A^{m-1}u \} = m,
\end{equation}
if and only if there exist at least~$m$ coefficients~$w_j\neq 0$
     with distinct~$\lambda_j$.
\end{proposition}
\begin{proof}
According to the eigendecomposition of~$A$ we have
$$
A^\ell u = \sum_{j=1}^n \lambda_j^\ell w_j q_j~~~\text{for $\ell\in\N_0$}.
$$
The matrix corresponding to the left-hand side of~\eqref{eq.polAfullrankX} takes the form of
a Vandermonde matrix,
\begin{equation}\label{eq.KryrankbybasisandVandermonde}
(u,Au,\ldots,A^{m-1}u) = (q_1 w_1,q_2w_2,\ldots,q_nw_n)
\begin{pmatrix}
1&\lambda_1&\lambda_1^2&\cdots&\lambda_1^{m-1}\\
1&\lambda_2&\lambda_2^2&\cdots&\lambda_2^{m-1}\\
\vdots&\vdots&\vdots&&\vdots\\
1&\lambda_n&\lambda_n^2&\cdots&\lambda_n^{m-1}
\end{pmatrix}\in\C^{n\times m}.
\end{equation}
Let $n_1\leq n$ be the number of nonzero coefficients~$w_j$,
thus, there exist indices $j(1),\ldots,j(n_1)$
with $w_{j(1)},\ldots,w_{j(n_1)}\neq 0$.
We define
$$
\Theta_1 = \big(q_{j(1)} w_{j(1)},q_{j(2)}w_{j(2)},\ldots,q_{j(n_1)}w_{j(n_1)}\big)\in\C^{n\times n_1}.
$$
The orthogonality properties of~$q_1,\ldots,q_n$ imply $\rank(\Theta_1)=n_1$.
For the corresponding rows of the Vandermonde matrix we introduce the notation
$$
\Theta_2 = \begin{pmatrix}
1&\lambda_{j(1)}&\lambda_{j(1)}^2&\cdots&\lambda_{j(1)}^{m-1}\\
1&\lambda_{j(2)}&\lambda_{j(2)}^2&\cdots&\lambda_{j(2)}^{m-1}\\
\vdots&\vdots&\vdots&&\vdots\\
1&\lambda_{j(n_1)}&\lambda_{j(n_1)}^2&\cdots&\lambda_{j(n_1)}^{m-1}
\end{pmatrix}\in\C^{n_1\times m}.
$$
The identity in~\eqref{eq.KryrankbybasisandVandermonde} can now be written as
\begin{equation}\label{eq.Kryrankproofeq1}
\big(u,Au,\ldots,A^{m-1}u\big) = \Theta_1 \Theta_2.
\end{equation}
With $\Theta_1\in\C^{n\times n_1}$ and $\rank(\Theta_1)=n_1$ we have $\rank(\Theta_1 \Theta_2)=\rank(\Theta_2)$.
Let $n_2\leq n_1$ be the number of distinct eigenvalues within $\lambda_{j(1)},\ldots,\lambda_{j(n_1)}$,
hence,~we have indices $\ell(1),\ldots,\ell(n_2)$ for which $\lambda_{j(\ell(1))},\ldots,\lambda_{j(\ell(n_2))}$ are distinct.
~Then the Vandermonde matrix $\Theta_2$ satisfies~$\rank(\Theta_2)=\min \{m,n_2\}$, hence,
\begin{equation}\label{eq.Kryrankproofeq2}
\rank(\Theta_1\Theta_2)=\rank(\Theta_2)=\min \{m,n_2\}.
\end{equation}
Combining~\eqref{eq.Kryrankproofeq1} with~\eqref{eq.Kryrankproofeq2} we conclude
\begin{equation*}
\rank\{u,Au,\ldots,A^{m-1}u \} = \min \{m,n_2\}.
\end{equation*}
We recall that $n_2$ is number of nonzero coefficients~$w_j$ with distinct $\lambda_j$,
and~\eqref{eq.polAfullrankX} holds if and only if $ n_2\geq m $ which completes the proof.
\end{proof}


\begin{proposition}\label{prop.app3.cjneqzero}
Let $w_1,\ldots,w_n\in\C$ with $w_j\neq 0$ and $\lambda_1<\ldots<\lambda_n$ be given.
Let $m<n$, and let $\theta_1<\ldots<\theta_m$ and $|c_1|^2,\ldots,|c_m|^2$ denote
quadrature nodes and quadrature weights, respectively,
and assume
\begin{equation}\label{eq.app3.idforpol}
\sum_{j=1}^n |w_j|^2 p(\lambda_j)=\sum_{j=1}^m |c_j|^2 p(\theta_j),~~~p\in\Pi_{2m-2}.
\end{equation}
Then $c_j\neq 0$ for $j=1,\ldots, m$.
\end{proposition}
\begin{proof}
We define $g_\ell(\lambda)=\prod_{j=1,j\neq \ell}^m (\lambda-\theta_j)^2 \in\Pi_{2m-2}$.
The polynomial $g_\ell$ is zero only at the nodes
$\theta_1,\ldots,\theta_{\ell-1},\theta_{\ell+1},\ldots,m$
and positive otherwise.
Due to $n>m$ at least one $\lambda_j$ is distinct to $\theta_1,\ldots,\theta_m$
and this yields
$$
\sum_{j=1}^n |w_j|^2 g_\ell(\lambda_j) >0.
$$
Making use of the identity~\eqref{eq.app3.idforpol} and
evaluating the right-hand side therein we conclude
$$
\sum_{j=1}^m |c_j|^2 g_\ell(\theta_j) = g_\ell(\theta_{\ell}) |c_\ell|^2 >0.
$$ 
With $g_\ell(\theta_{\ell})>0$ this concludes $|c_\ell|^2>0$.
\end{proof}

\begin{proposition}[\em Identities for rational functions in the rational Krylov subspace]
\label{prop.equalitiesratKry}

Let $U_m\in\C^{n\times m}$ with~$(U_m,U_m)_{\Minr}=I$
and~$\vspan\{U_m\}=\RKry_m(A,u)$ for the rational Krylov subspace
$\RKry_m(A,u)$ with denominator $q_{m-1}$.
Let~$A_m=(U_m,A\,U_m)_{\Minr}$ and $x = (U_m,u)_{\Minr}$.
\begin{enumerate}[label=(\roman*)]
\item\label{item.identratprop}
The following identities hold true,
\begin{equation*}
r(A)u = U_m\,r(A_m)\, x,~~~r\in\Pi_{m-1}/q_{m-1}.
\end{equation*}
\item\label{item.residualratprop}
Let~$r=p/q_{m-1}$ with~$p\in \Pi_{m}$ being a polynomial
of degree exactly~$m$, then
\begin{equation}\label{eq.pmqperpRKmA}
( U_m\,r(A_m)x - r(A)u) \perp_{\Minr} \vspan\{U_m\}.
\end{equation}
\end{enumerate}

\end{proposition}
\begin{proof}
We proceed similar to the proof of Proposition~\ref{prop.idforqorratKry}
in Subsection~\ref{sec.ratqor}.
We recall the identity $\RKry_m(A,u)=\Kry_m(A,u_q)$
with~$u_q=q_{m-1}^{-1}(A)u$.
Let $\zeta_0=\|u_q\|_{\Minr} $, let $V_m$ be the \tMinr-orthonormal basis of $\Kry_m(A,u_q)$,
and let $J_m$ be the respective Jacobi matrix.
Then the identity~\eqref{eq.pofAequalKm} w.r.t.\ $\Kry_m(A,u_q)$ implies
\begin{equation}\label{eq.ratkryqoridp1x}
p(A)u_q =\zeta_0\,V_m\,p(J_m) e_1,~~~p\in \Pi_{m-1}.
\end{equation}
This implies $q_{m-1}(A)u_q=\zeta_0\,V_m\,q_{m-1}(J_m) e_1$,
and with the identities $q_{m-1}(A) u_q = u $
and $(V_m, V_m)_{\Minr} = I$ we arrive at 
\begin{equation}\label{eq.ratkryqoridp2x}
\zeta_0\,e_1 = q_{m-1}^{-1}(J_m) (V_m,u)_{\Minr}.
\end{equation}
Let $r=p/q_{m-1}$ with $p\in\Pi_{m-1}$ then $r(A)u=p(A)u_q$, and with~\eqref{eq.ratkryqoridp1x} we have
\begin{equation}\label{eq.ratkryqoridp3x}
r(A)u = \zeta_0\,V_m\,p(J_m) e_1.
\end{equation}
Inserting~\eqref{eq.ratkryqoridp2x} into~\eqref{eq.ratkryqoridp3x} gives
\begin{equation}\label{eq.ratkryqoridpxA}
r(A)u = V_m\,p(J_m) q_{m-1}^{-1}(J_m) (V_m,u)_{\Minr}
= V_m\,r(J_m) (V_m,u)_{\Minr}.
\end{equation}
With the identity~$ K_m\,K_m^\Hast = I $ (see~\eqref{eq.Kmortho1})
and~\eqref{eq.RKAmtoJm} the matrix $A_m$
satisfies $r(J_m)=K_m\,r(A_m) K_m^\Hast$,
and together with~$V_m K_m = U_m$~\eqref{eq.Kmortho1} we have
\begin{equation}\label{eq.ratkryqoridp4x}
V_m\,r(J_m) (V_m,u)_{\Minr} = U_m\,r(A_m) (U_m,u)_{\Minr}.
\end{equation}
Combining~\eqref{eq.ratkryqoridpxA} with~\eqref{eq.ratkryqoridp4x}
completes the proof of~\ref{item.identratprop}.

For a polynomial $p$ of degree exactly~$m$
and w.r.t.~$\Kry_m(A,u_q)$ the property~\eqref{eq.pmperpKm} writes
\begin{equation}\label{eq.ratkryresp2residualofKryuq}
p(A) u_q - \zeta_0 V_m\,p(J_m) e_1
\perp_{\Minr} \vspan\{V_m\},~~~p\in\Pi_m .
\end{equation}
Let $r=p/q_{m-1}$,
then the identities in~\eqref{eq.ratkryqoridp2x}
and~\eqref{eq.ratkryqoridp4x} with $x = (U_m,u)_{\Minr}$ entail
$$
\zeta_0 V_m\,p(J_m) e_1 = V_m\,r(J_m) (V_m,u)_{\Minr}
 = U_m\,r(A_m) x.
$$
With $ r(A)u = p(A) u_q $ this yields
\begin{equation}\label{eq.ratkryresp2}
p(A) u_q - \zeta_0 V_m\,p(J_m) e_1
= r(A)u - U_m\,r(A_m) x.
\end{equation}
Making use of $\vspan\{V_m\}=\vspan\{U_m\}$ in~\eqref{eq.ratkryresp2residualofKryuq}
and substituting~\eqref{eq.ratkryresp2},
we conclude~\eqref{eq.pmqperpRKmA}.
\end{proof}

\begin{proposition}[\em The spectral coefficients $c_j$ for the rational
Krylov subspace and the choice of $U_m$]
\label{prop.cjindependentofUm}
The spectral coefficients $|c_j|$ of $x=(U_m,u)_{\Minr}\in\C^m$ in the
orthonormal eigenbasis of $A_m=(U_m,A\,U_m)_{\Minr}\in\C^{m\times m}$
are independent of the explicit choice of the underlying
orthonormal basis $U_m$ of~$\RKry_m(A,u)$.
\end{proposition}
\begin{proof}
We recall the representation of the spectral coefficients $c_j$
given in~\eqref{eq.cjforAm},
$$
c_j=(\widehat{q}_j, x)_2,~~~j=1,\ldots,m.
$$
Here $\widehat{q}_j\in\C^m$ refer to the
orthonormal eigenvectors of $A_m$.
We further recall that the rational Krylov subspace $\RKry_m(A,u)$
corresponds to the polynomial Krylov subspace $\Kry_m(A,u_q)$
with $u_q=q_{m-1}^{-1}(A)u$ for the denominator $q_{m-1}$.
Let us recall the notation $J_m$ and $V_m$ for the Jacobi matrix and
Krylov basis of~$\Kry_m(A,u_q)$.
Furthermore, we recall the orthonormal transformation
$K_m = (V_m,U_m)_{\Minr}\in\C^{m\times m}$ given in~\eqref{eq.defKm}.
With $U_m=V_mK_m$~\eqref{eq.UmtoVmKm}
the vector $x=(U_m,u)_{\Minr}$ corresponds to
\begin{equation}\label{eq.cjbyAmisJm}
x = K_m^\Hast (V_m,u)_{\Minr} =:  K_m^\Hast \xi,~~~
\text{and thus,}~~~c_j=(K_m \widehat{q}_j, \xi)_2.
\end{equation}
With the identity $ A_m = K_m^\Hast\,J_m\,K_m $~\eqref{eq.RKAmtoJm}
and $\widehat{q}_j$ being eigenvectors of $A_m$,
the vectors $K_m \widehat{q}_j$ for $j=1,\ldots,m$ correspond
to the orthonormal eigenvectors of $J_m$ (up to a complex phase).
Thus,~\eqref{eq.cjbyAmisJm} implies that the coefficients $c_j$
correspond to spectral coefficients of~$\xi=(V_m,u)_{\Minr}$ in the
orthonormal eigenbasis of $J_m$, and furthermore,
the coefficients~$c_j$ are independent of the explicit choice of~$U_m$.
More precisely, this applies to the absolute value~$|c_j|$ due to a potential
complex phase on orthonormal eigenvectors.
\end{proof}

\section{Auxiliary functions for the CMS Theorem}\label{appendix.A4}

\begin{proof}[Proof of \,\textbf{Proposition~\ref{prop.thatpol1}}]
We recapitulate arguments of~\cite{Sze85,Ak65} and others.
Let $\theta_1<\ldots<\theta_m$ and $k\in\{1,\ldots,m-1\}$ be given.
We first prove the existence of a polynomial $p_{\{+,k\}}$
of degree $2m-2$ which satisfies~\eqref{eq.ppmequal}
and~\eqref{eq.ppminequal}.
Let~$p$ be a polynomial of degree $2m-2$ subject
to the conditions
\vspace{.3cm}
\newline\noindent
\begin{minipage}{0.475\textwidth}
\centering
\begin{overpic}
[width=0.9\textwidth]{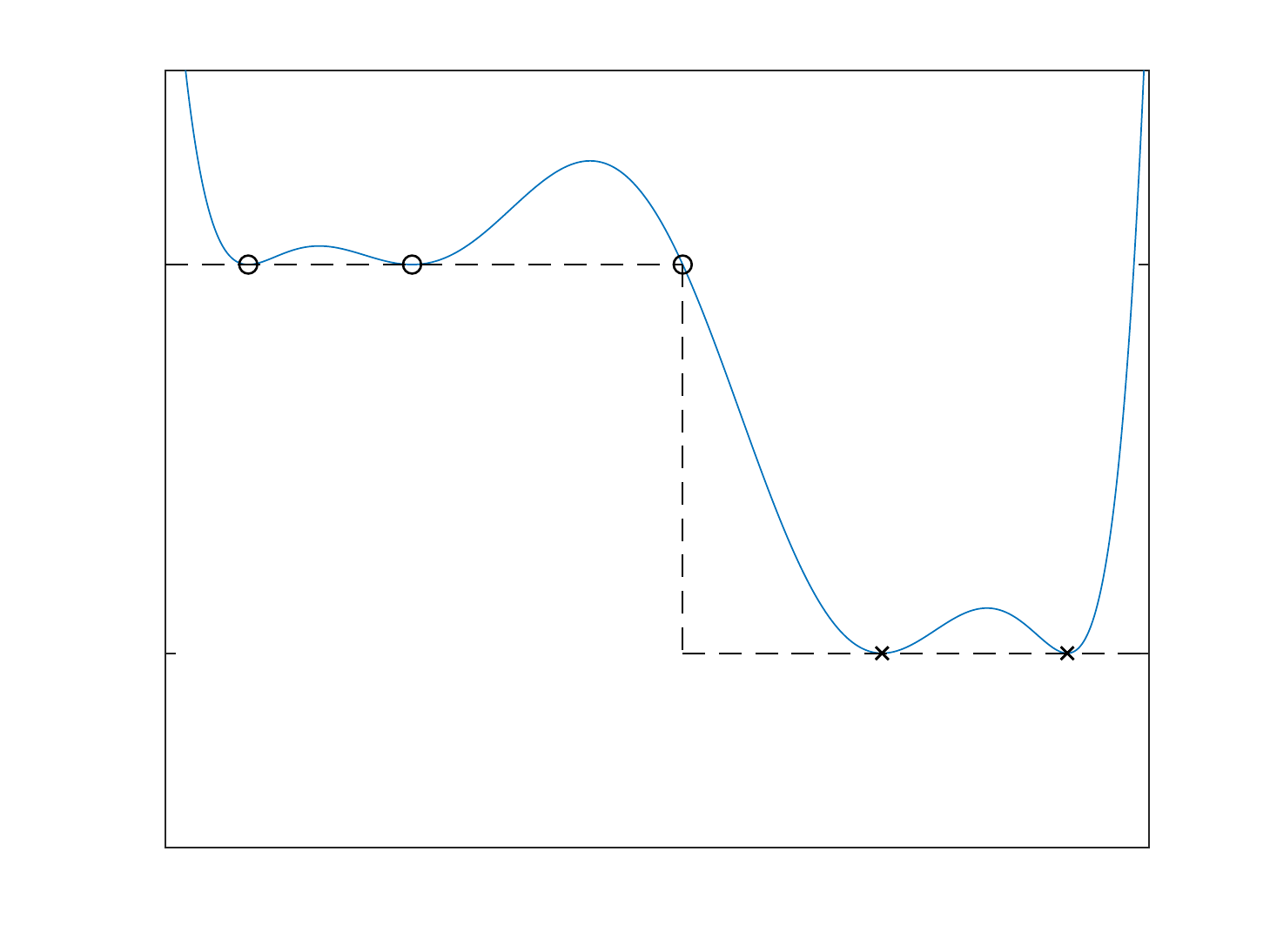}
\put(50,3){\small~$\lambda$}
\put(2,53){\small~$1.0$}
\put(2,22){\small~$0.0$}
\put(0,30){\rotatebox{90}{\small~$p(\lambda)$}}
\put(52,56){\small~$p(\theta_k)$}
\end{overpic}
\captionof{figure}{A numerical illustration of the
polynomial $p$ subject to the
conditions~\eqref{eq.pplus-conditions};
duplicated from~\figref{fig:pplusminusleft}.}
\end{minipage}
\hfill\vline\hfill
\begin{minipage}{0.475\textwidth}
\begin{equation}\label{eq.pplus-conditions}
\begin{array}{ll}
p(\theta_1)=1,~~&p'(\theta_1)=0,\\
~~~~\vdots &~~~~\vdots\\
p(\theta_{k-1})=1,~~&p'(\theta_{k-1})=0,\\
p(\theta_{k})=1,&\\
p(\theta_{k+1})=0,~~&p'(\theta_{k+1})=0,\\
~~~~\vdots &~~~~\vdots\\
p(\theta_{m})=0,~~&p'(\theta_{m})=0.
\end{array}
\end{equation}
Thus, we have $m$ many conditions for the polynomial $p$
and $m-1$ many conditions for its derivative,
and such a polynomial uniquely exists.
By the conditions~\eqref{eq.pplus-conditions}
the polynomial $p$ satisfies
the identities~\eqref{eq.ppmequal}.
\end{minipage}

\bigskip
\noindent
To show that $p$ satisfies the inequalities~\eqref{eq.ppminequal}
considering $p_{\{+,k\}}$, we proceed to locate the zeros of $p'$
which correspond to points of extreme values of $p$:
The derivative $p'$ is a polynomial of degree~$2m-3$,
and thus, has~$2m-3$ zeros.
By the conditions~\eqref{eq.pplus-conditions},
we have $m-1$ many zeros of $p'$ located at nodes.
For each pair of neighboring nodes in $\{\theta_1,\ldots,\theta_k\}$
and $\{\theta_{k+1},\ldots,\theta_m\}$
the conditions~\eqref{eq.pplus-conditions}
and Rolle's Theorem imply the existence of a
zero of~$p'$ between the respective nodes.
Thus, the derivative $p'$ has $m-1$ many simple zeros located at nodes and
$m-2$ many simple zeros located between nodes.
With~$p(\theta_{k})>p(\theta_{k+1})$ and
with the respective changes of sign for~$p'$
we conclude that $p$ satisfies the
inequalities for $p_{\{+,k\}}$ in~\eqref{eq.ppminequal}.

Furthermore,~we have $p(\lambda)>1$ for
$\lambda\in(\theta_j,\theta_{j+1})$ with $j=1,\ldots,k-1$
and~$\lambda<\theta_1$, and we have $p(\lambda)>0$
for $\lambda\in(\theta_j,\theta_{j+1})$ with $j=k,\ldots,m$
and~$\lambda>\theta_m$.
Thus, the inequalities for $p_{\{+,k\}}$ in~\eqref{eq.ppminequal}
are strict for $\lambda \notin \{\theta_1,\ldots,\theta_m\}$.

\bigskip
In a similar manner we conclude results for~$p_{\{-,k\}}$.
\newline\noindent
\begin{minipage}{0.475\textwidth}
\centering
\begin{overpic}
[width=0.9\textwidth]{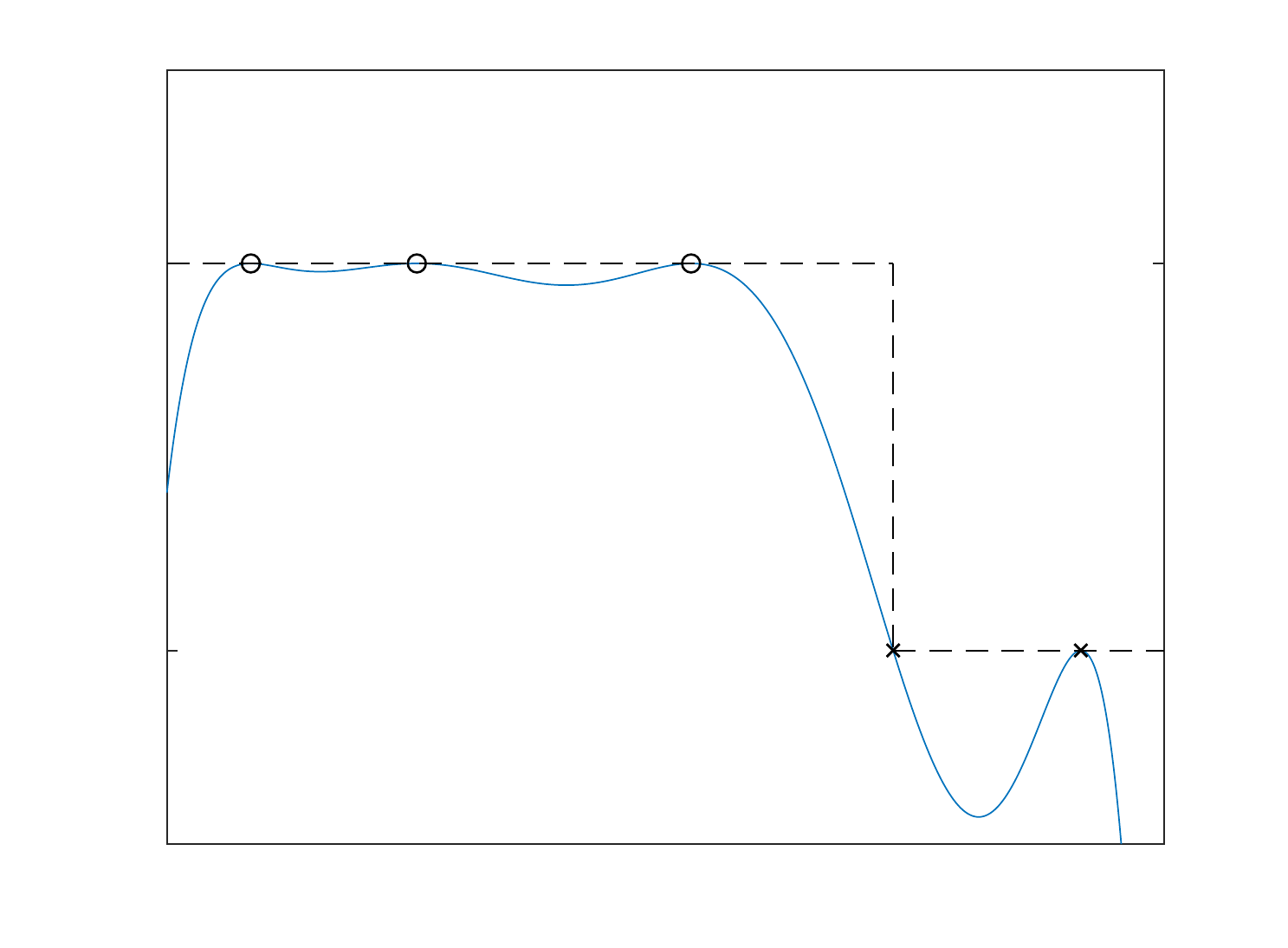}
\put(50,3){\small~$\lambda$}
\put(5,53){\tiny~$1.0$}
\put(5,22){\tiny~$0.0$}
\put(0,30){\rotatebox{90}{\small~$p(\lambda)$}}
\put(52,56){\small~$p(\theta_k)$}
\end{overpic}
\captionof{figure}{A numerical illustration of the
polynomial $p$ subject to the
conditions~\eqref{eq.pminus-conditions};
duplicated from~\figref{fig:pplusminusleft}.}
\end{minipage}
\hfill\vline\hfill
\begin{minipage}{0.475\textwidth}
Let $p$ be a polynomial of degree~$2m-2$ subject to the conditions
\begin{equation}\label{eq.pminus-conditions}
\begin{array}{ll}
p(\theta_1)=1,~~&p'(\theta_1)=0,\\
~~~~\vdots &~~~~\vdots\\
p(\theta_{k})=1,~~&p'(\theta_{k})=0,\\
p(\theta_{k+1})=0,&\\
p(\theta_{k+2})=0,~~&p'(\theta_{k+2})=0,\\
~~~~\vdots &~~~~\vdots\\
p(\theta_{m})=0,~~&p'(\theta_{m})=0.
\end{array}
\end{equation}
Then similar arguments as previously show that such a polynomial $p$
satisfies the identities~\eqref{eq.ppmequal}
and inequalities~\eqref{eq.ppminequal} associated with $p_{\{-,k\}}$.
\end{minipage}

\bigskip
Thus, the polynomials subject to the conditions~\eqref{eq.pplus-conditions}
and~\eqref{eq.pminus-conditions} satisfy the desired properties of
$p_{\{+,k\}}$ and $p_{\{-,k\}}$, respectively, which completes the proof.
\end{proof}

\begin{proof}[Proof of \,\textbf{Proposition~\ref{prop.thatr2}}]
Let $a=-\infty$ and $b=\infty$ to simplify the notation.

For the given pole $s\in\R$ we define the transformation
\begin{equation}\label{eq.pthatr1defx}
x\colon \R\setminus\{s\}\to\R\setminus\{0\},~~~~
x(\lambda):=(\lambda-s)^{-1}.
\end{equation}

For the case $\theta_1<s<\theta_m$ the indices $k_1>1$ and $k_m=k_1-1$
are given in~\eqref{eq.thatr2defk1ma} and
the values $x(\theta_j)$ 
satisfy the ordering
\begin{subequations}\label{eq.x1x2xm0ge}
\begin{equation}\label{eq.x1x2xm0}
x(\theta_{k_m}) < x(\theta_{k_m -1}) < \ldots < x(\theta_1) < 0 < x(\theta_m) < x(\theta_{m-1}) < \ldots < x(\theta_{k_1}).
\end{equation}
Otherwise, for $s<\theta_1$ (and $s>\theta_m$) we have
$k_1=1$ and $k_m=m$ as in~\eqref{eq.thatr2defk1mb}, and
\begin{equation}\label{eq.x1x2xm0x}
x(\theta_{k_m}) < \ldots < x(\theta_{k_1}) < 0, ~~~ s<\theta_1 ~~~~
\big(0 < x(\theta_{k_m}) < \ldots < x(\theta_{k_1}), ~~~ s>\theta_m \big).
\end{equation}
\end{subequations}
For the index $k_m$ we recall and highlight
\begin{equation}\label{eq.whatiskm}
k_m=k_1-1,~~~\text{for $\theta_1<s<\theta_m$,}~~~\text{and}~~~~
k_m=m,~~~\text{otherwise}.
\end{equation}

For any of these cases we define the index mapping
$$
\iota\colon \{1,\ldots,m\}\to \{1,\ldots,m\},~~~
\iota(j):=
\left\{
\begin{array}{ll}
k_1-j, &~~1\leq j <k_1,\\
m+k_1-j, &~~k_1\leq j \leq m.
\end{array}
\right.
$$
The action of $\iota$ is illustrated in the following table,
\begin{equation}\label{eq.prothatractioniota}
\begin{array}{|c|c|c|c|c|c|c|c|c|c|c|}
\hline
j &
1 & 2 &\cdots & k_1-2 & k_1-1 &
k_1 & k_1+1 & \cdots & m-1 & m\\\hline
\iota(j) &
k_1-1 & k_1-2& \cdots & 2 & 1 &
m & m-1 & \cdots & k_1+1 & k_1\\\hline
\end{array}\,,
\end{equation}
where $k_m=k_1-1$ or $k_m=m$ as specified in~\eqref{eq.whatiskm}.
(\,Thus, in the case of~$s<\theta_1$ or~$s>\theta_m$ this gives
\begin{equation}\label{eq.prothatractioniota1tom}
\begin{array}{|c|c|c|c|c|c|}
\hline
j &
1=k_1 & 2 & \cdots & m-1 & m=k_m\\\hline
\iota(j) &
m & m-1 & \cdots & 2 & 1\\\hline
\end{array}\,.~\big)
\end{equation}
We remark that $\iota$ is involutory with
$\iota(1)=k_m$ ($\iota(k_m)=1$),
and $\iota(m)=k_1$ ($\iota(k_1)=m$).
Thus, this mapping is bijective and with $\iota(k_m)=1$ we have
\begin{equation}\label{eq.prothatriiotanokm}
\iota(k)-1\in\{1,\ldots,m-1\},~~~\text{for}~~
k\in\{1,\ldots,m\}\setminus\{k_m\}.
\end{equation}
Let $\xi_1,\ldots,\xi_m$ denote the sequence of $x(\theta_j)$
arranged as in~\eqref{eq.x1x2xm0ge},~i.e.,
$$
\xi_j := x(\theta_{\iota(j)}), ~~~~\text{thus,} ~~~
\xi_1=x(\theta_{k_m}) <\ldots<\xi_m=x(\theta_{k_1}).
$$
We remark that $\iota$ being involutory implies
\begin{equation}\label{eq.prothatrdefxi2}
\xi_{\iota(j)}=x(\theta_j).
\end{equation}

We recall the definition of the index set $I_k$ given in~\eqref{eq.thatr2defRk},~i.e.,
\begin{equation}\label{eq.thatrIkrecall}
I_k=\left\{
\begin{array}{ll}
\{1,\ldots,k,k_1,\ldots,m \}, & ~~ 1 \leq k < k_1,\\
\{k_1,\ldots,k \}, & ~~ k_1\leq k \leq m.
\end{array}\right.
\end{equation}
The set $\{ x(\theta_j) \colon j\in I_k \}$ can be rewritten as follows:
With $\xi_{\iota(j)}=x(\theta_j)$~\eqref{eq.prothatrdefxi2} we have
\begin{equation}\label{eq.xofsetthetaktoxiiota1}
\{ x(\theta_j) \colon j\in I_k \}
=\left\{
\begin{array}{ll}
\{\xi_{\iota(1)},\ldots, \xi_{\iota(k)}\}
\cup \{\xi_{\iota(k_1)},\ldots, \xi_{\iota(m)}\},~~~ &1\leq k<k_1,~~~\text{and}\\
\{\xi_{\iota(k_1)},\ldots, \xi_{\iota(k)}\},~~~ &k_1\leq k \leq m.
\end{array}\right.
\end{equation}
We proceed to rewrite the indices of the sets
on the right-hand side of this equation using~\eqref{eq.prothatractioniota}.
In particularly, the identities $\iota(1)=k_1-1$, $\iota(k_1)=m$ and $\iota(m)=k_1$ 
imply
\begin{subequations}\label{eq.xofsetthetaktoxiiota2}
$$
\begin{aligned}
(\iota(1),\ldots, \iota(k))
&= (k_1-1, k_1-2,\ldots, \iota(k))~~~\text{for $k<k_1$ and}\\
(\iota(k_1),\ldots, \iota(m))
&= (m,m-1,\ldots,k_1).
\end{aligned}
$$
Thus,
\begin{equation}
\{\xi_{\iota(1)},\ldots, \xi_{\iota(k)}\}
\cup \{\xi_{\iota(k_1)},\ldots, \xi_{\iota(m)}\}
= \{ \xi_{\iota(k)},\xi_{m-1},\ldots, \xi_m \},~~~1\leq k<k_1.
\end{equation}
In a similar manner, identity~\eqref{eq.prothatractioniota} yields
$$
(\iota(k_1),\ldots, \iota(k))
= (m,m-1,\ldots,\iota(k)),~~~\text{for}~~k\geq k_1,
$$
which implies
\begin{equation}
\{\xi_{\iota(k_1)},\ldots, \xi_{\iota(k)}\}
= \{\xi_{\iota(k)},\ldots, \xi_m\},~~~k_1\leq k \leq m.
\end{equation}
\end{subequations}
The identity~\eqref{eq.xofsetthetaktoxiiota1}
together with~\eqref{eq.xofsetthetaktoxiiota2},
under consideration of the different cases for~$k$,
show
\begin{equation}\label{eq.r2Ikandxj}
\{ x(\theta_j)\colon j\in I_k \}
= \{\xi_{\iota(k)},\ldots,\xi_m \},
~~~ k=1,\ldots,m.
\end{equation}

In a similar manner the set~$R_k\subset\R$ in~\eqref{eq.thatr2defRk}
for $k=1,\ldots,m$ satisfies
\begin{subequations}\label{eq.r2Rkandxjboth}
\begin{equation}\label{eq.r2Rkandxj}
x(R_k) = [\xi_{\iota(k)},+\infty)\setminus \{0\},
~~~\text{and}~~
x(\R_s\setminus R_k) = (-\infty,\xi_{\iota(k)})\setminus \{0\},
\end{equation}
where $\R_s=\R\setminus\{s\}$.
The first identity in~\eqref{eq.r2Rkandxj}
is illustrated in \figref{fig:thatrproRk}.
Analogously, the interior of $R_k$ satisfies
\begin{equation}\label{eq.r2Rkandxjopen}
x(R_k^{\setopen}) = (\xi_{\iota(k)},+\infty)\setminus \{0\},
~~~\text{and}~~
x(\R_s\setminus R_k^{\setopen}) = (-\infty,\xi_{\iota(k)}]\setminus \{0\}.
\end{equation}
\end{subequations}

In the current setting we assume $k$
satisfies~$k\in\{1,\ldots,m\}\setminus\{k_m\}$,
thus, with~\eqref{eq.prothatriiotanokm} we have $\iota(k)-1\in\{1,\ldots,m-1\}$.
For the sequence $\xi_1<\ldots<\xi_m$ and
the index $\iota(k)-1$
we let~$p_{\{+,\iota(k)-1\}}$ and~$p_{\{-,\iota(k)-1\}}$ refer
to the polynomials introduced in Proposition~\ref{prop.thatpol1}.
Additionally, we define~$g_{\{\pm,k\}}$ by
\begin{equation}\label{eq.r2defgpmk}
g_{\{\pm,k\}}(y) :=
1-p_{\{\mp,\,\iota(k)-1\}}(y),~~~~k\in\{1,\ldots,m\}\setminus \{k_m\}.
\end{equation}
The identities~\eqref{eq.ppmequal} for $p_{\{\pm,\iota(k)-1\}}$
write
$$
p_{\{\pm,\iota(k)-1\}}(\xi_j)=\left\{
\begin{array}{ll}
1,&j=1,\ldots,\iota(k)-1,\\
0,&j=\iota(k),\ldots,m,
\end{array}
\right.
$$
and this entails the following identities for $g_{\{\pm,k\}}$,
$$
g_{\{\pm,k\}}(\xi_j)=\left\{
\begin{array}{ll}
0,&j=1,\ldots,\iota(k)-1,\\
1,&j=\iota(k),\ldots,m.
\end{array}
\right.
$$
With~\eqref{eq.r2Ikandxj} this conforms to
the following identities for the nodes $\theta_j$,
\begin{equation}\label{eq.ppmequalxx}
g_{\{\pm,k\}}(x(\theta_j))=\left\{
\begin{array}{ll}
1,&j \in I_k,\\
0,&\text{otherwise}.
\end{array}
\right.
\end{equation}

In a similar manner the inequalities~\eqref{eq.ppminequal}
for $p_{\{\pm,\iota(k)-1\}}$ read
\begin{equation}\label{eq.thatrproppminqubase}
p_{\{+,\iota(k)-1\}}(y) \geq \left\{
\begin{array}{ll}
1,& y \leq \xi_{\iota(k)-1},\\
0,& \xi_{\iota(k)-1} < y,
\end{array}
\right.
~~~\text{and}~~~
p_{\{-,\iota(k)-1\}}(y) \leq \left\{
\begin{array}{ll}
1,& y < \xi_{\iota(k)},\\
0,& \xi_{\iota(k)} \leq y,
\end{array}
\right.
\end{equation}
and this entails
\begin{equation}\label{eq.ppminequalxx}
g_{\{+,k\}}(y) \geq \left\{
\begin{array}{ll}
0,& y < \xi_{\iota(k)},\\
1,&  \xi_{\iota(k)} \leq y,
\end{array}
\right.
~~~\text{and}~~~
g_{\{-,k\}}(y) \leq \left\{
\begin{array}{ll}
0,& y \leq \xi_{\iota(k)-1},\\
1,& \xi_{\iota(k)-1} < y.
\end{array}
\right.
\end{equation}
With~\eqref{eq.r2Rkandxj} the
inequalities~\eqref{eq.ppminequalxx} for $g_{\{+,k\}}$
yield inequalities on the domain of $x$,
\begin{equation}\label{eq.ppminequalxx2a}
g_{\{+,k\}}( x(\lambda) ) \geq \left\{
\begin{array}{ll}
0,& \lambda \in \R_s \setminus R_k,\\
1,& \lambda \in R_k,
\end{array}
\right.
\end{equation}
To rewrite the inequalities~\eqref{eq.ppminequalxx} for $g_{\{-,k\}}$
we proceed in a similar manner:
We first consider the cases $s<\theta_1$ and $s>\theta_m$.
For these cases the action of the mapping~$\iota$
is illustrated in~\eqref{eq.prothatractioniota1tom} and we observe
$$
\iota(k)-1 = \iota(k+1),~~~k=1,\ldots,m-1,~~\text{and $s<\theta_1$ or $s>\theta_m$}.
$$
Thus, we have $\xi_{\iota(k)-1}=\xi_{\iota(k+1)}$ for these cases
and the identities~\eqref{eq.r2Rkandxjopen} imply
$$
x(R_{k+1}^{\setopen}) = (\xi_{\iota(k)-1},+\infty)\setminus \{0\},
~~~\text{and}~~
x(\R_s\setminus R_{k+1}^{\setopen}) = (-\infty,\xi_{\iota(k)-1}]\setminus \{0\}.
$$
Together with the inequalities for $g_{\{-,k\}}$ in~\eqref{eq.ppminequalxx},
this shows the following inequalities in the domain of $x$,
\begin{equation}\label{eq.prothatrineq21tm}
g_{\{-,k\}}(x(\lambda)) \leq \left\{
\begin{array}{ll}
0,& \lambda\in \R_s\setminus R_{k+1}^{\setopen},\\
1,& \lambda \in R_{k+1}^{\setopen},
\end{array}\right.
~~~~~~
\begin{array}{l}
\text{for $s<\lambda_1$ or $s>\lambda_m$, and}\\
~~k=1,\ldots,m-1.
\end{array}
\end{equation}
Similar results hold for the case $\theta_1<s<\theta_m$ (thus, $k_m<m$):
The illustration in~\eqref{eq.prothatractioniota} reveals
$$
\iota(k)-1 = \iota(k+1),~~~k\in\{1,\ldots,m-1\}\setminus\{ k_m\},
~~~~\text{and}~~~
\iota(m)-1 = \iota(1).
$$
Thus, with~\eqref{eq.r2Rkandxjopen} and the denotation $R_{m+1}=R_1$ we have
$$
x(R_{k+1}^{\setopen}) = (\xi_{\iota(k)-1},+\infty)\setminus \{0\},
~~~k\in\{1,\ldots,m\}\setminus\{ k_m\},
$$
with similar results considering $x(\R_s\setminus R_{k+1}^{\setopen})$.
With this identity, the inequalities for $g_{\{-,k\}}$
in~\eqref{eq.ppminequalxx} reveal inequalities
similar to~\eqref{eq.prothatrineq21tm} for the case $\theta_1<s<\theta_m$.
Together with~\eqref{eq.prothatrineq21tm}
for the case $s<\lambda_1$ or $s>\lambda_m$,
we conclude with
\begin{equation}\label{eq.prothatrineqminus}
g_{\{-,k\}}(x(\lambda)) \leq \left\{
\begin{array}{ll}
0,& \lambda\in \R_s\setminus R_{k+1}^{\setopen},\\
1,& \lambda \in R_{k+1}^{\setopen},
\end{array}
\right.
~~~~k\in\{ 1,\ldots,m\}\setminus\{k_m\},
\end{equation}

We define the rational function $r_{\{\pm,k\}}\in\Pi_{2m-2}/q_{m-1}^2$ by
$$
r_{\{\pm,k\}}(\lambda) := g_{\{\pm,k\}} (x(\lambda)) .
$$
Indeed, as demonstrated in Remark~\ref{rmk.gcircxisrat} further below,
such a function is rational.
In \figref{fig:thatrpro1} and~\ref{fig:thatrpro2}
we plot the rational function~$r_{\{\pm,m\}}$
and the respective auxiliary polynomial function $g_{\{\pm,m\}}$
for numerical examples.
For further illustrations of~$r_{\{-,m\}}$ we refer to \figref{fig:lemrpk}
in Subsection~\ref{subsec.SepTheoremRatKrylov}.

The rational functions $r_{\{\pm,k\}}$ satisfy the
identities~\eqref{eq.thatr2eq} and inequalities~\eqref{eq.thatr2ineq}
which concludes the proof of Proposition~\ref{prop.thatr2}:
The identities~\eqref{eq.ppmequalxx} conclude the
identities~\eqref{eq.thatr2eq} for
$ r_{\{\pm,k\}}(\theta_j) = g_{\{\pm,k\}} (x(\theta_j)) $.
Analogously,~\eqref{eq.ppminequalxx2a} and~\eqref{eq.prothatrineqminus}
entail the inequalities~\eqref{eq.thatr2ineq}.

\smallskip
Furthermore, we consider the inequalities~\eqref{eq.ppminequalxx2a}
and~\eqref{eq.prothatrineqminus}
to be strict for~$\lambda\neq\{\theta_1,\ldots,\theta_m\}$.
Indeed, for a given $\lambda$ with~$\lambda\neq\{\theta_1,\ldots,\theta_m\}$
we have $y=x(\lambda)\neq \{\xi_1,\ldots,\xi_m\}$ and the underlying
inequalities for $p_{\pm,\iota(k)-1}$ in~\eqref{eq.thatrproppminqubase}
are strict, which carries over to the inequalities~\eqref{eq.ppminequalxx2a}
and~\eqref{eq.prothatrineqminus}.
\end{proof}

\begin{figure}
\centering
\begin{overpic}
[width=0.5\textwidth]{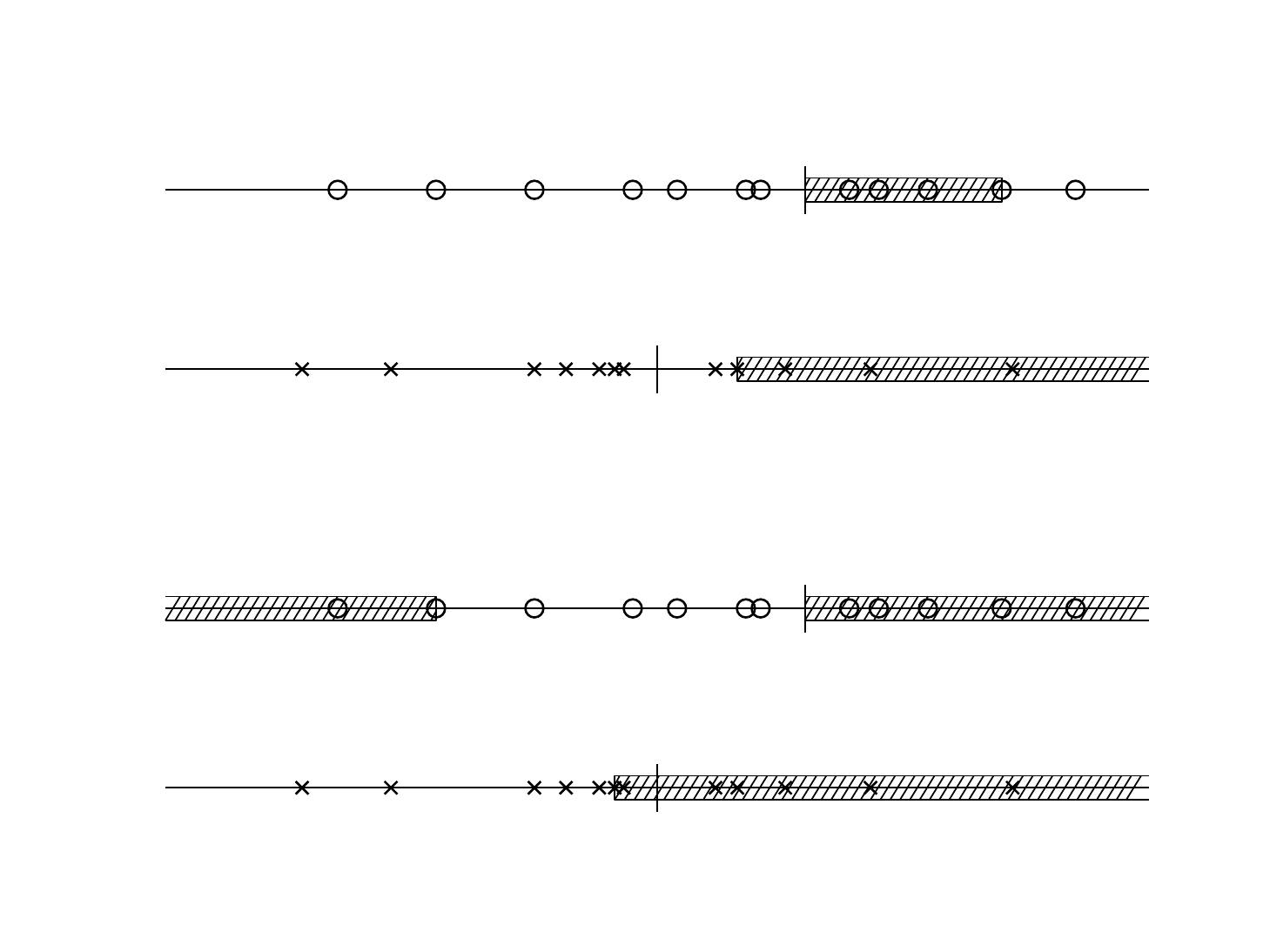}
\put(-40,59){\small \textbf{a)}}
\put(62.5,63){\small$ s $}   
\put(78,63){\small$\theta_k$}
\put(-30,59){\small {the case $\theta_k>s$:}}
\put(100,59){\small$R_k=(s,\theta_k]$}
\put(-40,45){\small \textbf{b)}}
\put(50.5,50){\small$ 0 $}
\put(57,50){\small$\xi_{\iota(k)}$}
\put(97,45){\small$x(R_k)=[\xi_{\iota(k)},+\infty)$}
\put(-40,25){\small \textbf{c)}}
\put(-30,25){\small {the case $\theta_k<s$:}}
\put(98,25){\small$R_k=(-\infty,\theta_k]\cup (s,+\infty)$}
\put(62.5,30){\small$ s $}   
\put(33,30){\small$\theta_k$}
\put(-40,11){\small \textbf{d)}}
\put(52,17){\small$ 0 $}
\put(42,17){\small$\xi_{\iota(k)}$}
\put(97,11){\small$x(R_k)=[\xi_{\iota(k)},+\infty)\setminus\{0\} $}
\put(80,8){\vector(1,0){10}}
\put(82,4){\small$ +\infty $}
\put(23,8){\vector(-1,0){10}}
\put(12,4){\small$ -\infty $}
\end{overpic}
\caption{In these figures we illustrate the
identity $x(R_k)=[\xi_{\iota(k)},+\infty)\setminus\{0\} $~\eqref{eq.r2Rkandxj}
for given nodes $\theta_1,\ldots,\theta_m$.
The pole $s$ is given and satisfies~$\theta_1<s<\theta_m$.
For the index $k$ we consider two different cases, namely,
we choose $k$ such that~$\theta_k>s$ in Figure~a) and~b),
and we choose $k$ such that~$\theta_k<s$ in Figure~c) and~d).
\newline\noindent--
Figure~\textbf{a)} and~\textbf{c)} show the real axis with the nodes
$\theta_1,\ldots,\theta_m$ ('$\circ$').
Furthermore, the set $R_k\subset\R$ given in~\eqref{eq.thatr2defRk}
is highlighted by a dashed area.
\newline\noindent--
Figure~\textbf{a)} and~\textbf{c)} show the real axis with
$\xi_1,\ldots,\xi_m$ ('$\times$'),~i.e.,~the
image of $\theta_1,\ldots,\theta_m$ under
the transformation $x$~\eqref{eq.pthatr1defx}
with $x(\theta_k)=\xi_{\iota(k)}$~\eqref{eq.prothatrdefxi2}.
Furthermore, the dashed area highlights $x(R_k)$
which satisfies the identity~\eqref{eq.r2Rkandxj}.
}
\label{fig:thatrproRk}
\end{figure}

\begin{remark}\label{rmk.gcircxisrat}
Let $g\in\Pi_{2m-2}$ and let $x(\lambda)=(\lambda-s)^{-1}$, then
\begin{equation}\label{eq.4apgpoltorat0}
r(\lambda)=g(x(\lambda))
\end{equation}
defines a rational function in $\lambda$, namely, $r\in\Pi_{2m-2}/q_{m-1}^2$
for $q_{m-1}(\lambda)=(\lambda-s)^{m-1}$.
To demonstrate this result we define
\begin{equation}\label{eq.4apgpoltorat}
\widehat{g}(\lambda) = g\big((\lambda-s)^{-1}\big) (\lambda-s)^{2m-2}.
\end{equation}
Expanding the right-hand side of~\eqref{eq.4apgpoltorat}
shows $\widehat{g}\in\Pi_{2m-2}$.
Substituting $x(\lambda)$ and $q_{m-1}(\lambda)$
in~\eqref{eq.4apgpoltorat}, and dividing by $q_{m-1}(\lambda)^2$
reveals the representation
$$
\widehat{g}(\lambda)/q_{m-1}(\lambda)^2 = g(x(\lambda)),
$$
and thus, with~\eqref{eq.4apgpoltorat0} we have
$ r(\lambda) = \widehat{g}(\lambda)/q_{m-1}(\lambda)^2 $.
This shows $r\in\Pi_{2m-2}/q_{m-1}^2$.
\end{remark}

\begin{figure}
\centering
\begin{tabular}{cc}
\begin{overpic}
[width=0.5\textwidth]{06-rpklem2}
\put(0,27){\rotatebox{90}{\small$r_{\{-,k\}}(\lambda)$}}
\put(-4,65){ \textbf{a)}}
\put(54,63){\small$s$}
\put(15,63){\small$k=8$}
\put(72,43){\scriptsize$r_{\{-,k\}}(\theta_k)$}
\put(83,48){\vector(0,1){5}}
\put(50,-1){\small$\lambda$}
\end{overpic}
&
\begin{overpic}
[width=0.5\textwidth]{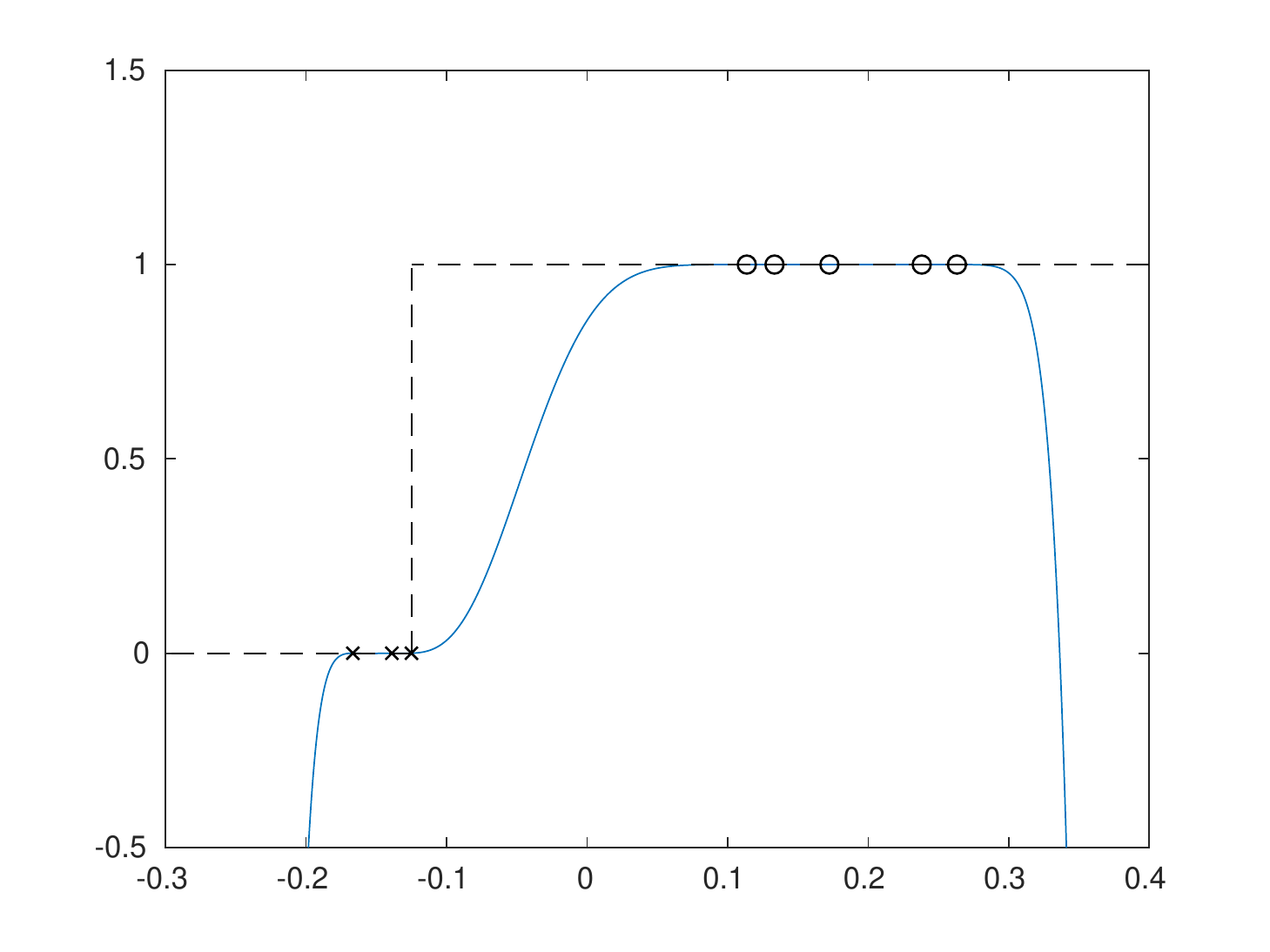}
\put(0,27){\rotatebox{90}{\small$g_{\{-,k\}}(x)$}}
\put(-4,65){ \textbf{b)}}
\put(15,63){\small$k=8, \,\iota(k)=4$}
\put(45,-1){\small$x=x(\lambda)$}
\put(47,43){\scriptsize$g_{\{-,k\}}(\xi_{\iota(k)})$}
\put(59,48){\vector(0,1){5}}
\end{overpic}
\end{tabular}
\caption{-- In Figure~\textbf{a)}
we plot the rational function~$r_{\{-,k\}}$
(introduced in Proposition~\ref{prop.thatr2})
for a numerical example;
we show~$r_{\{-,k\}}(\lambda)$
over~$\lambda$ for a given pole~$s=-3$,
and given nodes $\theta_1,\ldots,\theta_m$ with $m=8$.
We have $\theta_1<s<\theta_m$, namely, $\theta_{k_m}<s<\theta_{k_1}$
with $k_m=3$ and $k_1=4$.
For $j\in I_k$ we mark $r_{\{-,k\}}(\theta_j)$ by ('$\circ$'),
and for $j\notin I_k$ we mark $r_{\{-,k\}}(\theta_j)$ by ('$\times$').
The dashed lines illustrate the upper bounds of~$r_{\{-,k\}}$
given in~\eqref{eq.thatr2ineq}.
\newline\noindent--
In Figure~\textbf{b)} we show the auxiliary polynomial function~$g_{\{-,k\}}$
which appears in the proof of Proposition~\ref{prop.thatr2},
namely,~\eqref{eq.r2defgpmk} therein.
The nodes~$\xi_1,\ldots,\xi_m$
correspond to the image of the nodes~$\theta_j$ under~$x$,
namely,~$\xi_{\iota(j)}=x(\theta_j)$ for~$j=1,\ldots,m$.
For $k=8$ we have $\iota(k)=4$.
The symbols~('$\times$') and~('$\circ$') mark
$g_{\{-,m\}}(\xi_j)$ for~$j=\iota(k)-1,\ldots,m$
and~$j=\iota(k),\ldots,m$, respectively.
The dashed lines illustrate the upper bounds of~$g_{\{-,k\}}$
given in~\eqref{eq.ppminequalxx}.}
\label{fig:thatrpro1}
\end{figure}

\begin{figure}
\centering
\begin{tabular}{cc}
\begin{overpic}
[width=0.5\textwidth]{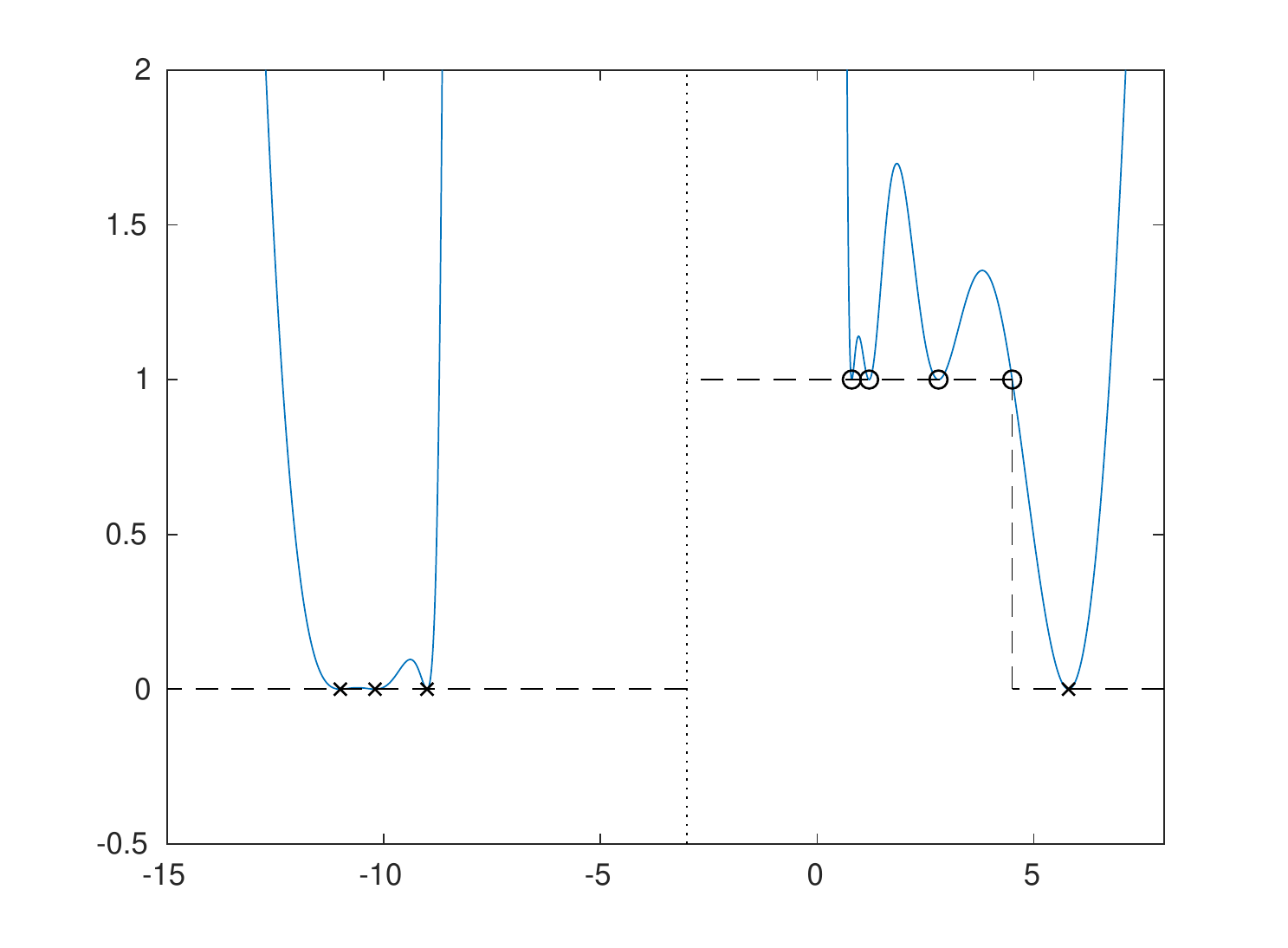}
\put(0,27){\rotatebox{90}{\small$r_{\{+,k\}}(\lambda)$}}
\put(-4,65){ \textbf{a)}}
\put(54,63){\small$s$}
\put(15,63){\small$k=7$}
\put(55,34){\scriptsize$r_{\{+,k\}}(\theta_k)$}
\put(72,38){\vector(1,1){5}}
\put(50,-1){\small$\lambda$}
\end{overpic}
&
\begin{overpic}
[width=0.5\textwidth]{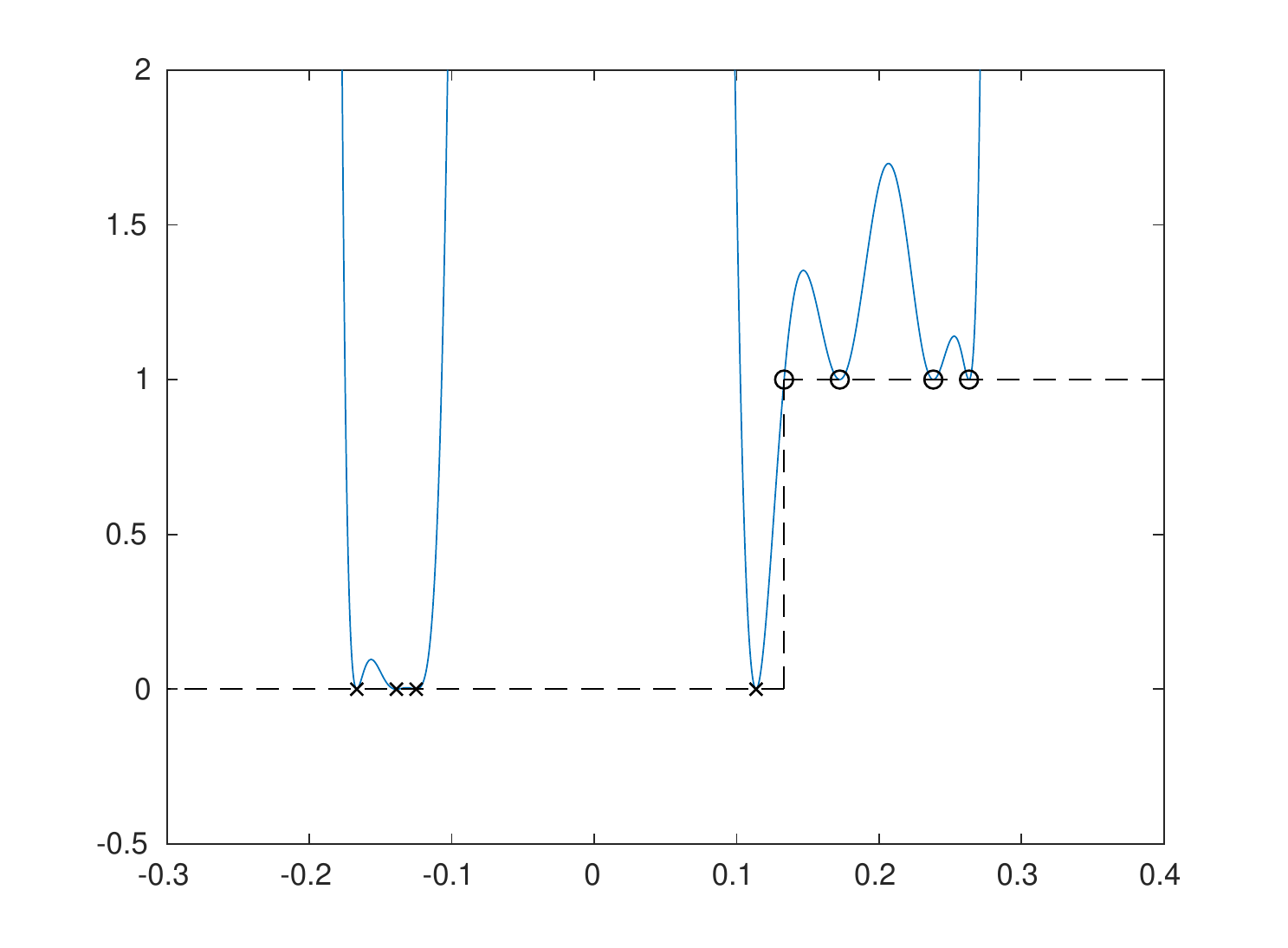}
\put(0,27){\rotatebox{90}{\small$g_{\{+,k\}}(x)$}}
\put(-4,65){ \textbf{b)}}
\put(54,63){\small$s$}
\put(15.3,63){\small$k=7,$}
\put(15,58){\small$\iota(k)=5$}
\put(45,-1){\small$x=x(\lambda)$}
\put(66,35){\scriptsize$g_{\{+,k\}}(\xi_{\iota(k)})$}
\put(67,38){\vector(-1,1){5}}
\end{overpic}
\\
\begin{overpic}
[width=0.5\textwidth]{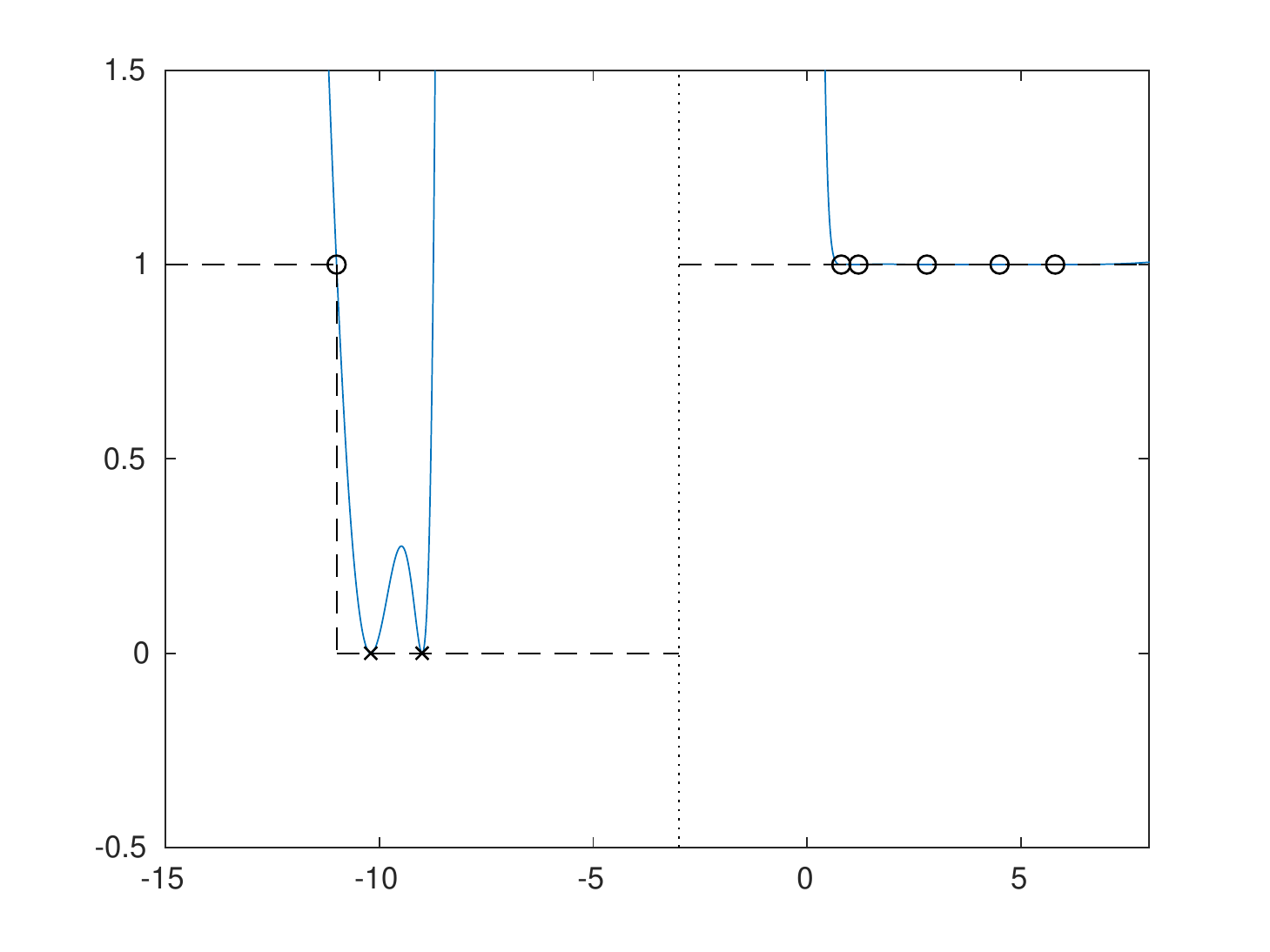}
\put(0,27){\rotatebox{90}{\small$r_{\{+,k\}}(\lambda)$}}
\put(-4,65){ \textbf{c)}}
\put(54,63){\small$s$}
\put(14,64){\small$k=1$}
\put(33,63){\scriptsize$r_{\{+,k\}}(\theta_k)$}
\put(33,61){\vector(-1,-1){5}}
\put(50,-1){\small$\lambda$}
\end{overpic}
&
\begin{overpic}
[width=0.5\textwidth]{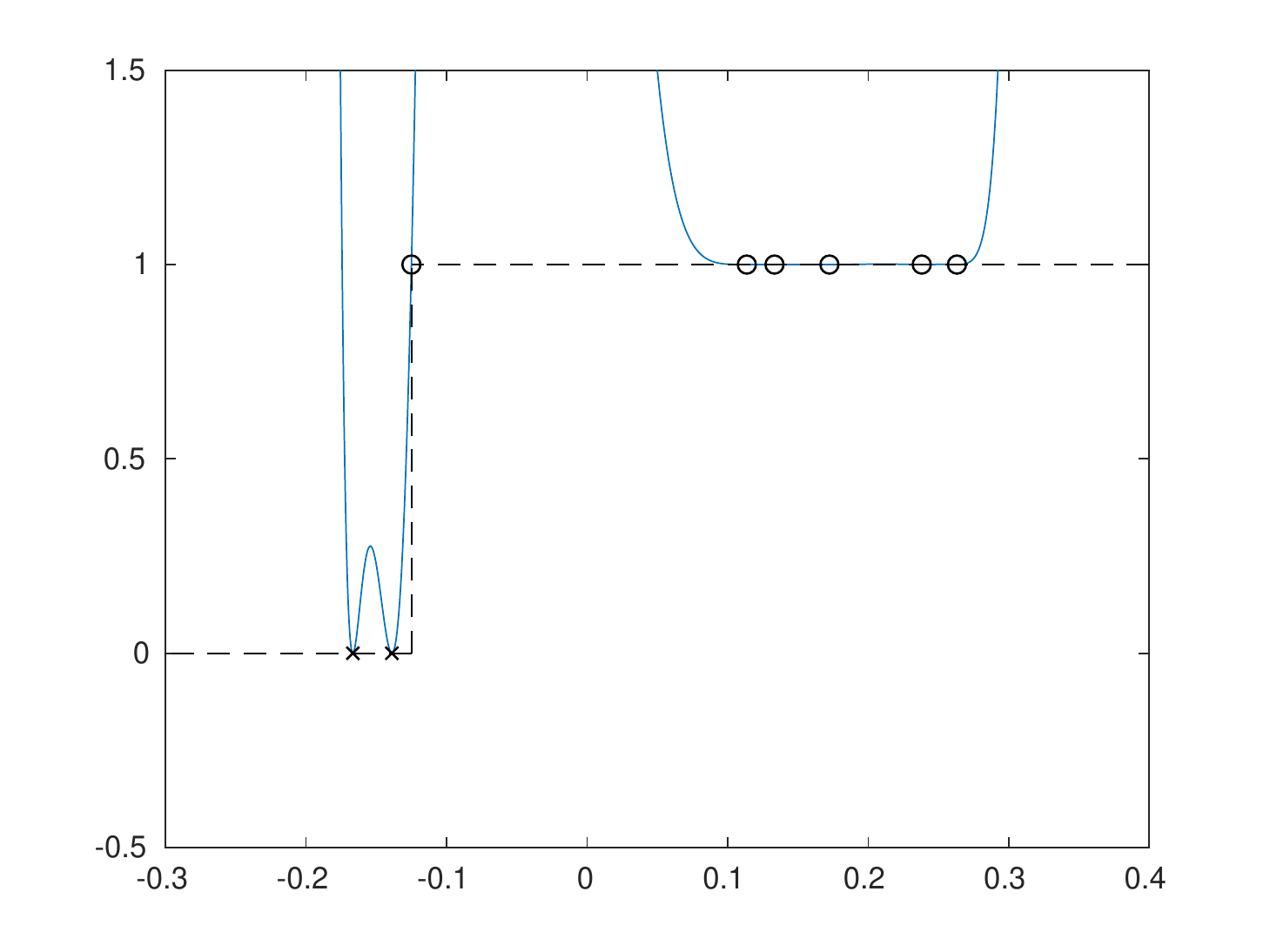}
\put(0,27){\rotatebox{90}{\small$g_{\{+,k\}}(x)$}}
\put(-4,65){ \textbf{d)}}
\put(14,64){\small$k=1,$}
\put(14,59){\small$\iota(k)=3$}
\put(45,-1){\small$x=x(\lambda)$}
\put(36,44){\scriptsize$g_{\{+,k\}}(\xi_{\iota(k)})$}
\put(39,48){\vector(-1,1){5}}
\end{overpic}
\end{tabular}
\caption{--\,In Figure~\textbf{a)} and~\textbf{c)}
we plot the rational function~$r_{\{+,k\}}$
(introduced in Proposition~\ref{prop.thatr2})
for numerical examples.
The nodes $\theta_1,\ldots,\theta_m$ with $m=8$
and the pole $s$ are given as in \figref{fig:thatrpro1} a).
For the index $k$ we choose $k=7$ in Figure~a)
and $k=1$ in Figure~c).
For $j\in I_k$ we mark $r_{\{+,k\}}(\theta_j)$ by~('$\circ$'),
and for $j\notin I_k$ we mark $r_{\{+,k\}}(\theta_j)$ by~('$\times$').
The dashed lines illustrate the upper bounds
of~$r_{\{+,k\}}$ given in~\eqref{eq.thatr2ineq}.
\newline\noindent --
In Figure~\textbf{b)} and~\textbf{c)}
we show the auxiliary polynomial function~$g_{\{+,k\}}$
which appears in the proof of Proposition~\ref{prop.thatr2},
namely,~\eqref{eq.r2defgpmk} therein.
The function~$g_{\{+,k\}}$ with $k=7$ plotted in Figure~b)
is associated with the function~$r_{\{+,k\}}$ in Figure~a),
and analogously, such a relation is given for Figure~d) and Figure~c).
For the index $\iota(k)$ which is relevant in the
proof of Proposition~\ref{prop.thatr2} we remark $\iota(k)=5$ for $k=7$
and $\iota(k)=3$ for $k=1$.
As in \figref{fig:thatrpro1} b),
the nodes~$\xi_1,\ldots,\xi_m$
correspond to the image of the nodes~$\theta_j$ under~$x$,
namely,~$\xi_{\iota(j)}=x(\theta_j)$ for~$j=1,\ldots,m$.
The symbols~('$\times$') and~('$\circ$') mark
$g_{\{+,m\}}(\xi_j)$ for~$j=\iota(k)-1,\ldots,m$
and~$j=\iota(k),\ldots,m$, respectively.
The dashed lines illustrate the upper bounds of~$g_{\{+,k\}}$
given in~\eqref{eq.ppminequalxx}.}
\label{fig:thatrpro2}
\end{figure}

\ifthesis
\end{subappendices}
\fi


\end{document}